\documentclass[aos]{imsart}

\RequirePackage{amsthm,amsmath,amsfonts,amssymb}
\RequirePackage[authoryear]{natbib}
\RequirePackage[colorlinks,citecolor=blue,urlcolor=blue]{hyperref}
\usepackage{dsfont,subcaption, mathtools, graphicx,color, bbm, mathrsfs}

\startlocaldefs

\newcommand{\dH}{\dot{H}}
\newcommand{\divH}{\dot{H}_{\diamond}}
\newcommand{\Hnull}{\divH}

\numberwithin{equation}{section}

\theoremstyle{plain}
\newtheorem{Theor}{Theorem}[section]
\newtheorem{Lem}[Theor]{Lemma}

\newtheorem{Prop}[Theor]{Proposition}
\newtheorem{Corol}[Theor]{Corollary}
\newtheorem{Condition}[Theor]{Condition}

\theoremstyle{definition}

\newtheorem{Remark}[Theor]{Remark}

\newcommand{\B}{\mathbb{B}}
\newcommand{\C}{\mathbb{C}}

\newcommand{\E}{\mathbb{E}}

\newcommand{\I}{\mathbb{I}}

\newcommand{\N}{\mathbb{N}}

\newcommand{\R}{\mathbb{R}}

\newcommand{\T}{\mathbb{T}}

\newcommand{\W}{\mathbb{W}}

\newcommand{\Z}{\mathbb{Z}}

\renewcommand{\AA}{\mathcal{A}}
\newcommand{\BB}{\mathcal{B}}

\newcommand{\HH}{\mathcal{H}}

\newcommand{\KK}{\mathcal{K}}
\newcommand{\LL}{\mathcal{L}}

\newcommand{\NN}{\mathcal{N}}

\newcommand{\XX}{\mathcal{X}}

\newcommand{\CCC}{\mathscr{C}}

\newcommand{\GGG}{\mathscr{G}}

\newcommand{\WWW}{\mathscr{W}}

\def\1{{\mathbb I}}

\newcommand{\ie}{\textit{i.e.} }

\newcommand{\ps}[2]{\left\langle #1,#2 \right\rangle}
\newcommand{\sm}{\setminus}

\newcommand{\ve}{\varepsilon}
\newcommand{\Om}{\Omega}

\newcommand{\spa}{\quad\quad}

\newcommand{\lad}{Ladyzhenskaya}

\newcommand{\embed}{\hookrightarrow}

\usepackage{scalerel,stackengine}
\stackMath
\newcommand\wwidehat[1]{%
	\savestack{\tmpbox}{\stretchto{%
			\scaleto{%
				\scalerel*[\widthof{\ensuremath{#1}}]{\kern-.6pt\bigwedge\kern-.6pt}%
				{\rule[-\textheight/2]{1ex}{\textheight}}
			}{\textheight}%
		}{0.5ex}}%
	\stackon[1pt]{#1}{\tmpbox}%
}

\newenvironment{Proof}[1]
{
	\begin{proof}[Proof of #1]
	}
	{
	\end{proof}
}

\allowdisplaybreaks

\endlocaldefs

\begin{document}

\begin{frontmatter}
\title{Data assimilation with the $2D$ Navier-Stokes equations: Optimal Gaussian asymptotics for the posterior measure}
\runtitle{Data assimilation with the $2D$ Navier-Stokes equations}

\begin{aug}
\author[A]{\fnms{Dimitri}~\snm{Konen}\ead[label=e1]{dk738@cam.ac.uk}}
\and
\author[A]{\fnms{Richard}~\snm{Nickl}\ead[label=e2]{nickl@maths.cam.ac.uk}}
\address[A]{Department of Pure Mathematics and Mathematical Statistics, University of Cambridge, Cambridge, CB3 0WA, United Kingdom\printead[presep={,\ }]{e1,e2}}

\end{aug}

\begin{abstract}
A functional Bernstein-von Mises theorem is proved for posterior measures arising in a data assimilation problem with the two-dimensional Navier-Stokes equation where a Gaussian process prior is assigned to the initial condition of the system. The posterior measure, which provides the update in the space of all trajectories arising from a discrete sample of the (deterministic) dynamics, is shown to be approximated by a Gaussian random vector field obtained as the solution to a linear parabolic PDE with Gaussian initial condition. The approximation holds in the strong sense of the supremum norm on the regression functions, showing that predicting future states of Navier-Stokes systems admits $1/\sqrt N$-consistent estimators even for commonly used nonparametric models. Consequences for coverage of credible bands and uncertainty quantification are discussed. A local asymptotic minimax theorem is derived that describes the lower bound for estimating the state of the nonlinear system, which is shown to be attained by the Bayesian data assimilation algorithm.
\end{abstract}



\end{frontmatter}


\setcounter{tocdepth}{2}
{
\hypersetup{linkcolor=black}
\tableofcontents
}

\section{Introduction}\label{sec:Intro}

\subsection{Bayesian data assimilation for a periodic Navier-Stokes model}

Many physical phenomena can be described by the nonlinear evolution of a dynamical system whose state at any point in time $t>0$ is described by an $m$-dimensional vector field $u(t, \cdot)$ over some Euclidean domain $\Omega$, and where the infinitesimal dynamics are governed by a partial differential equation (PDE). The state $u=u_\theta$ further depends on some generally unknown initial condition $\theta = u(0,\cdot)$. In `\textit{data assimilation}' (see, e.g., \cite{E09, LSZ15, RC15, EVvL22} and references therein), this initial condition is modeled by a Gaussian random field over $\Omega$, resulting in a random trajectory $u_\theta(t,\cdot)$ on the time-space cylinder. Given discrete measurements of the system $u_\theta$ corrupted by additive Gaussian noise $\varepsilon_i$, one then wishes to update the `prior trajectory' to the best `posterior' inference on the states of the dynamical system. This extends classical ideas in stochastic control and linear filtering due to  \cite{K60} to a nonlinear and infinite-dimensional setting, and can be regarded as a problem of Bayesian inference in a nonparametric regression context (see \cite{S10, GV17, NicklEMS}) where the regression function is $u_\theta$ and where the posterior arises from a Gaussian process prior $\Pi$ for $\theta$. In understanding the statistical properties of this methodology, it is natural to first study the so-called setting of `\textit{deterministic dynamics}' where the underlying PDE model is assumed to be correct (rather than subject to further, possibly stochastic, misspecification error).  In this case the posterior measure $\Pi(\cdot|Z^{(N)})$ arises from observations $Z^{(N)}=(Y_i, t_i, \omega_i)_{i=1}^N$ from the equation
\begin{equation}\label{model}
Y_i = u_\theta(t_i, \omega_i) + \varepsilon_i, ~~ i=1,\dots, N,
\end{equation}
where the $X_i=(t_i, \omega_i)$ are sampled uniformly at random from the time-space cylinder $[0,T] \times \Omega$. For simplicity, we assume that the error variables $\varepsilon_i \sim N(0,I_{\R^m})$ are also drawn iid, independently of the design -- this is not necessary but simplifies the theoretical development. Denote by $P_\theta^\N$ the infinite product probability law of $(Y_i, X_i)_{i=1}^\infty$, with corresponding expectation operator $E_{\theta}^\mathbb N$. While the posterior measure is initially supported in the space of initial conditions $\theta$, its push-forwards under the map $\theta \mapsto u_\theta(t)$ induce the updated posterior measures in the space of trajectories, see (\ref{imagemeasure}) below.

When the dynamics are non-linear, the posterior measures are not explicitly available and need to be computed by iterative methods such as MCMC, stochastic gradient descent, or filtering, including a numerical solver for the forward PDE. The hope that such methods produce reliable outputs is often related to an implicit hypothesis of a Gaussian approximation of the posterior, see e.g., \cite{CHSV24} for a recent reference. Even when computable (more on this below), it is highly unclear whether posterior inferences can be trusted to be statistically accurate. In this article, we develop some mathematical theory that justifies the use of such algorithms for the prototypical example of the incompressible $2D$-Navier-Stokes equations (see \cite{CF88, Robinson2001}) -- these are widely used in data assimilation tasks arising throughout geophysical sciences. Specifically, consider as spatial domain $\Om$ the two-dimensional torus $\T^2 =[0,1]^2$ with opposite points identified, as well as two-dimensional periodic vector fields $u_\theta = u:[0,T]\times \Om\to\R^2$ solving the PDE
 \begin{align}\label{eq:NSequations}\nonumber
    \frac{\partial u}{\partial t} - \nu \Delta u + (u\cdot \nabla) u  &= f - \nabla p
    ~~~ \textrm{on } (0,T]\times\Om, \\ 
    u(0)&=\theta
    ~~~ \textrm {on } \Om, \\ \nonumber
    \nabla\cdot u &= 0
    ~~~ \textrm {on } [0,T]\times \Om.
 \end{align}
    Here $\nu>0$ is a known `viscosity' term, $f:\Om\to\R^2$ is a given (time-independent) exterior forcing,  $p:[0,T]\times \Om\to \R$ is a scalar pressure term to be solved for alongside $u$, and we employ standard vector calculus notation (to be reviewed in detail later). As is common, we shall consider a suitably projected equation given in (\ref{eq:NS}) below, where $P$ is the Leray projector, $A=-P\Delta$ the Stokes operator, and $B[u,v]= P[(u\cdot \nabla) v]$ describes the nonlinearity of the system. The pressure term $\nabla p$ lies in the kernel of the Leray projector and can hence be disregarded when studying the projected equation, as we shall do here.

To introduce our results, we recall a basic Bayesian data assimilation model (`without model error') from \cite{CDRS09, S10, NicTit2024}. In the Navier-Stokes context, we will consider a Gaussian 
prior $\Pi={\rm Law}(\theta)$ for the initial condition
\begin{equation}\label{prior}
\theta \sim \mathcal N(0, \rho^2 A^{-\alpha}),~~ \alpha>3,\ \rho>0,
\end{equation} on the Hilbert space of vector fields
$$
\Hnull =\bigg\{ u\in L^2(\Om)^2 : \nabla\cdot u = 0,\ \int_{\Om} u = 0 \bigg\},
$$ 
where $$\rho=\rho_N=N^{-1/(2\alpha+2)}$$ is an appropriate variance scaling. In our periodic setting such priors can be constructed as Gaussian series expansions for the orthonormal basis functions of $\Hnull$ given by the vector fields 
$$
c_k(x)
\propto (k_2, -k_1) \cos (2\pi k \cdot x),\quad
s_k(x)
\propto (k_2, -k_1) \sin (2\pi k \cdot x),\quad 
x \in \Omega
$$ 
for $~k =(k_1, k_2) \in \Z^2 \setminus \{(0,0)\}$, arising as the eigenfunctions of $A$.  We can enumerate the square lattice by $k=k_j$, $j=1, 2,\dots$, in such a way that
\begin{equation}\label{eq:BasisDivH}
f_{2j-1} \equiv c_{k_j},\spa 
f_{2j} \equiv s_{k_j},\spa 
A f_j = \lambda_j f_j,\spa 
0<\lambda_j \simeq j,
\end{equation}
where the corresponding Weyl's law for the eigenvalues $\lambda_j$ follows as in Proposition 4.14 in \cite{CF88}. These eigenfunctions span the RKHS of $\theta$ and the Gaussian random series of vector fields
\begin{equation}\label{PriorSeries}
\theta(x_1, x_2) 
= \rho \sum_{j  \ge 1} g_j \lambda_j^{-\alpha/2} f_j(x_1, x_2),\quad (x_1,x_2) \in \Omega,\quad g_j \sim^{iid} N(0,1),
\end{equation}
converges almost surely in $\divH^\beta$ for any $1<\beta<\alpha-1$, where $\divH^\beta$ will be defined below as the closed subspace $\divH^\beta \equiv H^\beta(\Omega)^2 \cap \Hnull$ of the usual Sobolev space $H^\beta(\Omega)^2,  \beta \in \mathbb N$. The resulting Borel law of $\theta$ on $\divH^\beta \subset \Hnull \cap C(\Omega)^2$ induces the prior probability measure $\Pi$ alluded to above, using also a Sobolev imbedding into the space $C(\Omega)^2$ of continuous vector fields over~$\Omega$. For further details on standard results about Gaussian measures and processes, we refer the reader to Chapter 2, specifically Theorem 2.6.10, in \cite{GinNic2016}, and Section B.1.3 in \cite{NicklEMS}.

Given the prior $\Pi$ from (\ref{prior}) and data $Z^{(N)}$ from (\ref{model}) where $u_\theta$ solves the Navier-Stokes equation (\ref{eq:NS}), we obtain a posterior probability measure on $\Hnull$ as
 \begin{equation} \label{post}
d\Pi(\theta|Z^{(N)}) \propto \exp \Big\{-\frac{1}{2} \sum_{i=1}^N|Y_i-u_\theta(t_i, \omega_i)|^2 \Big\} d\Pi(\theta),~~\theta \in \divH,
\end{equation}
by applying standard arguments on conditional probabilities in dominated families, e.g., \cite{GV17} or Section 1.2.3 in \cite{NicklEMS}. The posterior expectation of $\|\theta\|_{\Hnull}$ is finite as is easy to show using the techniques we develop below, and we denote the posterior mean function in $\Hnull$ by $\tilde \theta_N = E^\Pi[\theta|Z^{(N)}]$. We also obtain an update for the whole dynamical system $u_\theta(t)$ whose marginal time distributions on $\Hnull$ are given by the image measures
\begin{equation} \label{imagemeasure}
\hat \Pi_{t,N} = {\rm Law} (u_\theta(t)), ~\theta \sim \Pi(\cdot|Z^{(N)}), ~t \ge 0.
\end{equation}
For the last observation times $t_{(N)} = \max_{i \le N} t_i \in (0,T],$ the preceding laws are sometimes called \textit{filtering distributions}. Since the prior draws take values in $\divH^\beta \subset C(\Omega)^2$ almost surely, so does any draw $\theta|Z^{(N)}$ from the posterior, and the same is then true also for the induced posterior law of $u_\theta(t), t>0,$ by Proposition \ref{PropNSBound}.

\subsection{Main contributions}

It was proved recently in \cite{NicTit2024} that under natural hypotheses on the ground truth initial condition $\theta_0$, the posterior measures $\hat \Pi_{t,N}, t \ge 0,$ and their mean vectors are statistically consistent for the recovery of the true state $u_{\theta_0}$ of the system in the large sample limit when $N \to \infty$. The convergence rates of $\tilde \theta_N$ in $L^2$-distance towards $\theta_0$ were shown to be generally of inverse logarithmic order in $N$, and it was further proved in Theorem 4 of \cite{NicTit2024} that no faster rate is possible when the $t_i$'s are sampled uniformly from $[T_0,T]$ for $T_0>0$, corresponding to a situation where one cannot measure at `the beginning of time'. In contrast it was recently shown in \cite{Nickl2024} that under general hypotheses on nonlinear dynamical systems (covering, for example, a closely related nonlinear reaction diffusion equation with `nice' reaction functions), one can obtain `parametric' statistical convergence rates and a Bernstein-von Mises (BvM) theorem for the measures $\hat \Pi_{t,N}$ in uniform norm topologies. This is true even when one maintains a `nonparametric' (infinite-dimensional) model for $\theta$ as long as the sampling process begins at $T_0=0$. In the present paper we develop techniques that allow to apply the general theory from \cite{Nickl2024} to the $2D$-Navier Stokes system. Our techniques further imply that such Bayesian data assimilation algorithms are  \textit{optimal} from an objective information-theoretic point of view in the sense that they attain the exact asymptotic local minimax constant for the state estimation problem. More generally our results inform other statistical inference problems in dissipative time evolution systems such as those studied recently in diffusion models in \cite{NR20, GR22, N24diff, GW25, NPR25, HR25, NS25}, parabolic Schr\"odinger equations \cite{K22, KSV24} or other data assimilation tasks  \cite{CHSV24, SAW24, RSGK24, DBBBBG25}.

\smallskip

The starting point of our analysis is initially to obtain `faster than logarithmic' posterior contraction rates for the initial conditions $\theta$, see Theorem \ref{TheorFastRateTheta} below, which requires novel arguments compared to those in \cite{Nickl2024} due to a lack of symmetry of the linearisation `score' operator of the flow in the Navier-Stokes model. This in turn allows to derive a functional Bernstein-von Mises theorem following ideas in \cite{CN14}, \cite{Nickl2020, Nickl2024}, if the underlying information equation is solvable, which we will show here to be the case for the Navier-Stokes system. Let us describe the limit process occurring in our main theorem, which is a Gaussian random field in the space of trajectories on the time-space cylinder. Let $u_{\theta_0}$ be the solution to (\ref{eq:NSequations}) with `sufficiently regular' initial condition $\theta_0$, and consider the Gaussian random vector field $U=(U_1, U_2)$ over $(0,T] \times \Omega$ obtained from the unique weak solution in $\Hnull$ of the following \textit{linear} time-dependent parabolic equation with random initial condition:
\begin{align} \label{linshow}
\frac{\partial}{\partial t} U(t,\cdot) + \nu A U(t, \cdot)  + B[u_{\theta_0}(t,\cdot), U(t, \cdot)] + B[U(t,\cdot), u_{\theta_0}(t, \cdot)]  &= 0 \text{ on } (0,\infty) \times \Omega
\notag \\
U(0, \cdot) &= \vartheta \sim \mathcal N_{\theta_0}.
\end{align}
The construction of the law $\mathcal N_{\theta_0}$ (given in Section \ref{limproc}) is delicate and requires the inversion of the `Fisher' information operator of the underlying statistical PDE model between suitable spaces. Even though the resulting Gaussian random field modelling the initial condition $\vartheta$ is not (almost surely) point-wise defined as a function, our proofs will imply that the weak solutions $U(t, \cdot)$ exist and almost surely define continuous functions on $\Omega$ for any $t >0$. 

\smallskip

Fix any time window $0<t_{\min} < t_{\max}<\infty$, define the separable Banach space 
\begin{equation} \label{calc}
\mathscr C \equiv C\big([t_{\min}, t_{\max}], C(\Omega)^2\big),~~\|v\|_{\mathscr C} \equiv \sup_{t\in[t_{\min}, t_{\max}], x \in \Omega, i=1,2} |v_i(t,x)|,
\end{equation}
of continuous $2$-dimensional vector fields on $[t_{\min}, t_{\max}] \times \Omega$. We also denote the Wasserstein distance between probability measures $\mu, \nu$ on a metric space $X=(X, d_X)$ by
\begin{equation}\label{wassdist}
    \mathscr W_{1,X}(\mu,\nu) = \sup_{H : X \to \R: ~\sup_{x \neq y \in X}|H(x)-H(y)|/d_X(x,y) \le 1}\Big|\int_X H(x) (d\mu(x)-d\nu(x))\Big|.
\end{equation} The symbol $\to^{P_{\theta_0}^\mathbb N}$ signifies convergence in probability under $P_{\theta_0}^\mathbb N$, while $\to^d_X$ denotes the corresponding notion of convergence in law of random variables in a metric space $X$. Here is the main result of this article,  which gives an infinite-dimensional Gaussian approximation for the updated posterior dynamical system.

\begin{Theor}\label{nstokbvm} 
Assume $\theta_0, f \in \divH^\infty, \nu>0$ and let $\alpha>12$.  Let $\mu_{N}= \mu(\cdot|Z^{(N)})$ be the conditional law in $\mathscr C$ of the stochastic process $$\Big\{\sqrt N (u_\theta(t,x) - u_{\tilde \theta_N}(t,x))|Z^{(N)}: t \in [t_{\min}, t_{\max}], x \in \Omega\Big\}$$ where $\theta \sim \Pi(\cdot|Z^{(N)})$ is drawn from the posterior measure (\ref{post}) with data $Z^{(N)}$ arising as in (\ref{model}) for $u_\theta$ solving the periodic $2D$ Navier-Stokes equations (\ref{eq:NS}), and where $\tilde \theta_N = E^\Pi[\theta|Z^{(N)}]$ is the posterior mean in $\Hnull$. Denote by $\mu$ the law in $\mathscr C$ of the Gaussian random vector field arising from the unique weak solution $U$ to the PDE (\ref{linshow}) with initial condition $\vartheta \sim \mathcal N_{\theta_0}$. Then we have as $N \to \infty$
$$\mathscr W_{1,\CCC}(\mu_{N}, \mu) \to^{P_{\theta_0}^\mathbb N} 0 ~~\text{ as well as }~~\sqrt N (u_{\tilde \theta_N} - u_{\theta_0}) \to_{\mathscr C}^d \mu.$$
\end{Theor}

 Theorem \ref{nstokbvm} implies that we can obtain `parametric' $1/\sqrt N$ convergence rates in supremum norm for inference on the Navier-Stokes regression function $u_\theta$ from standard regression data (\ref{model}) as long as the sampling horizon includes $T_0=0$. This is in surprising contrast to the situation of performing inference on the regression function $(u(t,x): t \in [0,T], x \in \Omega)$ \textit{without} the PDE constraint, where standard minimax theory (e.g., Chapter 6.3 in \cite{GinNic2016}) reveals that the rates are strictly slower than $1/\sqrt N$ over Sobolev balls. Related to this point let us remark that we have assumed that the true initial condition $\theta_0$ (but not the prior model) is smooth only to simplify some of the proofs involving the `inverse information operator'. This restriction is not necessary and can be weakened to $\alpha$-Sobolev-regularity of $\theta_0$, where $\alpha$ governs the smoothness of the prior in (\ref{prior}). In the theorem we have assumed $\alpha>12$, which demonstrates that our result holds over genuinely infinite-dimensional models but also warrants discussion of the particular condition we provide. More technical arguments likely allow to reduce the hypothesis on $\alpha$ significantly: if the prior arises from a truncated high-dimensional series expansion in a suitable basis, the posterior will inherit further regularisation properties that can be exploited in the proofs (but one looses the independence of our result on the discretisation scheme used). Moreover, if one deploys a proof strategy tailor-made to the Navier-Stokes equations instead of appealing to the general purpose Bernstein-von Mises theorem from \cite{Nickl2024}, and replacing energy estimates in Sobolev scales by Schauder type regularity theory for parabolic PDEs, a further reduction of the required bound on $\alpha$ can be attained. In this vein one can expect that the hypothesis $\alpha \ge 3$ is necessary and likely also sufficient in our setting, but the technical challenges involved in attaining this are well beyond the scope of the present paper. We believe that the main insights from our theorems and the consequences that can be drawn from them  extend to lower values of $\alpha$ than those covered by our theorem, although perhaps not to the very low regularity case.

\medskip
 
 Theorem \ref{nstokbvm} readily combines with the functional Delta-method to provide a variety of further limit theorems, for $U$ solving (\ref{linshow})
 $$
 \sqrt N (\Phi(u_{\tilde \theta_N}) - \Phi(u_{\theta_0})) 
 \stackrel{d}{\to}_\B
 \dot \Phi_{u_{\theta_0}}(U)
 ,\quad 
 N \to \infty
 ,
 $$ where $\Phi : \CCC \to \B$ is an appropriately differentiable map into a normed space $\B$ with diferential $\dot\Phi$ -- see Section \ref{functest} for details and examples.

 \smallskip

The second limit in Theorem \ref{nstokbvm} upgrades to convergence of all moments (as in \cite{Nickl2024}). Using our results and the general local asymptotic minimax theorem from mathematical statistics (Chapter 3.11 in \cite{WelVan1996}), it is further shown in the companion paper \cite{Kon25Minimax} that the limiting covariance attained in Theorem \ref{nstokbvm} is the optimal one (in the Gauss-Markov, or Cramer-Rao, sense), so that in particular no algorithm can outperform this method in a local asymptotic minimax sense.

\begin{Theor}
Consider data from (\ref{model}) where $u_\theta$ is the solution map of the periodic $2D$ Navier-Stokes equations (\ref{eq:NS}) with fixed $\theta_0, f \in \divH^\infty$ and $\nu>0$. Denote by $\mu$ the law in $\mathscr C$ of the Gaussian random vector field arising from the unique weak solution $U$ to the PDE (\ref{linshow}) with initial condition $\vartheta \sim \mathcal N_{\theta_0}$. Then,
\begin{equation}
\inf_{(\tilde u_N)} 
\sup_{\substack{J\subset \divH \cap C^\infty(\Om)^2\\ |J|<\infty}}
\liminf_{N \to \infty} 
\max_{h\in J}
N E_{\theta_0+\frac{h}{\sqrt N}}^N
\|\tilde u_N-u_{\theta_0+\frac{h}{\sqrt N}}\|^2_{\mathscr C} 
\ge 
E_{\mu} \|U\|_{\mathscr C}^2
\end{equation}
where the infimum ranges over all sequences $(\tilde u_N)$ of measurable functions $\tilde u_N$ of the observations~$Z^{(N)}$ taking values in $\mathscr C$, and the supremum ranges over all subsets $J\subset \divH^\infty$ with finite cardinality $|J|$.
\end{Theor}
 
Further details are provided in Section \ref{functest}. The estimator $u_{\tilde \theta_N}$ from Theorem \ref{nstokbvm} attains this lower bound (the required uniformity in $h$, while not stated in Theorem \ref{nstokbvm}, can be shown to follow from our proofs) and hence gives rigorous justification that the posterior measure performs optimal statistical inference in the Navier-Stokes model in the large sample limit $N \to \infty$, incorporating all relevant information provided by a noisy regression sample of the underlying nonlinear dynamics. The optimality holds for any model for the initial condition $\theta$ that is sufficiently rich to contain all smooth vector fields $\theta \in \divH^\infty$ (or in fact such that its $L^2$-closure is $\Hnull$). It remains an open question whether standard nonparametric methods that would directly estimate $u_\theta$ based on some linear smoothing method via the sample $(Y_i,t_i, X_i)_{i=1}^N$ are able to attain this asymptotic minimax constant or not.

The techniques we develop further allow to prove that posterior credible bands used in uncertainty quantification (UQ) for the filtering type distributions are proper asymptotic confidence bands: For a given significance level $0<1-\beta<1$, consider posterior credible sets for both vector entries of the velocity fields
$$
C_N(R)  
\equiv 
\Big\{u=(u_1, u_2): |u_i(t,x) - u_{\tilde \theta_N,i}(t,x))| \le R ~~\forall (t,x) \in [t_{\min}, t_{\max}] \times \Omega, i=1,2 \Big\},
$$ 
for $R>0$, and then let $C_N=C_N(R_N)$ where $R_N$ is chosen in such a way that $$\Pi(C_N(R_N)|Z^{(N)})=1-\beta.$$

\begin{Corol}\label{CorolUQ}
Under the assumptions of Theorem \ref{nstokbvm}, we have $\sqrt N R_N \to c$ in $P_{\theta_0}^\mathbb N$-probability for some $c>0$ as $N \to \infty$; in particular, ${\rm diam}(C_N)= O_{P^\mathbb N_{\theta_0}}(1/\sqrt N)$. Moreover, the credible band has exact asymptotic coverage 
$$P_{\theta_0}^\mathbb N(u_{\theta_0} \in C_N) \to 1-\beta,\spa  N \to \infty.$$
\end{Corol}

\vspace{1mm}

\vspace{3mm}

The proof follows from Theorem \ref{nstokbvm} and ideas from \cite{CN13}, and is given in Section \ref{corprf}. One can further show uniformity in $\theta_0$ in the above result but the details are technical and will not be demonstrated here.

\smallskip

Let us conclude with a remark about computation of the posterior measure. For this, one would typically discretise $\theta$ into a high-dimensional approximation space $E_D \simeq \R^D$ such as linear combintations of the first $D$ eigenfunctions of the Stokes operator. Provably reliable computation of the estimate $u_{\tilde \theta_N}$ can be expected to be possible in polynomial time in relevant parameters (such as $D,N$) by an MCMC -based algorithm, combining ideas developed in \cite{NW24} with the analytical estimates obtained here -- see also Chapter 5 in \cite{NicklEMS} for a more general framework. An initialiser for this algorithm can be obtained from the boundary trace $\hat u(0,\cdot)$ of a standard nonparametric smoothing estimator $\hat u = \hat u ((Y_i, t_i, X_i)_{i=1}^N)$ of $u_\theta$. The details will be studied elsewhere. Let us finally remark that the approach just described is based on log-concave (but non-Gaussian) approximations. In computational mathematics, `Laplace' approximations by normal distributions are often used; these can be justified, at least heuristically, from the Gaussian limits we obtain in the present paper, but their overall validity remains unclear.

\subsection{Notation and preliminaries}\label{sec:Notation}

Throughout, we take as spatial domain $\Om=\T^2$ the $2$-dimensional torus with unit length, \ie $\Om=[0,1]^2/\sim$ where~$\sim$ identifies opposite sides of $[0,1]^2$. The Euclidean inner product on $\R^2$ or $\T^2$ is denoted $u\cdot v=\sum_{i=1}^2 u_i v_i$ for vectors $u=(u_1,u_2)$, $v=(v_1,v_2)$, with corresponding Euclidean norm $|u|$. For a matrix $A=[A_{ij}]$ we write its transpose $A^T=[A_{ji}]$ and its norm $|A|^2=\sum_{i,j} |A_{ij}|^2$. 
We write $a_k\lesssim b_k$ (resp. $\gtrsim$) when there exists a positive constant $C$ such that $a_k\leq Cb_k$ (resp.~$\ge$) for all $k$. For any function $v:\Om\to\R^2$, we let $[v]_j$ or $v_j$ denote the $j$th component of $v$, $j=1,2$.

For any $k\in \N$, the space $C^k(\Om)^2$ is the collection of functions $v:\Om\to\R^2$ obtained as the restriction to $[0,1]^2$ of some $1$-periodic $\tilde{v}:\R^2\to\R^2\in C^k(\R^2)^2$ in the sense that $\tilde{v}(x+e_j)=\tilde{v}(x)$ for all $x\in\R^2$ and canonical direction $e_j\in\R^2$. We further define $C^\infty_{\rm}(\Om)^2\equiv \cap_{k\geq 0} C^k(\Om)^2$. For $v:\Om\to\R^a$, we write $\partial_i v$ for the partial derivative $(\partial v)/(\partial x_i)$ acting component-wise as $[\partial_i v]_j = \partial_i v_j$, and $\nabla v$ for the component-wise gradient of $v$: when $a=1$ then $\nabla v$ is the usual gradient of $v$, whereas if $a\geq 2$ then $\nabla v: \Om\to \R^{a\times 2}$ is the Jacobian matrix $[\nabla v]_{ij} = \partial_j v_i$, $i=1,2$ and $j=1,2,\ldots a$. Similarly, for $v:\Om\to\R^a$ then $\Delta v$ is Laplacian of $v$ given in components by $[\Delta v]_i=\Delta v_i = \sum_{j=1}^2 \partial^2_j v_i$. When $v:\Om\to\R^2$, the divergence of $v$ is $\nabla\cdot v = \sum_{i=1}^2 \partial_i v_i$.

We define the usual Lebesgue spaces $L^p(\Om)^2$, $1\leq p\leq +\infty$, for Lebesgue measure $dx$ which, by convention, consist of \textit{real-valued} functions unless otherwise specified. For complex-valued square-integrable vector-fields $u,v:\Om\to\C^2$, we let 
\begin{equation}\label{eq:complexL2}
\ps{u}{v}_{L^2}
=
\int_\Om u(x)\cdot \overline{v(x)}\, dx 
.
\end{equation}
Since the complex exponentials 
$$
h_k(x)\equiv e^{2i\pi k\cdot x}
,\spa 
k\in\Z^2
,\quad 
x\in\Om 
,
$$
form an orthonormal basis of \textit{complex-valued} periodic square integrable functions, then 
\begin{equation}\label{eq:L2Series}
L^2(\Om)^2 
=
\Big\{ \sum_{k\in\Z^2} (\alpha_k, \beta_k) h_k(x) : \alpha_{-k} = \overline \alpha_k,\ \beta_{-k} = \overline \beta_k,\ (\alpha_k), (\beta_k)\in\ell^2(\Z^2)\Big\} 
,
\end{equation}
where the conjugacy constraints mean that elements of $L^2(\Om)^2$ are real-valued. The complex exponentials $\{e^{2i\pi k x} : k\in\Z^2\}$ are the eigenfunctions of $-\Delta$ with corresponding eigenvalue $4\pi|k|^2$.
Let us enumerate $\Z^2$ by $\{k_j : j\geq 0\}$ in such a way that $j\mapsto |k_j|$ is nondecreasing and
\begin{equation}\label{eq:enum}
\lambda_0 \equiv 0,\spa 
\lambda_j \equiv
4\pi|k_j|^2
\asymp 
j
.
\end{equation}
Then, for $\alpha\in [0,\infty)$ define the usual Sobolev space 
$$
H^\alpha(\Om)^2
=
\Big\{ 
v:\Om\to\R^2 
: 
\|v\|^2_{H^\alpha}
<
\infty\Big\} 
,
$$
where the $\|\cdot\|_{H^\alpha}$-norm is the one induced by the inner product
\begin{equation}\label{SobMetricIntro}
\ps{u}{v}_{H^\alpha}
\equiv 
\sum_{j\geq 0} (1+\lambda_j^\alpha) \ps{u}{h_{k_j}}_{L^2}\cdot \overline{\ps{v}{h_{k_j}}}_{L^2}
,\spa 
\ps{u}{h_{k_j}}_{L^2}
\equiv 
\int_\Om u(x) e^{-2i\pi k_j\cdot x}\, dx 
\in\R^2
.
\end{equation}
We define the homogeneous Sobolev $\dH^\alpha(\Om)^2$ (semi) inner product 
\begin{equation}\label{dotSobIntro}
\ps{u}{v}_{\dH^\alpha}
\equiv 
\sum_{j\geq 1} 
\lambda_j^\alpha  
\ps{u}{h_{k_j}}_{L^2}\cdot \overline{\ps{v}{h_{k_j}}}_{L^2}
.
\end{equation}
The resulting $\|\cdot\|_{\dH^\alpha}$-norm is only a seminorm on $H^\alpha(\Om)^2$ since constants functions (\ie those associated with the eigenvalue $\lambda_0=0$) have zero $\|\cdot\|_{\dH^\alpha}$-norm. Notice that 
$$
\|u\|^2_{H^\alpha}
=
\|u\|^2_{L^2} + \|u\|^2_{\dH^\alpha}
,\spa 
\forall\ u\in H^\alpha(\Om)^2 
.
$$
When $\int_\Om u(x)\, dx = 0$, \ie $\ps{u}{h_{k_0}}_{L^2}=0$, then the $\|\cdot\|_{H^\alpha}$ and $\|\cdot\|_{\dH^\alpha}$ norms are equivalent: we have
\begin{equation}\label{equivSobIntro}
\|u\|_{\dH^\alpha}
\leq 
\|u\|_{H^\alpha}
\leq 
(1+\lambda_1^{-\alpha}) \|u\|_{\dH^\alpha}
;
\end{equation}
see Section \ref{sec:setting}. The natural space to study the system of equations (\ref{eq:NSequations}) will be the linear subspace of $L^2$ consisting of divergence-free and mean-zero vector fields, \ie 
$$
\Hnull
\equiv
\Big\{ 
u\in L^2(\Om)^2 
:
\nabla\cdot u=0,\ \int_\Om u = 
0
\Big\}
,
$$
where $\nabla\cdot u$ is to be understood in a distributional sense. Then $(\Hnull,\ps{\cdot}{\cdot}_{L^2})$ is a \textit{real} Hilbert space. An orthonormal basis for the complex-valued extension of $\Hnull$ is obtained (see Section~\ref{sec:setting}) as the subset $\{e_j : j\geq 1\}$ of eigenfunctions of the periodic Laplacian 
\begin{equation}\label{defej}
e_j(x) 
\equiv
c_j
h_{k_j}(x)
,\spa 
c_j 
=
\frac{(-k_{j,2},k_{j,1})}{|k_j|}
,\quad 
j\geq 1
.
\end{equation}
Then any $v\in \Hnull$ decomposes as 
$$
v(x)
=
\sum_{j\geq 1} v_j e_j(x) 
,\spa 
v_j = \ps{v}{e_j}_{L^2}
=
\int_\Om v(x) \cdot \overline{e_j(x)}\, dx 
,
$$
and the fact that $v$ is real-valued translates into conjugacy relations on the $v_j$'s as in (\ref{eq:L2Series}).
We finally define the associated Sobolev spaces $\divH^\alpha$ compatible with the structure of $\Hnull$ as the closed subspaces 
\begin{equation}\label{defdivHalpha}
\divH^\alpha 
\equiv 
H^\alpha(\Om)^2 \cap \Hnull
,\spa 
\alpha\in [0,\infty) 
.
\end{equation}
We endow $\divH^\alpha$ with the $\dH^\alpha$-topology or, equivalently, that of $H^\alpha(\Om)^2$ by virtue of (\ref{equivSobIntro}) since all functions in $\divH^\alpha$ have zero average. We further define the smooth functions of $\Hnull$ as 
$$
\divH^\infty 
\equiv
C^\infty(\Om)^2 \cap \Hnull
=
\cap_{\alpha>0} \divH^\alpha
.$$
Since we will always be working with $\R^2$-valued functions, from here on we will often write for simplicity 
$$
L^2 \equiv L^2(\Om)^2
,\spa
H^\alpha \equiv H^\alpha(\Om)^2
,\spa 
\dH^\alpha \equiv \dH^\alpha(\Om)^2 
,
$$
$$
C^\infty\equiv C^\infty(\Om)^2
,\spa 
H^\infty \equiv \cap_{\alpha>0} H^\alpha 
,\spa 
\dH^\infty \equiv \cap_{\alpha>0} \dH^\alpha
.
$$ For $\alpha<0$, we define  $\divH^\alpha$ as the completion of $\Hnull$ (or, equivalently, $\divH^\infty$) with respect to the norm induced by the inner product (\ref{dotSobIntro}). 

Since $(\Hnull,\ps{\cdot}{\cdot}_{L^2})$ is a closed linear space of $L^2$, one can define the so-called `Leray' projector, \ie the $L^2$ projection operator 
\begin{equation}\label{defP}
P: L^2\to \Hnull
.
\end{equation}
The unknown quantities to solve for in the Navier-Stokes equation (\ref{eq:NSequations}) are the velocity field $u:\Om\to\R^2$ and the scalar pressure $p:\Om\to \R$. As is common in the literature, we study the `projected' equation (\ref{eq:NSequations}) by applying $P$ to it. This leads to an equivalent formulation in functional form where the velocity field, as a map $u:[0,T]\to \Hnull$, is the unique solution $u=u_\theta$ to the time evolution equation
\begin{equation}\label{eq:NS}
    \frac{du}{dt} + \nu A u + B[u,u]  = f
    ,
    \spa 
    u(0)=\theta.
 \end{equation}
Here $A=-P\Delta = -\Delta$ is the Stokes operator (with the second identity holding on $\Hnull$ since $\Om$ is the two-dimensional torus, see (\ref{eq:PCommutesLaplacians})), and $B$ is the bilinear form 
 \begin{equation}\label{eq:defB}
 B[u,v]
 \equiv 
 P[(u\cdot \nabla) v]
 ,\spa 
 \forall\ u,v\in C^\infty
 ,
 \end{equation}
for $(u\cdot \nabla)v$ the matrix product $(\nabla u) v$ or, equivalently, the  vector-field with entries 
 $$
 [(u\cdot \nabla)v]_i
 =
 u\cdot \nabla u_i 
 =
 \sum_{j=1}^2 u_j \partial_j v_i 
 .
 $$ The existence of unique (strong) solutions to (\ref{eq:NS}) is reviewed in Proposition \ref{PropNSExist} below.

We finally introduce some more function spaces required below: For any measurable subset $U\subset\R^d$, $d\geq 1$, and Banach space $(\B,\|\cdot\|_\B)$ we let for any measurable $v:U\to \B$ and $p\in [1,\infty)$
$$
\|v\|_{L^p(U,\B)}
=
\bigg( \int_U \|v(x)\|^p_\B\, dx\bigg)^{1/p}
.
$$
We define $L^p(U,\B)$ as the collection of maps $v:U\to \B$ with finite $\|\cdot\|_{L^p(U,\B)}$-norm. When $(\B,\ps{\cdot}{\cdot}_{\B})$ is a Hilbert space, then $L^2(U, \B)$ is itself a Hilbert space with inner product 
$$
\ps{u}{v}_{L^2(U,\B)}
=
\int_U \ps{u(x)}{v(x)}_{\B}\, dx
.
$$ 
Similarly, one defines spaces $C(U, \B)$ of continuous maps $v: U \to \B$, normed by the supremum norm $\sup_{t \in U}\|v(t)\|_\B$.
For instance, we will often consider norms and spaces of the type
$$
\|v\|^2_{L^2([0,T],\dH^\alpha)} 
=
\int_0^T \|v(t)\|^2_{\dH^\alpha}\, dt
,~~ \|v\|_{C([0,T], \divH^\alpha)} = \sup_{0<t<T}\|v(t)\|_{\divH^a}.
$$ 


\section{Functional Bernstein-von Mises Theorem}
\label{sec:FunctionalBvM}

The proofs of Theorem \ref{nstokbvm} and its consequences are organised into several subsections, including an initial posterior contraction rate, the inversion of the information operator of the statistical model, the construction of the limit process and convergence in function space, as well as applications to functional estimation.


\subsection{Lipschitz stability estimate and faster posterior contraction}\label{sec:ContractionRates}


A key step in the proof of our main results is to  initially establish polynomial (i.e., faster than logarithmic) contraction rates for the initial condition, based on a new stability estimate for the forward map that exploits availability of `samples' near $t=0$. This in turn allows to use more local `BvM'-type techniques from \cite{CN14} in the proofs. The proof for this new stability estimate is different from the corresponding result for reaction diffusion equations (see (102) in \cite{Nickl2024}). There, one exploits the symmetry of the approximate linearisation (Schr\"odinger) operator but this is not possible in the Navier-Stokes setting. Rather, we use a suitable parabolic time-continuity argument at zero that relies on estimates for the bilinear term $B$. A similar result would hold true for general parabolic equations with a sufficiently `mild' nonlinearity.

For the theory that follows we will employ the Gaussian process prior from (\ref{prior}), which satisfies Condition 1 in \cite{NicTit2024} with RKHS $\mathcal H$ there equal to $\Hnull^\alpha$. The proof of Theorem 3 in \cite{NicTit2024} -- which is based on key ideas developed for a different PDE in \cite{MNP21a} as well as, by now well developed, techniques from Bayesian nonparametric statistics \cite{GV17} -- together with Theorem 2.2.2 in \cite{NicklEMS} with $d=2$, $\alpha \geq 2$, any $\beta<\alpha-1$ and $\kappa=0$ (see Proposition 1 in \cite{NicTit2024}) then entails that for 
    \begin{equation}\label{RateDeltaN}
    \delta_N=N^{-\frac{\alpha}{2\alpha+2}}
    \end{equation}
    and all $M$ large enough we have
    \begin{equation}\label{fwdrat}
    \Pi\Big( \theta \in \divH^\beta : \|\theta\|_{\dH^\beta}\leq M,\ \int_0^T\|u_{\theta}(t)-u_{\theta_0}(t)\|^2_{L^2}\, dt \le M \delta_N\ |\ Z^{(N)} \Big) 
    \to 
    1
    \end{equation}
    in $P_{\theta_0}^\N$-probability as $N\to\infty$ (which is a slightly upgraded version of (39) in \cite{NicTit2024}). Combined with the following stability estimate, this yields a fast contraction rate at the $\theta$-level in Theorem \ref{TheorFastRateTheta} below.

    \begin{Prop}\label{PropStabForward}
        For any $f\in\divH^1$, $T>0, \nu>0,L>0$ and $\theta_0, \theta\in\divH^2$, let $u_\theta, u_{\theta_0}$ denote the  corresponding solutions  to the periodic $2D$ Navier-Stokes equations (\ref{eq:NS}). Then there exists $C=C(\nu, T, \|f\|_{\dH^1}, L)>0$ such that if $\|\theta\|_{\dH^2}+\|\theta_0\|_{\dH^2}\le L$  we have 
        $$
        \|\theta-\theta_0\|^2_{\dH^{-1}}
        \leq 
        C \int_0^T \|u_\theta(t)-u_{\theta_0}(t)\|^2_{L^2} dt
        .
        $$
    \end{Prop}

    \vspace{1mm}

    \begin{Proof}{Proposition \ref{PropStabForward}}
        Let $U\equiv u_{\theta}-u_{\theta_0} \in \Hnull$ and observe that $U$ satisfies the equation 
        $$
        \frac{dU}{dt} - \nu \Delta U + B[U, u_\theta] + B[u_{\theta_0}, U]
        =
        0
        ,\spa 
        U(0)=\theta-\theta_0
        .
        $$
        For $a=-1$ and $a^*=2$ as in (\ref{eq:DefinAStar}), since $f\in \divH^1=\divH^{a^*-1}$ and $\theta_0,\theta\in\divH^{a^*}$, Proposition \ref{PropLinPStab} with $a=-1$, $\xi=\theta-\theta_0$, and $g=0$ entails that 
        $$
        \|\theta-\theta_0\|^2_{\dH^{-1}}
        \le
        C\textbf{}\int_0^T \|U(t)\|^2_{L^2}\, dt
        ,
        $$
        for some constant $C=(\nu, T,\|f\|_{\dH^1},L)>0$.
    \end{Proof}

    \vspace{3mm}

     We see that when the sampling time horizon is $[0,T]$ (rather than $[T_0,T]$ for $T_0>0$), one can obtain Lipschitz stability estimates that improve on the lower bound in Theorem 2 of \cite{NicTit2024} for $T_0>0$. The next result leverages this and shows that faster than inverse logarithmic in $N$ contraction rates for $\theta$ can be obtained in this case. When compared to the minimax lower bound Theorem 4 in \cite{NicTit2024}, this reveals how crucial statistical information is `lost' when no samples near time $t=0$ are available -- this can be explained by the dissipative nature of the underlying Navier-Stokes system.

\begin{Theor}\label{TheorFastRateTheta}
		Let $\Pi$ be a prior as in (\ref{prior}), let $\beta \in (0, \alpha-1)$, let $Z^{(N)}$ be data obtained as in (\ref{model}) where $u_\theta$ solves (\ref{eq:NS}) for some $f \in \divH^1, \nu>0$, and denote by $\Pi(\cdot|Z^{(N)})$ the resulting posterior distribution. Assume that $\theta_0\in \divH^\alpha$. Then, for all $M>0$ large enough we have as $N \to \infty$ and in $P_{\theta_0}^\N$-probability that
		$$
		\Pi\Big(\theta\in  \divH^\beta : 
        \|\theta\|_{\dH^\beta}\leq M,\ 
        \|\theta-\theta_0\|_{L^2}\leq M\delta_N^{\frac{\beta}{\beta+1}}\ |\ Z^{(N)} \Big)
		\to 1.$$
        
    \end{Theor}

    The proof follows from (\ref{fwdrat}),  Proposition \ref{PropStabForward}, and interpolation of $L^2$ between $
\divH^{-1}$ and $\divH^\beta$, and in fact works for $\alpha \ge 2$.

\subsection{Inverting the information operator}\label{sec:InvertI}

In this section we establish the key result that the  \textit{(Fisher) information operator} underlying our measurement model (\ref{model}), 
\begin{equation} \label{infoop}
\frac{1}{T}\I^*_{\theta} \I_\theta : \Hnull\mapsto \Hnull,
\end{equation}
is a homeomorphism between $\divH^a$ and $\divH^{a+2}$ for any integer $a\ge -1$ (Theorem \ref{TheorInverseI} below). To define it, we follow Definition 3.1.2 in \cite{NicklEMS} and note first that, for smooth enough $\theta$, the solution map $\theta\mapsto u_\theta$ of the PDE (\ref{eq:NS}) admits a continuous linearisation
$$
\I_{\theta}:\Hnull\mapsto L^2([0,T],\Hnull)
$$
in the sense that $u_{\theta+h}-u_\theta - \mathbb I_\theta[h]$ is $o(\|h\|)$ in appropriate norms for any direction $h \in \Hnull$ -- see Proposition \ref{PropLinQuad} in the Appendix. This linear map is hence the usual `score operator' in our statistical regression model (\ref{model}) occurring also in its LAN expansion  (see, e.g., Section~3.1.1 in \cite{NicklEMS} and Section 2 in \cite{Kon25Minimax}). Specifically, for any $h\in \divH^1$, the function $\I_{\theta}[h]$ is the solution $U:[0,T]\to \Hnull$ to the linear PDE
\begin{align}\label{eq:LinPDE}
\frac{dU}{dt}
-
\nu \Delta U(t)
+B[u_\theta(t), U(t)] + B[U(t),u_\theta(t)]
=
0
,\spa 
U(0)=h
,
\end{align}
where we recall that, for $P$ the Leray projector from (\ref{defP}), we have $B[u,v]=P[(u\cdot \nabla)v]$. For smooth inital conditions $h\in \divH^\infty\equiv C^\infty \cap \Hnull$, parabolic regularity estimates (see Proposition \ref{PropLinSob}) for (\ref{eq:LinPDE}) provide refined mapping properties for $\I_\theta$: we have
\begin{equation}\label{mapI}
\|\I_{\theta}[h]\|_{L^2([0,T],\divH^{a+1})}
\lesssim 
\|h\|_{\dH^a}
,\spa 
\forall\ h\in \divH^\infty
,\quad 
\forall\ a\in\Z
.
\end{equation}
In particular, since $\divH^\infty$ is dense in $\divH^a$ for any $a\in\Z$, one extends $\I_\theta$ to a continuous linear $L^2([0,T], \divH^{a+1})$-valued map defined on all of $\divH^a$. In addition,  when $h\in\divH^{-a}$ for some $a>0$, then $\I_\theta[h]$ is a weak solution to (\ref{eq:LinPDE}) in the sense that the PDE is satisfied for almost all $t>0$ when tested against an arbitrary $\varphi\in \divH^a$. 

Denote by $\I_{\theta}^*$ the adjoint operator of $\I_{\theta}:\Hnull\to L^2([0,T],\Hnull)$. The mapping properties (\ref{mapI}) of $\I_\theta$ transfer to $\I_\theta^*$; specifically, $\I_\theta^*$ extends to a continuous operator $\I_\theta^* : L^2([0,T],\divH^a)\to \divH^{a+1}$ for any $a\in\Z$ (arguing, e.g., as in the proof of Lemma \ref{PropCompactPerturbation} below).  It is clear from the previous remarks that $\I_\theta^* \I_\theta : \divH^a\to \divH^{a+2}$ defines a continuous linear operator for any $a\in\Z$, and our main result can thus be considered as an infinite-dimensional analogue of the existence of the inverse Fisher information of the underlying statistical model. This will in particular allow us to demonstrate the existence and tightness of the law $\mathcal N_{\theta_0}$ from~(\ref{linshow}), which plays a pivotal role in Theorem \ref{nstokbvm}. It adds to a still relatively short list of examples of PDE models where the information operator can be inverted between relevant function spaces, see \cite{MNP19, Nickl2020, NR20, MNP21, NicklEMS, Nickl2024}.

To prove the result, we follow ideas laid out in Section 3.4 in \cite{Nickl2024} by comparing the solution $\I_\theta[h]$ of (\ref{eq:LinPDE}) to a solution of the standard heat equation in $\Hnull$: the difference of these operators will be seen to define a compact operator so that Fredholm theory can be applied. To this end, for any $h\in \Hnull$ let $\LL[h]$ denote the solution $U$ to the PDE 
\begin{equation}\label{eq:heat}
\frac{dU}{dt}
-
\nu \Delta U(t)
=
0
,\spa 
U(0)=h 
.
\end{equation}
The representation formula (\ref{repheat}) together with the parabolic estimate (\ref{heatestimate}) with $g=0$, entail that $\LL[h]$ satisfies the divergence-free and mean-zero constraints when $h$ does and, furthermore, defines a continuous linear operator 
$$
\LL:\divH^a\to L^2([0,T], \divH^{a+1})
,\spa 
\forall\ a\in \Z
.
$$
In particular, $\LL$ is a bounded linear operator $\LL: \Hnull\mapsto L^2([0,T],\Hnull)$, with associated adjoint 
$$
\LL^* 
: 
L^2([0,T],\Hnull)
\to 
\Hnull
.
$$
As previously, the mapping properties of $\LL$ transfer to $\LL^*$ so that $\LL^*$ defines a continuous linear operator $\LL^* : L^2([0,T], \divH^{a})\to \divH^{a+1}$ for any $a\in\Z$. In particular, $\LL^* \LL : \divH^a\to \divH^{a+2}$ defines a continuous linear operator for any $a\in\Z$.

The strategy of the proof consists in decomposing the information operator $\I_\theta^* \I_\theta$ as 
\begin{align}\label{eq:IstarI}
\I_\theta^* \I_\theta 
=
\LL^* \LL + \big\{ \I_\theta^* (\I_\theta-\LL) + (\I_\theta-\LL)^* \LL \big\}
\equiv
\LL^* \LL + \KK
.
\end{align}
In what follows, we will show that $\LL^* \LL$ is a compact perturbation of a constant multiple of $(-\Delta)^{-1}$ and that the bracketed terms in (\ref{eq:IstarI}), \ie the linear operator $\KK$, map $\divH^a$ compactly into $\divH^{a+2}$. Since $\I_\theta^* \I_\theta$ will also be shown to be injective, this entails, by Fredholm theory, that $\I_\theta^* \I_\theta$ has the same mapping properties as $(-\Delta)^{-1}$ and, in particular, that it is also surjective and hence invertible between appropriate spaces.

\smallskip

As a last preparation, define the heat-semigroup $\{S_t : t\geq 0\}$ associated with (\ref{eq:heat}) 
$$
S_t[h]
\equiv
\LL[h](t)
,\spa 
\forall\ h\in \divH^\infty
.
$$
Recall from (\ref{defej}) that $\{e_j : j\geq 1\}$ is an $L^2$-orthonormal basis of (the complex extension of) $\Hnull$ consisting of eigenfunctions of $-\Delta$ with corresponding eigenvalues $\lambda_j\asymp j$. Then, a standard Fourier decomposition of (\ref{eq:heat}) yields the well-known spectral representation
\begin{align}\label{eq:specrep}
    S_t[h] 
    \equiv 
    \sum_{j \geq 1} e^{-\nu\lambda_jt} \ps{h}{e_j}_{L^2} e_j
    ,\spa 
    \forall\ h\in \divH^\infty
    ,\ t\geq 
    0
    .
\end{align}
It follows that for $h\in\divH^a$, $a\in\Z$, and positive $t>0$, $S_t[h]$ is well-defined. Moreover, $S_t$ is self-adjoint on $\Hnull$, \ie 
$$
\ps{S_t[h]}{\phi}_{L^2}
=
\ps{h}{S_t[\phi]}_{L^2}
,\spa 
\forall\ h,\phi\in \Hnull 
.
$$
Given the spectral definition of the Sobolev spaces $\divH^a$, $a\in\R$, in (\ref{defdivHalpha}) and (\ref{dotSobIntro}), one easily sees that, for any \textit{positive} time $t>0$ and $a,b\in\R$, we have 
\begin{equation}\label{eq:HeatSobolev}
\|S_t[h]\|_{\dH^b} 
\lesssim 
\|h\|_{\dH^a}
,\spa 
\forall\ h\in\divH^a
.
\end{equation}
Therefore, as soon as $t>0$, then $S_t$ maps any $\divH^a$ with $a\in\R$ into $\divH^\infty$; in particular, any restriction $S_t : \divH^a\to \divH^b$, for some $a,b\in\R$, defines a compact operator. 

\begin{Prop}\label{PropLstarL}
    Fix $a\in\R$. Then $\LL^* \LL$ maps $\divH^a$ continuously into $\divH^{a+2}$ and we have the equality
    $$
     \Delta \LL^* \LL = \frac{1}{2\nu}\big( S_{2T} - {\rm Id}_{\divH^a}\big)
    $$
    as linear operators on $\divH^a$.
\end{Prop}

\begin{Proof}{Proposition \ref{PropLstarL}}
By (\ref{eq:specrep}), we have $S_t[e_j]=e^{-\nu\lambda_jt}e_j$ for all $t\geq 0$ and $j\geq 1$. Also recall from the comments preceding (\ref{eq:IstarI}) that $\LL^* \LL : \divH^a \to \divH^{a+2}$ defines a continuous linear operator. Using that $S_t$ is self-adjoint on $\Hnull$ and the semi-group property $S_t\circ S_t = S_{2t}$, we then have for all $h\in \divH^\infty$
\begin{eqnarray*}
    \lefteqn{
        \LL^* \LL[h]
        =
        \sum_{j\geq 1} \ps{\LL^* \LL [h]}{e_j}_{L^2}e_j 
        =
        \sum_{j\geq 1} \ps{\LL [h]}{\LL[e_j]}_{L^2([0,T],\Hnull)}e_j  
    }
    \\&&
    =
    \sum_{j\geq 1} \Big(\int_0^T \ps{S_t[h]}{S_t[e_j]}_{L^2}\, dt \Big) e_j
    =
    \sum_{j\geq 1} \Big(\int_0^T e^{-2\nu\lambda_jt} \ps{h}{e_j}_{L^2}\, dt \Big) e_j
    \\&&
    =
    \frac{1}{2\nu} \sum_{j\geq 1} \lambda_j^{-1} (1-e^{-2\nu\lambda_j T}) \ps{h}{e_j}_{L^2} e_j
    =
    \frac{1}{2\nu} (-\Delta)^{-1}\big({\rm Id} - S_{2T}\big)[h]
    ,
\end{eqnarray*}
where we used Fubini's theorem, the identity holding in $\Hnull^{a+2}$. Since $(-\Delta): \divH^{a+2}\to \divH^a$ is a homeomorphism for all $a\in\R$, since $S_{2T} : \divH^a\to \divH^b$ is continuous for any $a,b\in\R$ by (\ref{eq:HeatSobolev}), and since $\divH^\infty$ is dense in any $\divH^a$, the result follows.
\end{Proof}
\vspace{3mm}

We now turn to mapping properties of the bracketed terms in (\ref{eq:IstarI}). The next lemma establishes that, although $\I_\theta$ and $\LL$ are each $1$-smoothing, their difference, however, is $2$-smoothing, a result that follows from parabolic regularity estimates.

\begin{Lem}\label{LemDiffUI}
    Fix $a\in\Z$ and let $b=\max\{2,|a|+1\}$.
    Assume that $\theta\in\divH^{b}$ and $f\in\divH^{b-1}$.
    Then the linear operator $\I_\theta-\LL$ maps $\divH^a$ continuously into $L^2([0,T],\divH^{a+2})$.
\end{Lem}

\begin{Proof}{Lemma \ref{LemDiffUI}}
Fix $h \in \divH^\infty$ and define $w\coloneqq (\LL - \I_\theta)[h]$. Then, $w$ solves the PDE
$$
\frac{dw}{dt} - \nu \Delta w
+ B\big[u_\theta,\I_\theta[h]\big] + B\big[\I_\theta[h],u_\theta\big]
=0,\spa 
w(0) = 0
,
$$
by Proposition \ref{PropExistLin}(ii).
Standard regularity estimates for the periodic heat equation (see~(\ref{heatestimate})) yield
$$
\|w\|_{L^2([0,T], \dH^{a+2})}
\lesssim 
\big\|B\big[u_\theta,\I_\theta[h]\big] + B\big[\I_\theta[h],u_\theta\big]\big\|_{L^2([0,T], \dH^a)}
.
$$
Now, using Proposition \ref{PropProdSob} for the bilinear form $B[\cdot,\cdot]$, the Sobolev bound on $u_\theta$ from Proposition \ref{PropNSBound}, and the Sobolev estimate for $\I_\theta$ in Proposition \ref{PropLinSob} (noticing that $b\ge |a|+1\geq a^*=|a|+\mathbbm{1}_{\{|a|\leq 1\}}$), we have
\begin{eqnarray*}
\lefteqn{
\|w\|^2_{L^2([0,T], \dH^{a+2})}
\lesssim
\int_0^T \big\|B\big[u_\theta,\I_\theta[h]\big](t) + B\big[\I_\theta[h],u_\theta\big](t)\big\|^2_{\dH^a}\, dt
}
\\
&&
\lesssim 
\int_0^T \|u_\theta(t)\|^2_{\dH^b} \|\I_\theta[h](t)\|^2_{\dH^{a+1}}\, dt 
\lesssim
\big(1+\|\theta\|^{2b}_{\dH^b}\big)\int_0^T \|\I_\theta[h](t)\|^2_{\dH^{a+1}}\, dt 
\lesssim 
\|h\|^2_{\dH^a}
.
\end{eqnarray*}
We deduce that 
$$
\|(\LL-\I_\theta)[h]\|_{L^2([0,T], \dH^{a+2})}
\lesssim 
\|h\|^2_{\dH^a}
,\spa 
\forall\ h\in \divH^\infty
.
$$
The result follows by approximation.
\end{Proof}
\vspace{3mm}

    

We can now use the preceding lemma to establish that the perturbation $\KK$ in (\ref{eq:IstarI}) maps $\divH^a$ into $\divH^{a+2}$ compactly.

\begin{Lem}\label{PropCompactPerturbation}
    Fix $a\in\Z$ and let $b=\max\{a+4, |a|+1\}$. Assume that $\theta\in\divH^b$ and $f\in\divH^{b-1}$. Then the linear map $\KK$ from (\ref{eq:IstarI}) defines a compact operator $\KK:\divH^a\to \divH^{a+2}$.
\end{Lem}

\begin{Proof}{Lemma \ref{PropCompactPerturbation}}
    Let us first observe that, since $a\in\Z$, then 
    $$
    b
    =
    \max\{a+4, |a|+1\} 
    =
    \max\{2, |a|+1, |a+3|+1\} 
    .
    $$
    Since $b\geq \max\{2,|a|+1\}$, Lemma \ref{LemDiffUI} entails that $\I_\theta-\LL$ defines a continuous operator 
    $$
    \I_\theta - \LL : \divH^a \to L^2([0,T],\divH^{a+2})
    .
    $$
    Now define $a_1=-(a+3)$ and $a_1^*$ as in (\ref{eq:DefinAStar}). Because $b\geq |a+3|+1\geq a_1^*$, then Proposition \ref{PropLinSob} implies that $\I_\theta$ maps $\divH^{-(a+3)}$ continuously into $L^2([0,T],\divH^{-(a+2)})$. Consequently, $\I_\theta^*$ defines a continuous operator
    $$
    \I_\theta^* : 
    L^2([0,T],\divH^{a+2})\to \divH^{a+3}
    ,
    $$
    since we have for any $h\in \divH^\infty$
    \begin{eqnarray*}
    \lefteqn{
    \|\I_\theta^*[h]\|_{\dH^{a+3}}
    =
    \sup_\psi \ps{\I_\theta^*[h]}{\psi}_{L^2(\Om)^2}
    =
    \sup_\psi \ps{h}{\I_\theta[\psi]}_{L^2([0,T],\Hnull)}
    }
    \\&&
    \leq 
    \sup_\psi 
    \int_0^T \ps{h(t)}{\I_\theta[\psi](t)}_{L^2(\Om)^2}\, dt
    \leq 
    \sup_\psi
    \int_0^T \|h(t)\|_{\dH^{a+2}} \|\I_\theta[\psi](t)\|_{\dH^{-(a+2)}}\, dt
    \\&&
    \leq 
    \sup_\psi
    \Big(\int_0^T \|h(t)\|^2_{\dH^{a+2}}\, dt\Big)^{1/2} 
    \Big(\int_0^T \|\I_\theta[\psi](t)\|^2_{\dH^{-(a+2)}}\, dt\Big)^{1/2}
    \\&&
    =
    \|h\|_{L^2([0,T], \dH^{a+2})} \sup_\psi \|\I_\theta[\psi]\|_{L^2([0,T], \dH^{-(a+2)})}
    \\&&
    \lesssim 
    \|h\|_{L^2([0,T], \dH^{a+2})}
    \sup_\psi \|\psi\|_{\dH^{-(a+3)}}
    =
    \|h\|_{L^2([0,T], \dH^{a+2})}
    ,
    \end{eqnarray*}
    where the supremum is taken over all $\psi\in \divH^\infty$ such that $\|\psi\|_{\dH^{-(a+3)}}\leq 1$. It follows that 
    $$
    \I_\theta^*(\I_\theta-\LL)
    :
    \divH^a
    \to 
    \divH^{a+3}
    $$
    is a continuous linear operator. Recalling that $a_1=-(a+3)$ and since $b\geq \max\{2, |a_1|+1\}$, Lemma \ref{LemDiffUI} entails that $\I_\theta-\LL$ defines a continuous operator 
    $$
    \I_\theta-\LL 
    :
    \divH^{-(a+3)} 
    \to 
    L^2([0,T],\divH^{-(a+1)})
    .
    $$
    Reasoning as previously yields the well-posedness and continuity of
    $$
    (\I_\theta-\LL)^*
    :
    L^2([0,T],\divH^{a+1}) 
    \to 
    \divH^{a+3}
    .
    $$
    Recalling that $\LL$ maps $\divH^a$ continuously into $L^2([0,T],\divH^{a+1})$ by (\ref{heatestimate}), we see that 
    $$
    (\I_\theta-\LL)^* \LL 
    :
    \divH^a \mapsto \divH^{a+3}
    ,
    $$
    is a continuous operator. Consequently, the operator $\KK$ maps $\divH^a$ continuously into $\divH^{a+3}$, which is compactly embedded into $\divH^{a+2}$ by Rellich's theorem (e.g., p.330 in \cite{Taylor2011}).
\end{Proof}
\vspace{3mm}


\begin{Theor}\label{TheorInverseI}
    Fix $a\in\Z$ with $a\geq -1$, and let $b=\max\{a+4, |a|+1\}$. Assume that $\theta\in\divH^b$ and $f\in \divH^{b-1}$. Then $\Delta\I_\theta^* \I_\theta$ is a homeomorphism of $\divH^a$.
\end{Theor}

\begin{Proof}{Theorem \ref{TheorInverseI}}
    By the decomposition $\I_\theta^* \I_\theta=\LL^* \LL + \KK$ from (\ref{eq:IstarI}), and putting together Proposition \ref{PropLstarL} and Lemma \ref{PropCompactPerturbation} entails that $\I_\theta^* \I_\theta$ maps $\divH^a$ continuously into $\divH^{a+2}$. Consequently, $\Delta \I_\theta^* \I_\theta : \divH^a\to \divH^a$ is a linear and continuous mapping, with
    \begin{align*}
    \Delta \I_\theta^* \I_\theta 
    &=
    \Delta \LL^* \LL 
    +
    \Delta \KK 
    \\ 
    &=
    \frac{1}{2\nu}\big(S_{2T}-{\rm Id}_{\divH^a}\big) + \Delta \KK 
    .
    \end{align*}
    Since $\KK : \divH^a\to \divH^{a+2}$ is compact by Lemma \ref{PropCompactPerturbation}, and $\Delta : \divH^{a+2} \to \divH^a$ is continuous, we deduce that $\Delta \KK : \divH^a\to \divH^a$ is a compact operator. We also have that $S_{2T} : \divH^a\to \divH^a$ is compact by virtue of (\ref{eq:HeatSobolev}). Consequently, $\Delta \I_\theta^* \I_\theta$ writes as 
    $$
    \Delta \I_\theta^* \I_\theta 
    =
    -\frac{1}{2\nu} {\rm Id}_{\divH^a}
    +
    \tilde{\KK}
    ,
    $$
    where $\tilde{\KK}$ is a compact operator on $\divH^a$. Consequently, Fredholm's alternative (see Theorem 5 of Section D.5 in \cite{Evans1998}) implies that $\Delta \I_\theta^* \I_\theta$ is a homeomorphism of $\divH^a$, provided $\Delta \I_\theta^* \I_\theta$ is injective on $\divH^a$. To see this let $h\in\divH^a$ be such that $\Delta \I_\theta^* \I_\theta[h] = 0$. Since $\Delta$ is injective on $\divH^{a+2}$, then $\I_\theta^* \I_\theta[h]=0$ in $\divH^{a+2}$. Since $a\ge -1$, we have $\divH^{a+2} \subset \divH^{-a}$ so that we may write the $H^a$-$H^{-a}$ dual pairing as
    $$
    0 
    =
    \ps{\I_\theta^* \I_\theta[h]}{h}_{L^2}
    =
    \ps{\I_\theta[h]}{\I_{\theta}[h]}_{L^2([0,T],\Hnull)}
    =
    \|\I_\theta[h]\|^2_{L^2([0,T],\Hnull)}
    .
    $$
    Since $b\geq 3$ we have $\theta\in \divH^2$ and $f\in \divH^1$, so that the stability estimate 
    $$
    \|\I_\theta[\phi]\|_{L^2([0,T],\Hnull)}
    \gtrsim 
    \|\phi\|_{\dH^{-1}}
    ,\spa 
    \forall\ \phi\in\divH^{-1}
    ,
    $$
    from Proposition \ref{PropLinStab} with $a=-1$, entails that $h=0$ in $\divH^{-1}\supset \divH^{\alpha}$. We deduce that $\Delta \I_\theta^* \I_\theta : \divH^a \to \divH^a$ is injective, hence concludes the proof.
\end{Proof}

\subsection{The limit process and convergence in function space}

\subsubsection{The law $\mathcal N_{\theta_0}$} \label{limproc}

Given $h \in \Hnull \cap C^\infty$, let us denote by $\bar h$ the solution to the equation $\Delta \I^*_{\theta_0} \I_{\theta_0}\bar h = \Delta h$. In view of Theorem \ref{TheorInverseI} and since $\Delta$ is a homeomorphism from $\divH^{a+2} \to \divH^a$, this solution is unique and lies in $\divH^a$ whenever $h \in \divH^{a+2}$. In particular we can write $\bar h = (\I^*_{\theta_0} \I_{\theta_0})^{-1}h \in \divH^a$ and $(\I^*_{\theta_0} \I_{\theta_0})^{-1}$ maps $\divH^{a+2}$ homeomorphically into $\divH^a$ for any $a\in\Z$ with $a\geq -1$. With these preparations,  consider a centred Gaussian process $(\mathbb W(g): g \in \Hnull \cap C^\infty)$ with covariance 
$$
E\big[\W(g) \W(h)\big]
= 
\langle g, (\I_{\theta_0}^*\I_{\theta_0})^{-1}h\rangle_{L^2}
,\spa 
g,h \in \Hnull \cap C^\infty.
$$
Restricting $\W$ to the eigenfunctions $\{ e_j : j\geq 1\}$ of the Stokes operator $A=-\Delta$ defines a cylindrical probability measure on $\R^\N$ as the law $\NN_{\theta_0}$ of $Z=\{\W(e_j) : j\geq 1\}$. We further have
\begin{align} \label{expnorm}
\E[\|Z\|^2_{\dH^{-\beta}}]
&=
\E\Big[\sum_{j\geq 1}\lambda_j^{-\beta} \lvert\W(e_j)\rvert^2\Big] 
=
\sum_{j\geq 1} \lambda_j^{-\beta} \ps{e_j}{(\I_{\theta_0}^* \I_{\theta_0})^{-1} e_j}_{L^2}
\lesssim 
\sum_{j\geq 1} \lambda_j^{-\beta}\|e_j\|^2_{\dH^1}
,
\end{align}
where the last inequality follows from 
$$
\ps{e_j}{(\I_{\theta_0}^* \I_{\theta_0})^{-1}e_j}_{L^2}
\leq 
\|e_j\|_{\dH^1}\|(\I_{\theta_0}^* \I_{\theta_0})^{-1}e_j\|_{\dH^{-1}}
\lesssim 
\|e_j\|^2_{\dH^1}
,
$$
since $\I_{\theta_0}^* \I_{\theta_0} : \divH^{-1}\to \divH^1$ is a homeomorphism by virtue of Theorem \ref{TheorInverseI} with $a=-1$. Since $\|e_j\|_{\dH^1} = \lambda_j^{1/2}$ by (\ref{eq:LapNorm}), and $\lambda_j \asymp j$ by (\ref{eq:enum}), we deduce that
$$
\E[\|Z\|^2_{\dH^{-\beta}}]
\lesssim 
\sum_{j\geq 1} \lambda_j^{-(\beta-1)}
\lesssim 
\sum_{j\geq 1} j^{-(\beta-1)/2}
.
$$
This last series converges if and only $\beta>2$. Therefore $Z\equiv \sum_{j\geq 1} \W(e_j)e_j$ defines a $\divH^{-\beta}$-valued random variable for any such $\beta$, and the Oxtoby-Ulam theorem (see Proposition 2.1.4 in \cite{GinNic2016}) entails the sufficiency part  of the following result:
\begin{Prop}The law $\NN_{\theta_0} \equiv {\rm Law} (\sqrt T Z) = N\big(0, T(\mathbb I_{\theta_0}^* \mathbb I_{\theta_0})^{-1}\big)$ defines a tight Borel probability measure on the separable space $\divH^{-\beta}$ if and only if $\beta>2$.
\end{Prop}

Necessity of $\beta>2$ will be established below. Note also the necessary rescaling of the information operator by the observation time horizon $T$ in (\ref{infoop}), and that the dependence of $\mathbb I_{\theta_0}^* \mathbb I_{\theta_0}$ on $T$ is supressed in the notation.

The reproducing kernel Hilbert space (RKHS) of this process is equal to $\divH^{-1}$: Indeed, Proposition 2.6.8 in \cite{GinNic2016} with Banach space $\B=\divH^{-\beta}$ for some $\beta>2$, and $L^2$-dual $\B^*=\divH^\beta$ (see (\ref{dualSob2}) and (\ref{dualSob3})), entails that the RKHS of $\W$ is the $L^2$-dual space to the completion of $\Hnull\cap C^\infty$ with respect to the norm $\|f\|^2_*\equiv \ps{f}{(\I_{\theta_0}^* \I_{\theta_0})^{-1}f}_{L^2}$. Therefore, it is enough to establish that $\|\cdot\|_{*}$ and $\|\cdot\|_{\dH^1}$ define equivalent norms on $\divH\cap C^\infty$. For this purpose, observe, on the one hand, that Theorem \ref{TheorInverseI} with $a=-1$ entails that $\I_{\theta_0}^* \I_{\theta_0} : \divH^{-1}\to \divH^1$ is a homeomorphism, which yields 
$$
\ps{f}{(\I_{\theta_0}^* \I_{\theta_0})^{-1}f}_{L^2}
\leq 
\|f\|_{\dH^1}\|(\I_{\theta_0}^* \I_{\theta_0})^{-1}f\|_{\dH^{-1}}
\lesssim 
\|f\|^2_{\dH^1}
.
$$
On the other hand, Proposition \ref{PropLinStab} with $a=-1$ provides
$$
\ps{f}{(\I_\theta^* \I_\theta)^{-1}f}_{L^2}
=
\|\I_\theta (\I_{\theta_0}^* \I_{\theta_0})^{-1}f\|^2_{L^2([0,T],\Hnull)}
\gtrsim 
\|(\I_{\theta_0}^* \I_{\theta_0})^{-1}f\|^2_{\dH^{-1}}
\gtrsim 
\|f\|^2_{\dH^1}
,
$$
from which we deduce that the RKHS of $\W$ is $\divH^{-1}$. The last display further gives the converse inequality in (\ref{expnorm}) and hence, in view of Theorem 2.1.20 (and Example 2.1.6) in \cite{GinNic2016}, implies that any $X \sim \mathcal N_{\theta_0}$ satisfies $\Pr(\|X\|_{\dH^{-\beta}}<\infty)=0$ for all $\beta \le  2$, so that the hypothesis  $\beta>2$ in the preceding proposition is also necessary.



\subsubsection{Proof of Theorem \ref{nstokbvm}}

The first main result about asymptotic behaviour of posterior measures is the following limit theorem for the posterior initial conditions in a `weak' Schwartz-type topology of a sufficiently large negative Sobolev space. We did not attempt to optimise the value of $\sigma$ in the next theorem, mostly because this is not relevant for the proof of Theorem \ref{nstokbvm}, but also because a sharp bound $\sigma>2$, required for existence of the limiting process, can typically only be attained for information operators that are approximately diagonal in the basis used to expand the prior (see, e.g., \cite{CN13, CN14, NR20}), and finding such a basis  in the model studied here is a very challenging task.
   
\begin{Theor}\label{TheorBvMTheta} Assume $\theta_0, f\in  \divH^\infty, \nu>0$.
Let $\Pi$ be the prior defined in (\ref{prior}) with $\alpha>12$, and let $\theta|Z^{(N)}\sim \Pi(\cdot|Z^{(N)})$ be the posterior measure (\ref{post}) arising from data (\ref{model}) where $u_\theta$ solves the Navier-Stokes equations (\ref{eq:NS}), with corresponding posterior mean $\tilde{\theta}_N=\E[\theta|Z^{(N)}]$. Then, for any $\sigma>9$,
   $$
   \WWW_{1,\divH^{-\sigma}}\Big({\rm Law}\big(\sqrt{N}(\theta-\tilde{\theta}_N)\big), \NN_{\theta_0}\Big)
   \to 
   0\quad
   \text{in }\ P_{\theta_0}^{\N}\text{-probability}
   ,
   $$
    and, as $N \to \infty$,
   $$
   \sqrt{N}(\tilde{\theta}_N-\theta_0)
   \to^d
   \NN_{\theta_0}\quad
   \text{in }\ \divH^{-\sigma}.
   $$
\end{Theor}

 The proof of Theorem \ref{TheorBvMTheta} consists of an application of Theorem 2 in \cite{Nickl2024}, and is deferred to Section~\ref{bvmap} below. We now deduce Theorem \ref{nstokbvm} from Theorem \ref{TheorBvMTheta}. Recall the separable Banach space
$$\CCC=C([t_{\min}, t_{\max}], C(\Om)^2).$$ To ease notation, we will write $P_N$ and $E[\cdot]$  to denote, respectively, the posterior measure $\Pi(\cdot|Z^{(N)})$ and the expectation under the posterior measure $\E^\Pi[\cdot | Z^{(N)}]$. 

\begin{Proof}{Theorem \ref{nstokbvm}}
    From Theorem \ref{postcont} below we know that on events of $P_{\theta_0}^\N$-probability approaching $1$ as $N\to\infty$, the posterior mean $\tilde \theta_N=E[\theta]$ lies in $\divH^{\beta}$ for any $\beta\in [0,\alpha-1)$ and satisfies, for $\delta_N$ as in (\ref{RateDeltaN}) and all $\xi\in [0,\beta]$, $N\ge 1$ and $1\le q<\infty$
    \begin{equation}\label{PostMeanBound}
    \|\tilde\theta_N-\theta_0\|_{\dH^\xi}
    \leq 
    E\big[\|\theta-\theta_0\|^q_{\dH^\xi}]^{1/q}
    \lesssim
    \delta_N^{(\beta - \xi)/(\beta+1)}
    \equiv
    \tilde \delta_N(\xi)
    .
    \end{equation}
    In the rest of this proof, we will thus tacitly restrict to such frequentist events. Then, fix $\beta\in (1,\alpha-1)$ (to be determined later), and define the random variables
    $$
    X_N 
    =
    \sqrt{N}(\theta-\tilde{\theta}_N)|Z^{(N)}
    ,\spa 
    Y_N 
    =
    \sqrt{N}(u_{\theta}-u_{\tilde{\theta}_N}) | Z^{(N)}
    $$
   where $\theta\sim P_N$ is a posterior draw (conditional on $Z^{(N)}$). Since $\theta\in\divH^\beta$ with $P_N$-probability $1$ we know that~$X_N$ is an $\divH^\beta$-valued random variable, and Proposition \ref{PropNSBound} entails that $Y_N$ defines a random element of $C([0,T],\divH^\beta)$. Since $\beta>1$, the Sobolev embedding $\divH^\beta\embed H^\beta\embed C(\Om)^2$ further implies that $Y_N$ defines a $\CCC$-valued random variable.

    Propositions \ref{PropLinSmooth} and \ref{PropLinSob} with starting time $t_{\min}$ instead of $0$  imply that
    \begin{equation}\label{map}
    \I_{\theta_0}
    :
    \divH^{-a} 
    \to 
    \CCC
    ,\spa 
    \forall\ a\geq 0
    ,
    \end{equation}
    defines a continuous linear mapping. In particular if $X\sim \NN_{\theta_0}$ is the Gaussian random variable on $\divH^{-a}, a>2,$ constructed in Section \ref{limproc}, then $\I_{\theta_0}[X]$ and $\I_{\theta_0}[X_N]$ define random variables in $\CCC$, whose law coincides with that of the solution of the linear parabolic PDE (\ref{linshow}) with random initial condition $X$ and $X_N$, respectively. Proposition 2.1.4 in \cite{GinNic2016} further implies that  the law $\mu$ of $\I_{\theta_0}[X]$ is a tight Gaussian Borel probability measure on the separable Banach space $\CCC$. Therefore, recalling (\ref{wassdist}) and denoting by $\rm Lip$ the collection of $1$-Lipschitz maps $h:\CCC\to\R$, we have (with all expectations conditional on $Z^{(N)}$)	
\begin{eqnarray*}
    \lefteqn{
        \hspace{4mm}
        \WWW_{1,\CCC}({\rm Law} (Y_N), \mu) 
        =
        \WWW_{1,\CCC}\big({\rm Law}(Y_N), {\rm Law}(\I_{\theta_0}[X])\big)
        =
        \sup_{h\in {\rm Lip}} \big| E h(Y_N) - Eh(\I_{\theta_0}[X]) \big|
    }
    \\
    &&
    \hspace{12mm}
    =
    \sup_{h\in {\rm Lip}} 
    \Big| 
    E\Big[ h\big( \sqrt{N}(u_{\theta}-u_{\tilde{\theta}_N}) - 
    \I_{\theta_0}[X_N] + \I_{\theta_0}[X_N] \big)
    - h\big(\I_{\theta_0}[X_N]\big) \Big] 
    \\
    &&
    \hspace{32mm}
    +\
    E\big[ h\big(\I_{\theta_0}[X_N]\big) \big] - E\big[h\big(\I_{\theta_0}[X]\big)\big]
    \Big| 
    \\
    &&
    \hspace{12mm} 
    \leq 
    \sqrt{N}\ 
    E\Big[ \big\|u_{\theta}-u_{\tilde{\theta}_N}-\I_{\theta_0}[\theta-\tilde{\theta}_N]\big\|_{\CCC} \Big]
    +
    \sup_{h\in{\rm Lip}} 
    \big|
    E\big[ h\big(\I_{\theta_0}[X_N]\big) \big] - E\big[h\big(\I_{\theta_0}[X]\big)\big]
    \big|
    \\[2mm]
    &&
    \hspace{12mm}
    =
    I + II
    .
\end{eqnarray*}	
Let $\sigma>0$ be as in Theorem \ref{TheorBvMTheta}. Since $\I_{\theta_0} : \divH^{-\sigma}\to \CCC$ is linear and continuous it is  Lipschitz and then $h\circ \I_{\theta_0}:\divH^{-\sigma}\to \R$ is Lipschitz as well for all $h\in{\rm Lip}$ with Lipschitz-modulus upper-bounded by that of $\I_{\theta_0}$. Consequently, taking the supremum over all $1$-Lipschitz maps $g:\divH^{-\sigma}\to \R$, we have as $N\to\infty$
$$
II
\lesssim 
\sup_{g} 
\big|
E[ g(X_N) ] - E[g(X)]
\big|
=
\WWW_{1,\divH^{-\sigma}}\Big({\rm Law}\big( \sqrt{N}(\theta-\tilde{\theta}_N)\big), \NN_{\theta_0}\Big)
=
o_{P_{\theta_0}^\N}(1)
$$
in view of Theorem \ref{TheorBvMTheta}.

\smallskip

To show that $I$ converges to $0$ in $P_{\theta_0}^N$-probability, observe that
$$
u_{\theta}(t)-u_{\tilde{\theta}_N}(t)-\I_{\theta_0}[\theta-\tilde{\theta}_N](t)
=
R_{\theta}(t)-R_{\tilde \theta_N}(t)
,
$$
where we let $R_{\theta'}(t)\equiv u_{\theta'}(t)-u_{\theta_0}(t)-\I_{\theta_0}[\theta'-\theta_0](t)$ for any $\theta'$ and $t$. Since $\beta>1$, then for any $L>0$, a periodic version of Agmon's inequality (see Lemma 4.9 and Lemma 4.10 in \cite{CF88}), together with Proposition \ref{PropExistLin}(iii) yields that
$$
\|R_{\theta'}(t)\|_{L^\infty}
\lesssim 
\|R_{\theta'}(t)\|^{(\beta-1)/\beta}_{L^2}\|R_{\theta'}(t)\|^{1/\beta}_{\dH^\beta}
\lesssim 
\|\theta'-\theta_0\|^{2(\beta-1)/\beta}_{L^2}\|R_{\theta'}(t)\|^{1/\beta}_{\dH^\beta}
,
$$
holds uniformly in $t\in [t_{\min},t_{\max}]$ and $\theta'\in\divH^\beta$ with $\|\theta'\|_{\dH^\beta}\le L$.
We further have by Proposition \ref{PropNSBound} with $a=\beta$ and Proposition \ref{PropLinSob} with $a=\beta$
\begin{eqnarray} \label{betbd}
\lefteqn{
\hspace{-10mm}
\sup_{t\in [t_{\min}, t_{\max}]}
\|R_{\theta'}(t)\|_{\dH^\beta}
\leq 
\sup_{t\in [t_{\min}, t_{\max}]}
\big( \|u_{\theta'}(t)\|_{\dH^\beta} + \|u_{\theta_0}(t)\|_{\dH^\beta} + \|\I_{\theta_0}[\theta'-\theta_0](t)\|_{\dH^\beta}
\big)
} \notag
\\[2mm]
&&\lesssim 
1+\|\theta'\|^{2\beta}_{\dH^\beta} + \|\theta_0\|^{2\beta}_{\dH^\beta} + \|\theta'-\theta_0\|_{\dH^\beta}
\lesssim 
1+\|\theta'\|^{2\beta}_{\dH^\beta}
.
\end{eqnarray}
With these preparations, fix $M>0$ and define 
$$
S_M
=
\big\{\theta'\in \divH^{\beta} : \|\theta'\|_{\dH^\beta}\leq \|\theta_0\|_{\dH^\beta} + M \big\}
.
$$
Since $\tilde\theta_N\in S_M$ for large enough $M$ by virtue of (\ref{PostMeanBound}) with $\xi=\beta$, we deduce from what precedes that 
\begin{equation*}
\|u_{\theta}-u_{\tilde{\theta}_N}-\I_{\theta_0}[\theta-\tilde{\theta}_N]\|_{\CCC} 1_{\{\theta\in S_M\}}
\lesssim 
\|\theta-\theta_0\|^{2(\beta-1)/\beta}_{L^2}
+
{\tilde\delta_N(0)}^{2(\beta-1)/\beta},
\end{equation*}
where we used again (\ref{PostMeanBound}) with $\xi=0$ in the last inequality. We then have, again by (\ref{PostMeanBound}),
\begin{align} \label{1stbd}
E\Big[\|u_{\theta}-u_{\tilde{\theta}_N}-\I_{\theta_0}[\theta-\tilde{\theta}_N]\|_{\CCC} 1_{\{\theta\in S_M\}}\Big]
\lesssim 
E\big[\|\theta-\theta_0\|^{2(\beta-1)/\beta}_{L^2}\big]
+
\tilde \delta_N(0)^{2(\beta-1)/\beta}
\lesssim 
\delta_N^m 
,
\end{align}
where $m=\frac{2(\beta-1)\beta}{\beta(\beta+1)}=2(\beta-1)/(\beta+1)$ (all with $P_{\theta_0}^\mathbb N$- probability approaching one). We obtain
$$
\delta_N^m 
=
N^{-\frac{\alpha}{2\alpha+2} \frac{2(\beta-1)}{\beta+1}}
=
o\Big(\frac{1}{\sqrt{N}}\Big)
,
$$
provided one can pick $\beta\in (1,\alpha-1)$ such that 
$$
\frac{\alpha}{2\alpha+2} \frac{2(\beta-1)}{\beta+1}
>
\frac12 
.
$$
This happens as soon as $(3\alpha+2)/(\alpha+2)<\beta<\alpha-1$ which is possible since $\alpha>1+\sqrt{5}$. Consequently, for such a $\beta$ we deduce that the bound in (\ref{1stbd}) is $o_{P_{\theta_0}^\N}(1)$, as required.

Now let us deal with the term involving $S_M^c.$ Using (\ref{betbd}), the Sobolev embedding 
$\divH^{\beta}\embed H^\beta \embed C(\Om)^2$  we have
\begin{eqnarray*}
\lefteqn{
\|u_{\theta}-u_{\tilde{\theta}_N}-\I_{\theta_0}[\theta-\tilde{\theta}_N]\|_{\CCC}
\leq 
\|R_\theta\|_\CCC + \|R_{\tilde\theta_N}\|_\CCC 
}
\\[2mm]
&&
\lesssim 
\sup_{t\in [t_{\min}, t_{\max}]}\Big(\|R_\theta(t)\|_{\dH^\beta} + \|R_{\tilde\theta_N}(t)\|_{\dH^\beta} \Big)
\lesssim 
1+\|\theta\|^{2\beta}_{\dH^\beta} + \|\tilde \theta_N\|^{2\beta}_{\dH^\beta}
.
\end{eqnarray*}
It follows from the Cauchy-Schwarz inequality that 
$$
E\Big[\|u_{\theta}-u_{\tilde{\theta}_N}-\I_{\theta_0}[\theta-\tilde{\theta}_N]\|_{\CCC} 1_{\{\theta\in S_M^c\}}\Big]
\lesssim 
\Big(1+\|\tilde\theta_N\|^{2\beta}_{\dH^\beta}
+E\Big[\|\theta\|^{2\beta}_{\dH^\beta}\Big]\Big)^{1/2} P_N[S_M^c]^{1/2}
.
$$
But Theorem \ref{postcont} entails that, for some $b>0$, we have
\begin{equation*}
    P_{\theta_0}^N\Big( P_N\big(S_M^c\big)
    \le 
    e^{-bN\delta_N^2} 
    \Big)
    \to 
    1
    , 
\end{equation*}
so that, up to further restricting ourselves to events of $P_{\theta_0}^N$-probability approaching $1$ as $N\to\infty$, we may assume that 
\begin{equation}\label{pEMcomp}
P_N(S_M^c\big)
\le 
e^{-bN\delta_N^2} 
,\spa 
\forall\ N\geq 1
.
\end{equation}
Combining (\ref{pEMcomp}), (\ref{PostMeanBound}) with $\xi=\beta$ for $\tilde\theta_N$, and (\ref{PostMeanBound}) with $\xi=\beta$ and $q=2\beta$ for $\theta$, yields 
$$
\sqrt{N}\ E\Big[\|u_{\theta}-u_{\tilde{\theta}_N}-\I_{\theta_0}[\theta-\tilde{\theta}_N]\|_{\CCC} 1_{\{\theta\in S_M^c\}}\Big]
=
O_{P_{\theta_0}^N}\big(\sqrt{N} e^{-\frac{b}{2}N\delta_N^2}\big)
=
o_{P_{\theta_0}^N}(1)
,
$$
which concludes the proof.
\end{Proof}

\begin{Remark} \normalfont \label{wow}
The proof shows that the limit theorem holds in $C([t_{\min}, t_{\max}], \divH^\zeta)$ for some $\zeta>1$ (one just uses a slightly stronger version of the forward smoothing estimate in (\ref{map}) and later uses a Sobolev imbedding in place of Agmon's inequality).
\end{Remark}

\subsection{Delta method and exact asymptotic minimax optimality}\label{functest}

We established in Theorem~\ref{nstokbvm} that, for $\tilde\theta_N=\E^\Pi[\theta|Z^{(N)}]$ the posterior mean of $\Pi(\cdot | Z^{(N)})$ arising from data $Z^{(N)}$ in (\ref{model}) and prior $\Pi$ in (\ref{prior}) with $\alpha>12$, then the random process $\sqrt{N}(u_{\tilde \theta_N}-u_{\theta_0})$ converges weakly in $\CCC$ to the law $\mu$ of the solution $U$ of the linear parabolic equation (\ref{linshow}) with Gaussian random initial condition $\vartheta\sim \NN_{\theta_0}$. In particular, the functional Delta-method (see, e.g., Section 3.10 in \cite{WelVan1996}) entails that 
$$
\sqrt N (\Phi(u_{\tilde \theta_N}) - \Phi(u_{\theta_0})) 
\stackrel{d}{\to}
\dot\Phi_{u_{\theta_0}}(U)
\quad \textrm{in}\quad \B,\quad 
N \to \infty
,
$$ 
for any Hadamard-differentiable map $\Phi : \CCC \to \B$ taking values in the normed space $\B$, with Hadamard derivative $\dot \Phi$; that is, if $\Phi: \CCC \to \B$ is such that for all $v\in\CCC$ there exists a bounded linear operator $\dot \Phi_v : \CCC\to \B$ such that 
$$
\frac{\Phi(v+t_n h_n)-\Phi(v)}{t_n}
\to 
\dot \Phi_v[h]
\spa 
\textrm{in}\quad \B 
,
$$
for all $h\in\CCC$, and sequences $(t_n)\subset \R$ and $(h_n)\subset \CCC$ such that $t_n\to 0$ and $h_n\to h$ in $\CCC$. For instance, one shows without difficulty (using estimates for the `nonlinear term' in the appendix) that the `advection' map $\Phi_t(u) = (u(t) \cdot \nabla) u(t), t>0,$ is a Fr\'echet- (and hence Hadamard) differentiable map from $\CCC$ into $\divH^{-1}$ with derivative at $u$ given by $$\dot \Phi_{t,u}[h] = (u(t) \cdot \nabla)h(t) + (h(t) \cdot \nabla)u(t).$$ The resulting functional limit theorem for $$\sqrt N \big((u_{\tilde \theta_N}(t) \cdot \nabla) u_{\tilde \theta_N}(t) - (u_{\theta_0}(t) \cdot \nabla) u_{\theta_0}(t)\big) \to \dot \Phi_{t,u_{\theta_0}}(U(t))$$ can be strengthened to hold in $\divH$ rather than just in $\divH^{-1}$ by appealing to Remark \ref{wow}. 

\smallskip

We thus obtain $\sqrt N$-consistent estimators with mean zero Gaussian limits in $\B$ for a large class of functionals $\Phi$, and just as in the classical Gauss-Markov, Cram\'er-Rao, or local asymptotic minimax theorems, one can ask whether the limiting covariance is \textit{minimal} in a certain sense, including the case of state prediction itself where $\Phi \equiv {\rm Id}$ is the identity on $\CCC$. A result of this type entails that Bayesian data assimilation methods, which employ the initial condition $\theta$ as the underlying parameter space, do in fact (asymptotically) extract \textit{all} relevant information to estimate/predict $u_\theta, \Phi(u_\theta)$ from the sample $Z^{(N)}$, and therefore provides a benchmark for any other algorithm attempting to recover these parameters. Just as in Section 7.5 of \cite{Nickl2020} in a slightly easier measurement and PDE setting, the analytical techniques in the present paper and the infinite-dimensional minimax theorem (e.g., Section 3.12 of \cite{WelVan1996}) entail a general minimax lower bound for estimating functionals $F(\theta)=\Phi(u_\theta)$, and this is investigated in full detail in the companion paper \cite{Kon25Minimax}, covering two concrete applications in the setting of Navier-Stokes equations: $F(\theta)=u_\theta$ itself as well as the advection term $F(\theta)=(u_{\theta}\cdot \nabla) u_\theta$, both at strictly positive times $t>0$. Let us state, for further reference, the main result relevant for us from \cite{Kon25Minimax}. It establishes a lower bound for the local asymptotic minimax risk of \textit{any} estimator $T_N$, which is attained by the estimator $\Phi(u_{\tilde \theta_N})$ provided in Theorem \ref{nstokbvm}.


\begin{Theor}\label{TheorMinimaxU}
    Fix $b\geq 2$. Let $\theta_0\in\divH^b$ and $f\in\divH^{b-1}$. Let $\B$ and $\tilde\B$ be Banach spaces, and $\Phi:\tilde\B\to\B$ be a Hadamard differentiable map with Hadamard derivative $\dot\Phi$. Fix $0<t_{\min}<t_{\max}<\infty$, and further assume that
    the continuous embedding 
    $$
    C([t_{\min},t_{\max}], \divH^{b})
    \embed 
    \tilde\B 
    $$
    holds (including the choice $\tilde \B = \CCC$). Consider data $Z^{(N)}=(t_i, \omega_i, Y_i)_{i=1}^N$ as in (\ref{model}) with $u_\theta$ solving the periodic $2D$ Navier-Stokes equations (\ref{eq:NS}), and a random field $U \sim \mu $ drawn as the solution of the linear PDE (\ref{linshow}) with Gaussian random initial condition $\vartheta\sim\NN_{\theta_0}$. 
    Then, we have 
\begin{equation}
\inf_{(\tilde u_N)} 
\sup_{\substack{J\subset \divH^\infty\\ |J|<\infty}}
\liminf_{N \to \infty} 
\max_{h\in J}
N E_{\theta_0+\frac{h}{\sqrt N}}^N
\|\tilde u_N-\Phi(u_{\theta_0+h})\|^2_{\B}  
\ge 
E_{\mu}\|\dot\Phi_{u_{\theta_0}}\big[U]\big]\|^2_{\B}
\end{equation}
where the infimum ranges over all sequences $(\tilde u_N:([0,T]\times \Om \times \R^2)^N\to \B)$ of measurable functions of the observations~$Z^{(N)}$, and the supremum ranges over all subsets $J\subset \divH^\infty$ with finite cardinality $|J|$.
\end{Theor}

\subsection{Proof of Theorem \ref{TheorBvMTheta}} \label{bvmap}

Our proof is based on verifying the conditions of the general Bernstein-von Mises theorem (Theorem 2) in \cite{Nickl2024} with $d=2$, $\Om=\T^2$ and $W=\R^2$. This requires some minor notational adjustments to our `divergence-free vector field' setting. In essence, the proof in \cite{Nickl2024} only requires an appropriate Sobolev scale with basis functions arising as eigenfunctions of a Laplace type operator with corresponding spectral asymptotics. For the convenience of the reader we restate here all the relevant results in \cite{Nickl2024} and provide some remarks about modifications in the proofs when necessary.

\smallskip


The prior we consider---see (\ref{prior}) and (\ref{PriorSeries})---is obtained as the law $\Pi = \Pi_{\alpha, N}$ of 
\begin{equation}\label{eq:Prior}
\theta 
= 
\frac{1}{\sqrt N \delta_N} \theta' \text{ where }\delta_N 
= 
N^{-\frac{\alpha}{2\alpha+2}},
\end{equation}
for $\theta'\sim \Pi'$, where $\Pi'$ is the centred Gaussian Borel probability measure on $\divH$ with RKHS $\HH\equiv \divH^\alpha$ for some $\alpha>2$. The RKHS $\HH_N$ of $\Pi_{\alpha,N}$ is also given by $\divH^\alpha$ but with RKHS norm 
$$
\|\cdot\|_{\HH_N} 
= 
\sqrt{N} \delta_N \|\cdot\|_\HH
,
$$
in view of the rescaling in (\ref{eq:Prior}). Then by the Karhunen-Loève theorem,  $\Pi$ is the law of the Gaussian random series of vector fields
\begin{equation}
\theta(x_1, x_2) 
= 
\sum_{j  \ge 1} g_j \lambda_j^{-\alpha/2} f_j(x_1, x_2),\quad (x_1,x_2) \in \Omega,\quad g_j \stackrel{iid}{\sim} N(0,1),
\end{equation}
where $(f_j)$ is the orthonormal basis of $\divH$ from (\ref{eq:BasisDivH}). The previous series converges almost surely in $\divH^\beta$ for any $\beta<\alpha-1$ (see Section \ref{sec:Intro} for details), hence $\Pi_{\alpha,N}$ is supported on the separable Hilbert space $\divH^\beta$ for any such $\beta$. Note, however, that this construction does technically not allow for a direct use of Theorem 2 in \cite{Nickl2024} as our $L^2(\Om, \R^2)$-valued prior is supported on a closed (strict) subspace of $L^2(\Om,\R^2)$ which cannot be expressed as coordinate-wise copies of a random series prior expanded in a basis of $L^2(\Om, \R)$. Nonetheless, inspection of Theorem 2 in \cite{Nickl2024} shows that their proof applies to our setting by replacing the (real-valued) orthonormal basis $(h_j)$ in eq.~(23) there by our (vector-field) basis $(f_j)$, noting also that the behaviour of their spectral coefficients $(\lambda_j)$ is the same as ours by virtue of (\ref{eq:BasisDivH}).\smallskip


As in (\ref{post}), the posterior measure $\theta|Z^{(N)}$ arising from data $Z^{(N)}$ via (\ref{model}) is then given by
\begin{equation}\label{postlog}
d\Pi(\theta|Z^{(N)}) \propto e^{\ell_N(\theta)} d\Pi(\theta),~~\ell_N(\theta) = -\frac{1}{2} \sum_{i=1}^N |Y_i - \GGG(\theta)(X_i)|_W^2,~~\theta \in \divH
,
\end{equation}
where we write $X_i = (t_i, \omega_i)$, and where the forward map 
\begin{equation}\label{gmap}
\GGG : \divH^1 \to L^2([0,T], \divH),\ \theta\mapsto u_{\theta}
,
\end{equation}
is given by the solution $u_\theta$ to the Navier-Stokes system (\ref{eq:NSequations}). The main result is a Bernstein-von Mises theorem for the posterior on $\theta$ in the `weak' norm topology of the space $\divH^{-k-2}$ with $k$ as in (\ref{eq:constants}). The posterior mean vector-field $\tilde \theta_N= E^\Pi[\theta|Z^{(N)}]$ exists as a Bochner integral in $\divH$ and hence also in $\divH^{-k-2}$, while the limiting Gaussian measure~$\mathcal N_{\theta_0}$ is constructed in Section \ref{limproc} as a tight Borel law on $\divH^{-k-2}$.


\begin{Theor}\label{TheorBvMTheta2} Assume $\theta_0, f\in  \divH^\infty$ and $\nu>0$, and let $\alpha>12, k>7$.
Let $\Pi$ be the prior defined in (\ref{prior}), and let $\theta|Z^{(N)}\sim \Pi(\cdot|Z^{(N)})$ be the posterior measure (\ref{post}) arising from data (\ref{model}) where $u_\theta$ solves the Navier-Stokes equations (\ref{eq:NS}), with corresponding posterior mean $\tilde{\theta}_N=\E[\theta|Z^{(N)}]$. Then, we have as $N\to\infty$
   $$
   \WWW_{1,\divH^{-k-2}}\Big({\rm Law}\big(\sqrt{N}(\theta-\tilde{\theta}_N)\big), \NN_{\theta_0}\Big)
   \to 
   0\quad
   \text{in }\ P_{\theta_0}^{\N}\text{-probability}
   ,
   $$
    and, under $P_{\theta_0}^\N$,
   $$
   \sqrt{N}(\tilde{\theta}_N-\theta_0)
   \to^d 
   \NN_{\theta_0}\quad
   \text{in }\ \divH^{-k-2}.
   $$
\end{Theor}

We establish Theorem \ref{TheorBvMTheta2}, which implies Theorem \ref{TheorBvMTheta} with $\sigma=k+2$, by following the same steps used in \cite{Nickl2024} to prove Theorem~2 there. This proof is designed for general regression models and requires the verification of analytical conditions on the forward map $\theta\mapsto \GGG(\theta)\equiv u_{\theta}$ from (\ref{model}), which we now recall and then verify in our setting. In doing so, a number of constants need to be chosen.
 In view of the preceeding remarks, we replace the balls $U(r,L)$ in (25) from \cite{Nickl2024} with the corresponding balls
\begin{equation} \label{rball}
U(r, L)
=
\big\{u \in \divH^r: \|u\|_{\dH^r} \le L \big\}
,\spa L>0
,
\end{equation}
of mean-zero and divergence-free periodic vector fields. We take integers $\beta, \zeta\in\N$, and real numbers $\alpha$, $\kappa$ and $k$ satisfying 
\begin{equation}\label{eq:constants}
\zeta = 2,\spa 
\beta=11,\spa
\alpha > 12,\spa 
6 < \kappa < k-1
.
\end{equation}

\begin{Condition}\label{gemol}
Fix constants $\alpha,\beta, \zeta$ satisfying (\ref{eq:constants}), and let $\theta_0 \in \divH^\alpha$. Then, the forward map $\GGG$ from (\ref{gmap}) satisfies the following:

\smallskip

A) For every $L>0$, there exists $c=c(L)>0$ such that
\begin{equation} \label{ubd}
\sup_{\substack{\theta \in U(\beta, L)\\0 < t \le T\\ x \in \Omega}} |\GGG(\theta)(t,x)| \le c <\infty.
\end{equation}

\smallskip

B) For every $L>0$, there exists $c=c(L)>0$ such that  
$$
\|\GGG(\theta)-\GGG(\vartheta)\|_{L^2([0,T], L^2(\Omega))}
\le 
c \|\theta - \vartheta\|_{L^2(\Omega)}
,\spa 
\forall\ \theta, \vartheta \in U(1, L)
.
$$ 

\smallskip

C) There exists a continuous linear operator 
\begin{equation}\label{linmap}
\mathbb I_{\theta_0} \equiv D\GGG_{\theta_0} : \divH \to L^2([0,T], \divH)
,
\end{equation}
such that, for every $B>0$ and $L>0$, there exists $c=c(\theta_0, B,L)>0$ such that 
$$\|\GGG(\theta_0+h) - \GGG(\theta_0) - \mathbb I_{\theta_0}[h]\|_{L^2([0,T], L^2(\Omega))} \le c \|h\|^2_{L^2(\Omega)}
$$
holds for all $h \in U(\beta, L)$ with $\|h\|_{L^2(\Om)} \le B$.

\smallskip

D) For every $L>0$, there exists $c=c(\theta_0, L)>0$ such that
\begin{equation}\label{zetabus}
\|\GGG(\theta) - \GGG(\theta_0)\|_{L^\infty([0,T]), L^\infty(\Omega))} \le c \|\theta - \theta_0\|_{\dH^{\zeta}},~~\forall\ \theta \in U(\zeta, L) 
.
\end{equation}
In addition, the linearization $\mathbb I_{\theta_0}$ from (\ref{linmap}) is continuous from $\divH^{\zeta} \to L^\infty([0,T],L^\infty(\Omega))$.

\smallskip

E) For every $L>0$, there exists $c=c(\theta_0, L)>0$ such that 
$$
\|\GGG(\theta) - \GGG(\theta_0)\|_{L^2([0,T], L^2(\Omega))} 
\ge 
c\|\theta-\theta_0\|_{\dH^{-1}},\spa 
\forall\ \theta \in U(\beta, L)
.
$$

\smallskip

F) Consider the adjoint $\I_{\theta_0}^* : L^2([0,T], \divH) \to \divH$ and the information operator $\I_{\theta_0}^*\I_{\theta_0}$ acting on $\divH^\eta \subset\divH$ by restriction, with $\eta\ge 0$. For $\Delta$ the Laplacian and all $\eta \ge 0$, the operator
\begin{equation}\label{ical}
\mathcal I \equiv \Delta \mathbb I_{\theta_0}^*\mathbb I_{\theta_0}
\end{equation}
defines a continuous linear homeomorphism from $\divH^\eta \to \divH^\eta$.
\end{Condition}

\begin{proof}{[Verification of Condition \ref{gemol} for the Navier-Stokes system:]}
We start by verifying the basic forward regularity estimates: Condition $A$ follows from the Sobolev embedding $\divH^2\embed C(\Om)^2$ and Proposition \ref{PropNSBound}(i) with $a=\beta$. Condition B follows from $L^2$-Lipschitz continuity of the forward map; see, e.g., Proposition 1B in \cite{NicTit2024} or Proposition \ref{PropNSExist}(ii). 

\smallskip

Now about the linearization: The existence part of Condition C follows from Proposition \ref{PropExistLin}(i) with $\I_{\theta_0}[\xi](t,\cdot)\equiv\Lambda_{\theta_0,t}[\xi](\cdot)$, see~also (\ref{defI}). The fact that $\I_{\theta_0}$ maps $\Hnull$ into $L^2([0,T],\Hnull)$ follows from Proposition \ref{PropLinSob} with $a=0$, and the quadratic approximation condition follows from Proposition \ref{PropExistLin}(iii).

\smallskip

Condition D follows from the Sobolev embedding $\divH^2\embed C(\Omega)^2$ and the Sobolev continuity of the forward map and its linearization, established in Proposition \ref{PropNSLip}
and Proposition \ref{PropLinSob}, respectively, with $a=\zeta(\ge 2)$ in both cases.

\smallskip

We finally turn to the crucial quantitative injectivity estimates: Condition E follows from the stability estimate in  Proposition~\ref{PropStabForward}, while Condition F holds for any $\eta_0\geq -1$ on the closed linear subspaces $\divH^\eta$ of $H^\eta_0(\Om)^2$, by virtue of Theorem \ref{TheorInverseI}. Let us remark that, formally, Condition F in \cite{Nickl2024} is required for $\tilde {\mathbb I}$ defined on the larger space $L^2$, whose restriction to $\Hnull$ coincides with $\mathbb I_{\theta_0}$. If $P$ is the Leray projector and $\tilde {\mathbb I}^*$ the adjoint of $\tilde {\mathbb I}$, then the information operator restricted to $\divH$ satisfies $P\tilde {\mathbb I}^*\tilde {\mathbb I}_{|\divH} = \mathbb I_{\theta_0}^* \mathbb I_{\theta_0}$ so that Theorem 2 in \cite{Nickl2024} still applies.
\end{proof}








To proceed define the span in $\divH$ of the trigonometric basis of $\divH$ from (\ref{eq:BasisDivH}) as
\begin{equation}\label{ej}
E_J = {\rm span}\{f_j: 1 \le j \le J\},\spa 
\forall\ J\ge 1,
\end{equation}
with associated $L^2$-projector $P_{E_J}: \divH \to E_J$. For $0\le\xi\le\beta$, define sequences
\begin{equation}\label{deltan}
\tilde \delta_N(\xi) = \delta_N^{(\beta-\xi)/(\beta+1)}
,\spa 
\tilde \delta_N(0) \equiv \tilde \delta_N,
\end{equation}
and
\begin{equation} \label{approxrkhs}
J_N \simeq N \delta_N^2,\spa 
J_N\in\N 
.
\end{equation}
For $\psi \in \divH$, define $L^2$-projections onto $E_{J_N}$ from (\ref{ej}) as
\begin{equation}\label{projpsi}
p_N(\psi)
\equiv 
P_{E_{J_N}}[\psi]
.
\end{equation}
The consistency result below is expressed in terms of measurable sets defined, for any $L>0$, by
\begin{align} \label{barthetan}
\bar \Theta_N &=\Bar \Theta_{N,L}
\equiv 
\Big\{
\theta \in \divH^1 : 
\|\theta\|_{\dH^{\beta}} \le L,~ 
\sup_{\psi \in U(\kappa,1)} |\langle \theta, p_N(\psi) \rangle_{\mathcal H_N}| \le L \sqrt N \delta_N K_N  
\Big\}  \\
&~~~~ \cap ~\Big\{
\|\GGG(\theta)-\GGG(\theta_0)\|_{L^2([0,T], L^2(\Omega))} \le L \delta_N,~ 
\|\theta - \theta_0\|_{\dH^\xi} \le L \tilde \delta_N(\xi)\ ~\forall\ 0 \le \xi <\beta 
\Big\}, \notag
\end{align}
and the proof of the following theorem is the same as that of Theorem 3 in \cite{Nickl2024}.

\begin{Theor} \label{postcont} 
Fix $\eta>0$ arbitrarily small and let 
$$
K_N = \sqrt N \delta_N\max(1,J_N^{\frac{\alpha- (\kappa-1)+\eta}{2}})
.
$$
Then, for all $b>0$ there exists $L>0$ such that
\begin{equation} \label{ratekey}
P_{\theta_0}^\mathbb N \big(\Pi(\bar \Theta_{N,L} |Z^{(N)}) \le 1 - e^{-bN\delta_N^2} \big) \to_{N \to \infty} 0.
\end{equation}
For any $1 \le q<\infty$ and $0 \le \xi \le \beta$, we further have 
\begin{align}\label{meanreg}
\|E^\Pi[\theta|Z^{(N)}] - \theta_0\|_{\dH^\xi} &\le (E^\Pi[\|\theta-\theta_0\|^q_{\dH^\xi}|Z^{(N)}])^{1/q} = O_{P_{\theta_0}^\mathbb N}(\tilde \delta_N(\xi))
\end{align}
\end{Theor}

The previous result allows to show that the posterior measure $\Pi(\cdot|Z^{(N)})$ is asymptotically equivalent to a hypothetical posterior $\Pi^{\bar \Theta_N}(\cdot|Z^{(N)})$ arising from (\ref{postlog}) with `localised' prior 
\begin{equation} \label{respr}
\Pi^{\bar \Theta_N} = \frac{\Pi(\cdot \cap \bar \Theta_N)}{\Pi(\bar \Theta_N)}
\end{equation}
 restricted to the regularisation set $\bar \Theta_N=\bar \Theta_{N,L}$ from ({\ref{barthetan}). 
If we denote by $\tau_N$ the conditional law
$$
\tau_N 
= 
{\rm Law}(\sqrt N (\theta - T))
$$ 
for $\theta \sim \Pi(\cdot|Z^{(N)})$ and any fixed re-centring $T \in \divH^{-k-2}$, and if we denote by $\bar \tau_N$ the corresponding law where $\theta$ is replaced by a draw $\theta \sim \Pi^{\bar \Theta_N}(\cdot|Z^{(N)})$, then the proof of Proposition~1 in \cite{Nickl2024} yields the following approximation in the Wasserstein distance $\WWW_{1, \divH^{-k-2}}$ featuring in Theorem \ref{TheorBvMTheta2}.

\begin{Prop} \label{wassconv}
For any $T \in \divH^{-k-2}$ and some $c>0$, we have as $N \to \infty$
$$\WWW_{1, \divH^{-k-2}}(\tau_N, \bar \tau_N) 
= 
O_{P_{\theta_0}^\N}(e^{-cN\delta_N^2})
.
$$
\end{Prop}

Consequently, it suffices to prove the first limit in Theorem \ref{TheorBvMTheta} for $\bar \tau_N$ rather than $\tau_N$. The second and most important step consists of an asymptotic approximation of the Laplace transform of $\Pi^{\bar \Theta_N}(\cdot|Z^{(N)})$. For any fixed $\psi \in \divH^{\kappa+2}$ with $\kappa$ as in (\ref{eq:constants}), let 
$$
\bar \psi
=
\bar \psi_{\theta_0} 
= 
(\I_{\theta_0}^*\I_{\theta_0})^{-1}[\psi]
,
$$ 
which, by virtue of Theorem \ref{TheorInverseI}, defines an element 
$\bar \psi_{\theta_0} 
\in 
\divH^{\kappa}$. For any $\theta \in \divH^{\beta}$ consider a perturbation 
\begin{equation}\label{perturbio}
\theta_{(t, \psi)} 
\equiv 
\theta - \frac{t}{\sqrt N} p_N(\bar \psi_{\theta_0}), ~t \in \R,
\end{equation}
which, as a projection onto $E_{J_N}$, lies in the support $\divH^\beta$ of the prior. An application of the central limit theorem shows that the random variable
\begin{equation}\label{hatpsi}
\hat \Psi_{N} = \langle \psi, \theta_0 \rangle_{L^2} + \frac{1}{N} \sum_{i=1}^N \langle \mathbb I_{\theta_0}[ p_N (\bar \psi_{\theta_0})](X_i), \varepsilon_i \rangle_{\R^2}.
\end{equation}
satisfies the convergence in distribution $\sqrt N (\hat \Psi_{N} - \langle \theta_0, \psi \rangle_{L^2}) \to^d \mathcal N(0, \|\mathbb I_{\theta_0}[\bar \psi_{\theta_0}]\|_{L^2(\mathcal X)}^2)$ for each $\psi \in \divH\cap C^\infty(\Om)^2$, where we recall that $\XX$ is the time-space cylinder $[0,T]\times\Om$. We then prove as in \cite{Nickl2024} that:

\begin{Theor}\label{laplaceapp}
For all $\psi \in \divH^{\kappa+2}$ and $t\in\R$ fixed, the localized posterior from (\ref{respr}) satisfies 
\begin{eqnarray*}
\lefteqn{
E^{\Pi^{\bar \Theta_N}}\big[\exp\big\{t \sqrt N \big(\langle \theta, \psi\rangle_{L^2} - \hat \Psi_{N} \big)\big\}|Z^{(N)}\big]
}
\\[2mm]
&&
\hspace{10mm}
= 
e^{\frac{t^2}{2} \|\mathbb I_{\theta_0}[\bar \psi]\|_{L^2(\mathcal X)}^2} 
\times 
\frac{\int_{\bar \Theta_N} e^{\ell_N(\theta_{(t, \psi)})}d\Pi(\theta)}{\int_{\bar \Theta_N} e^{\ell_N(\theta)}d\Pi(\theta)} 
\times 
e^{r_N(t, \psi)}
,
\end{eqnarray*}
where $\ell_N$ is as in (\ref{postlog}), and where for every $L>0$ and $t \in \R$ fixed, $r_N$ is such that
$$
\sup_{\|\psi\|_{\divH^{\kappa+2}} \le L}
r_N(t,\psi)=o_{P_{\theta_0}^\mathbb N}(1),~~\text{ as } N \to \infty
.
$$
\end{Theor}

 Next, arguing as in \cite{Nickl2024}, an application of the Cameron-Martin theorem entails that the ratio 
\begin{equation}\label{postbus00}
\frac{\int_{\bar \Theta_N} e^{\ell_N(\theta_{(t, \psi)})}d\Pi(\theta)}{\int_{\bar \Theta_N} e^{\ell_N(\theta)}d\Pi(\theta)}  
\end{equation}
featuring in Theorem \ref{laplaceapp} above equals 
\begin{align} \label{postbus0}
(1+o(1))\frac{\int_{\bar \Theta_{N, t, \psi}} e^{\ell_N(\vartheta)}d\Pi(\vartheta)}{\int_{\bar \Theta_N} e^{\ell_N(\theta)}d\Pi(\theta)}
&=
(1 +o(1)) \frac{\Pi(\bar \Theta_{N,t, \psi}|Z^{(N)})}{\Pi(\bar \Theta_{N}|Z^{(N)})}
,
\end{align}
for shifted integration regions 
$$\bar \Theta_{N,t, \psi} = \{\theta_{(t,\psi)} : \theta \in \bar \Theta_N\}
.
$$ 
One then shows that, for each $\psi \in \divH\cap  C^\infty(\Om)^2$ fixed, the ratio in (\ref{postbus0}) satisfies
\begin{equation} \label{postlim}
\frac{\Pi(\bar \Theta_{N,t, \psi}|Z^{(N)})}{\Pi(\bar \Theta_{N}|Z^{(N)})} \to^{P_{\theta_0}^\N} 1,
\end{equation}
which shows that the Laplace transform on the l.h.s.~in Theorem \ref{laplaceapp} converges to the Laplace transform of a centred normal distribution with variance $\|\mathbb I_{\theta_0}[\bar \psi]\|_{L^2(\mathcal X)}^2$ for every smooth $\psi$, from which yields the convergence of finite-dimensional distributions used to prove Proposition~\ref{PropApproxBvM} below. The ratio (\ref{postbus00}) is further upper-bounded, uniformly in $\psi$ and for each $t \in \R$, by
\begin{align}\label{postbus000}
(1 +o(1)) \frac{1}{\Pi(\bar \Theta_{N}|Z^{(N)})} 
=
(1+o_{P_{\theta_0}^\mathbb N}(1)) 
\end{align}
by virtue of Theorem \ref{postcont}, which gives a sub-Gaussian bound on the posterior marginals.\medskip

One proves Theorem \ref{TheorBvMTheta2} initially with a centring $\hat \theta_N$ that is different from $\tilde \theta_N=E^\Pi[\theta|Z^{(N)}]$: taking $\psi=f_j$, the orthonormal basis of $\divH$ from (\ref{eq:BasisDivH}) with corresponding $\hat \Psi_N=\hat \Psi_{N,j}$, we obtain the (via $Z^{(N)}$ random) Fourier series 
\begin{equation} \label{hatter}
\hat \theta_N = \sum_{j \ge 1}\hat \Psi_{N,j} f_j.
\end{equation}
This defines a random variable in $\divH^{-k-2}$ with 
\begin{align}\label{momhat}
NE_{\theta_0}^\mathbb N \|\hat \theta_N - \theta_0\|^2_{H^{-k-2}} 
\le 
C
,
\end{align}
for all $N\ge 1$ and a universal constant $C>0$.
The next result establishes convergence of the localised posterior from (\ref{respr}) centred at $\hat \theta_N$ towards the desired Gaussian $\NN_{\theta_0}$ constructed in Section \ref{limproc}. For this purpose, define conditional laws in $\divH^{-k-2}$
\begin{equation}\label{bartau}
\bar \tau_N = {\rm Law} \big(\sqrt N (\theta - \hat \theta_N)|Z^{(N)}\big)
,\spa 
\theta \sim \Pi^{\bar \Theta_N}(\cdot|Z^{(N)})
.
\end{equation}

 \begin{Prop}\label{PropApproxBvM}
 We have, as $N\to\infty$
 $$
 \WWW_{1,\divH^{-k-2}}(\bar \tau_N, \mathcal N_{\theta_0}) \to^{P_{\theta_0}^\mathbb N} 0
 .
 $$
 \end{Prop}

The proof is again as in \cite{Nickl2024} using Theorem \ref{laplaceapp}, combining a uniform tightness estimate with convergence of finite-dimensional distributions.  The preceding proposition together with Proposition \ref{wassconv} for $T=\hat \theta_N$ implies, for unrestricted posterior draws $\theta \sim \Pi(\cdot|Z^{(N)})$, that 
\begin{equation}\label{hatlim}
\WWW_{1,\divH^{-k-2}}({\rm Law}(\sqrt N (\theta - \hat \theta_N), \mathcal N_{\theta_0}) \to^{P_{\theta_0}^\mathbb N} 0
\end{equation}
as $N \to \infty$, that is, Theorem \ref{TheorBvMTheta} holds for the centering $\hat \theta_N$ instead of $\tilde \theta_N=E^\Pi[\theta|Z^{(N)}]$. Since convergence in $\WWW_1$ implies convergence of first moments, we obtain for $Z\sim \NN_{\theta_0}$
$$
\sqrt N (E^\Pi[\theta |Z^{(N)}] -\hat \theta_N)  \to^{P_{\theta_0}^\mathbb N} E_{\mathcal N_{\theta_0}}Z =0 \text{ in } H^{-k-2},$$ 
and then also
\begin{eqnarray*}
\lefteqn{
\WWW_{1,\divH^{-k-2}}\big({\rm Law}(\sqrt N (\theta - \hat \theta_N)), {\rm Law}(\sqrt N(\theta- \tilde \theta_N))\big)
}
\\[2mm]
&&\hspace{10mm}
\le 
\sqrt N \|E^\Pi[\theta|Z^{(N)}] - \hat \theta_N\|_{H^{-k-2}} = o_{P_{\theta_0}^\mathbb N}(1),
\end{eqnarray*}
which, combined with (\ref{hatlim}), yields the first limit in Theorem~\ref{TheorBvMTheta2}. The second limit in Theorem~\ref{TheorBvMTheta}, namely the asymptotics for the posterior mean, is obtained as in \cite{Nickl2024}.

\subsection{Proof of Corollary \ref{CorolUQ}} \label{corprf}

 We note first that Exercise 2.6.5 in \cite{GinNic2016}, the discussion in Section \ref{limproc} as well as Proposition \ref{PropLinSmooth} imply that the RKHS $\mathcal H_\mu$ of the Gaussian measure $\mu={\rm Law}(U)$ on the separable space $\mathscr C$ equals the image of $\divH^{-1}$ under the linear solution map (\ref{linshow}). One can apply Theorem 11.10 in \cite{Robinson2001} with $h=0$ and appropriate $A$ and combine it with Proposition \ref{PropProdSob} to show that $\mathcal H_\mu$ is a nontrivial linear subspace of $\mathscr C$, in particular $\mu$ is non-degenerate. The proof of the corollary then proceeds as in Theorem 7.3.23 in \cite{GinNic2016}, noting that the distribution function
$\Phi(t)=\mu(\|U\|_\mathscr C \le t)$ is uniformly continuous (Exercise 2.4.4 in \cite{GinNic2016}) and \textit{strictly} increasing on $[0,\infty)$. The last claim follows from Theorem 2.4.5 and Corollary 2.6.18 in \cite{GinNic2016} and since any shell $\{s<\|u\|_\mathscr C <t\}, s<t,$ contains a $\CCC$-ball of radius $(t-s)/3$ centred at an element of $\mathcal H_\mu$.

\section{Theory for the periodic 2D Navier-Stokes equations}\label{sec:TheoryNSE}

\subsection{Functional setting}\label{sec:setting}

\subsubsection{Usual Sobolev spaces}

Let us recall some notation from Section \ref{sec:Notation}. We have an orthonormal basis for the complex-valued square-integrable vector fields $v:\Om\to\R^2$, given by 
$$
h_k(x) 
= 
e^{2i\pi k\cdot x} 
,\spa 
k\in\Z^2
,\quad 
x\in\Om 
.
$$
For $\{k_j = (k_{j,1}, k_{j,2}) : j\geq 0 \}$ the enumeration of $\Z^2$ from (\ref{eq:enum}), the $h_{k_j}$'s form an $L^2(\Om)$-orthonormal basis of eigenfunctions of $-\Delta$ with corresponding eigenvalues $\lambda_0=0$ and $0<\lambda_{j-1}\leq \lambda_j\asymp j$ for $j\geq 2$. For $\alpha\in [0,\infty)$, we have the Sobolev space $H^\alpha$ of vector fields, with inner product
\begin{equation}\label{eq:SobMetric2}
    \ps{u}{v}_{H^\alpha}
    =
    \sum_{j\geq 0} (1+\lambda_j^\alpha) \big\{ \ps{u}{h_{k_j}}_{L^2} \cdot \overline{\ps{v}{h_{k_j}}}_{L^2}\big\} 
    ,
\end{equation}
and the so-called `homogeneous' Sobolev space $\dH^\alpha$ obtained as the closed subspace of vector fields $u\in H^\alpha$ such that $\int_\Om u = 0$, \ie $\ps{u}{h_{k_0}}_{L^2}=0$, endowed with the inner product 
\begin{equation}\label{eq:dotSobMetric}
    \ps{u}{v}_{\dH^\alpha}
    =
    \sum_{j\geq 1} \lambda_j^\alpha \ps{u}{h_{k_j}}_{L^2}\cdot \overline{\ps{v}{h_{k_j}}}_{L^2}
    .
\end{equation}
The quantity $\|u\|_{\dH^\alpha}$ is well-defined for any $u\in H^\alpha$, but only defines a seminorm on $H^\alpha$ since elements of the (non-empty) eigenspace of the periodic Laplacian associated with $\lambda_0=0$ have zero $\dH^\alpha$-norm. For $\alpha>0$, we define $\dH^{-\alpha}$ as the completion of $L^2$ (or, equivalently, $C^\infty$) with respect to the $\dH^{-\alpha}$-norm obtained from (\ref{eq:dotSobMetric}) by replacing $\alpha$ by $-\alpha$. Since $\lambda_j>0$ for all $j\geq 1$, it is straightforward to see that 
\begin{equation}\label{dualSob1}
(\dH^\alpha)^* 
\cong 
\dH^{-\alpha}
,\spa  \alpha\in\R
,
\end{equation}
with
\begin{equation}\label{dualSob12}
\|u\|_{\dH^{-\alpha}}
=
\sup_{\substack{\psi\in \dH^{\alpha}\\ \|\psi\|_{\dH^{\alpha}}\leq 1}} \ps{u}{\psi}_{L^2}
,\spa 
\forall\ u\in \dH^{-\alpha}
,\quad \alpha\in\R
.
\end{equation}
By approximation, the above supremum can be restricted to $\dH^\infty = \dH^\alpha\cap C^\infty$. It follows that, for all $\alpha\in\R$, we have the Cauchy-Schwarz inequality
\begin{equation}\label{eq:dualProd}
\lvert\ps{u}{v}_{L^2}\rvert
\leq 
\|u\|_{\dH^{-\alpha}} \|v\|_{\dH^\alpha}
,\spa 
\forall\ u\in \dH^{-\alpha},\ v\in \dH^{\alpha}
.
\end{equation}
It further follows from the definitions that 
$$
\|u\|^2_{H^\alpha}
=
\|u\|^2_{L^2}
+
\|u\|^2_{\dH^\alpha}
,\spa 
\forall\ u\in H^\alpha
,\quad 
\alpha\geq 0
.
$$
For all $\alpha\in\R$, we define the diagonal operator $(-\Delta)^\alpha$ by 
$$
(-\Delta)^\alpha u
=
\sum_{j\geq 1} \lambda_j^\alpha \ps{u}{h_{k_j}}_{L^2} h_{k_j}
,\spa 
\forall\ u\in H^{2\alpha}
.
$$
In particular, we have 
$$
(-\Delta)^\alpha h_{k_j}
=
\lambda_j^\alpha h_{k_j}
,\spa 
j\geq 0 
,
$$
and
\begin{equation}\label{eq:LapNorm}
\|u\|^2_{\dH^\alpha}
=
\|(-\Delta)^{\alpha/2}u\|^2_{L^2}
,\spa 
\forall\ u\in H^\alpha
.
\end{equation}
For all $\alpha,\beta,\gamma \in \R$, and $u,v\in C^\infty$, our periodic setting yields the following integration by parts formulae
\begin{equation}\label{IppLaplace}
\langle(-\Delta)^\alpha u , (-\Delta)^\beta v \rangle_{\dH^\gamma}
=
\langle u , (-\Delta)^{\alpha+\beta} v\rangle_{\dH^\gamma}
=
\langle u, (-\Delta)^{\alpha+\beta+\gamma} v \rangle_{L^2}
.
\end{equation}
In particular, 
\begin{equation}\label{IppSob1}
    \|(-\Delta)^\alpha u\|_{\dH^\beta}
    =
    \|u\|_{\dH^{\beta+2\alpha}}
    ,
\end{equation}
and
\begin{equation}\label{IPPSob2}
\ps{u}{-\Delta u}_{\dH^\alpha}
=
\langle (-\Delta)^{1/2} u , (-\Delta)^{1/2}u \rangle_{\dH^\alpha}
=
\|(-\Delta)^{1/2} u\|^2_{\dH^\alpha}
=
\|u\|^2_{\dH^{\alpha+1}}
.
\end{equation}
We further have 
\begin{equation}\label{IPPSob3}
    \|\nabla u\|_{\dH^\alpha}
    =
    \|u\|_{\dH^{\alpha+1}}
    ,
\end{equation}
since, by standard integration by parts and since $\Om$ has no boundary,
\begin{eqnarray*}
    \lefteqn{
    \|\nabla u\|^2_{\dH^\alpha}
    =
    \sum_{i=1}^2 \|\partial_i u\|^2_{\dH^\alpha}
    =
    \sum_{i=1}^2 \sum_{j\geq 1} \lambda_j^\alpha \lvert \ps{\partial_i u}{h_{k_j}}_{L^2}\rvert^2
    =
    \sum_{i=1}^2 \sum_{j\geq 1} \lambda_j^\alpha \lvert \ps{u}{\partial_i h_{k_j}}_{L^2}\rvert^2
    }
    \\[2mm]
    &&
    =
    \sum_{j\geq 1} \lambda_j^\alpha \lvert \ps{u}{h_{k_j}}_{L^2}\rvert^2 \sum_{i=1}^2 (2\pi [k_j]_i)^2 
    =
    \sum_{j\geq 1}  \lambda_j^{\alpha+1} \lvert \ps{u}{h_{k_j}}_{L^2}\rvert^2
    =
    \|u\|^2_{\dH^{a+1}}
    .
\end{eqnarray*}
For all $\alpha,\beta\in\R$ with $\alpha\leq \beta$, and $u\in C^\infty$ we have the Poincaré type inequality
\begin{equation}\label{Poincare}
\|u\|^2_{\dH^\alpha}
=
\sum_{j\geq 1} \lambda_j^\alpha \lvert\ps{u}{h_{k_j}}_{L^2}\rvert^2
\leq 
\lambda_1^{-(\beta-\alpha)} \sum_{j\geq 1} \lambda_j^\beta \lvert\ps{u}{h_{k_j}}_{L^2}\rvert^2 
=
\lambda_1^{-(\beta-\alpha)} \|u\|^2_{\dH^\beta}
.
\end{equation}
As a consequence, on $\dH^\alpha$ with $\alpha\geq 0$, we have equivalence of the $H^\alpha$ and $\dH^\alpha$ norms since 
\begin{equation}\label{equivSob}
\|u\|^2_{\dH^\alpha}
\leq 
\|u\|^2_{H^\alpha}
\leq 
(1+\lambda_1^{-\alpha})\|u\|^2_{\dH^\alpha}
,\spa 
\forall\ u\in \dH^\alpha
,\quad 
\alpha \ge 0
.
\end{equation}
Finally, recalling the notation in Section \ref{sec:Notation}, for any $a\in\Z$ and $u\in L^2([0,T], \dH^{a+1})$ such that $du/dt\in L^2([0,T], \dH^{a-1})$, we can prove as in Corollary 7.3 of \cite{Robinson2001} that
\begin{equation}\label{IPPdt}
\Big\langle \frac{du}{dt}, u(t)\Big\rangle_{\dH^a}
=
\frac12 \frac{d}{dt}\|u(t)\|^2_{\dH^a}
.
\end{equation}

\subsubsection{The homogeneous divergence-free Sobolev scale}

We also consider the space of divergence-free and mean-zero vector fields 
$$
\Hnull 
=
\bigg\{ u \in L^2 : \nabla\cdot u = 0,\ \int_\Om u = 0 
\bigg\} 
,
$$
endowed with the $L^2$-inner product. Consider an arbitrary smooth function $u\in \Hnull\cap C^\infty$. Since $u\in L^2$, we have the orthogonal decomposition 
$$
u(x)
=
\sum_{j\geq 0} u_j h_{k_j}(x)
,\spa 
u_j
=
\ps{u}{h_{k_j}}_{L^2}
\in \C^2
,
$$ 
with $(u_j)\in \ell^2(\N)$. Noting that $k_0 = (0,0)$, we then have 
$$
(\nabla\cdot u)(x)
= 2i \pi \sum_{j\geq 1} (u_j\cdot k_j) h_{k_j}(x).$$
In particular, the condition $\nabla\cdot u=0$ writes $ u_j\cdot k_j=0$ for all $j\geq 1$, which implies that
$$
u_j 
=
\tilde{u}_j \frac{(k_{j,2},-k_{j,1})}{|k_j|}
,\spa 
j\geq 1
,
$$
for some sequence $(\tilde u_j)_{j\geq 1}\in\ell^2(\N\sm\{0\})$. Since the condition $\int_{\Om} u=0$ corresponds to $u_0=0$, we deduce that 
$$
u
=
\sum_{j\geq 1} \tilde{u}_j c_j h_{k_j}
,\spa 
c_j 
=
\frac{(k_{j,2}, -k_{j,1})}{|k_j|}
,\quad 
\tilde{u}_j
=
c_j\cdot \int_\Om u(x) \overline{h_{k_j}(x)}\, dx
.
$$
Since any such $u$ also defines an element of the complex-valued extension of $\Hnull$, we obtain the $L^2$-orthonormal basis $\{e_j : j\geq 1\}$ of the complex extension of $\Hnull$, given by 
$$
e_j(x)
=
c_j
h_{k_j}(x)
,\spa  
j\geq 1
,
$$ 
so that any $u\in \Hnull$ writes 
$$
u
=
\sum_{j\geq 1} \ps{u}{e_j}_{L^2}e_j
,\spa 
\ps{u}{e_j}_{L^2}
\equiv 
\int_\Om u(x)\cdot \overline{e_j(x)}\, dx 
,
$$
and $u$ will be real-valued provided conjugacy relations as in (\ref{eq:L2Series}) are satisfied. We define the corresponding Sobolev spaces 
$$
(\divH^\alpha, \ps{\cdot}{\cdot}_{\divH^\alpha})
\equiv 
(\dH^\alpha\cap \Hnull, \ps{\cdot}{\cdot}_{\dH^a})
,\spa 
\alpha\in [0,\infty)
.
$$
For $\alpha>0$, we define $\divH^{-\alpha}$ as the completion of $\Hnull=\divH^0$ (or, equivalently, $\divH^\infty\equiv C^\infty\cap \Hnull$) with respect to the $\dH^{-\alpha}$ norm. 

Let us remark that the equality $\ps{e_j}{h_{k_\ell}}_{L^2}=c_\ell 1\{j = \ell\}$ for all $j\geq 1$ entails that for any $u\in\Hnull$, with $u=\sum_{j\geq 1} \ps{u}{e_j}_{L^2}e_j$, we have
$$
\ps{u}{h_{k_\ell}}_{L^2} h_{k_\ell}
=
\Big(\sum_{j\geq 1} \ps{u}{e_j}_{L^2} \overline{\ps{e_j}{h_{k_\ell}}}_{L^2}\Big) h_{k_\ell}
=
\ps{u}{e_\ell}_{L^2} c_\ell  h_{k_\ell}
=
\ps{u}{e_\ell}_{L^2} e_\ell
.
$$
In particular, we have for all $\alpha\in \R$
\begin{equation}\label{eq:defLap}
(-\Delta)^\alpha u 
=
\sum_{j\geq 1} \lambda_j^\alpha \ps{u}{h_{k_j}}_{L^2} h_{k_j}
=
\sum_{j\geq 1} \lambda_j^\alpha \ps{u}{e_j}_{L^2} e_j 
,\spa 
\forall\ u\in \divH^\alpha
,
\end{equation}
and 
$$
(-\Delta)^\alpha e_j 
=
\lambda_j e_j 
,\spa j\geq 1,$$
so that, as a consequence, the $\divH^\alpha$-norms can be computed in either basis.
Finally, since $\lambda_j>0$ for all $j\geq 1$, it is straightforward to see that 
\begin{equation}\label{dualSob2}
(\divH^\alpha)^* 
\cong 
\divH^{-\alpha}
,\spa 
\forall\ \alpha\in\R
,
\end{equation}
with
\begin{equation}\label{dualSob3}
\|u\|_{\dH^{-\alpha}}
=
\sup_{\substack{\psi\in \divH^{\alpha}\\ \|\psi\|_{\dH^{\alpha}}\leq 1}} \ps{u}{\psi}_{L^2}
=
\sup_{\substack{\psi\in \dH^{\alpha}\\ \|\psi\|_{\dH^{\alpha}}\leq 1}} \ps{u}{\psi}_{L^2}
,\spa 
\forall\ u\in \divH^{-\alpha}
.
\end{equation}
By approximation, the above supremum can be restricted to $\divH^\infty$. It follows as in (\ref{eq:dualProd}) that, for all $\alpha\in\R$, we have the Cauchy-Schwarz inequality
\begin{equation}\label{eq:dualProd2}
\lvert\ps{u}{v}_{L^2}\rvert
\leq 
\|u\|_{\dH^{-\alpha}} \|v\|_{\dH^\alpha}
,\spa 
\forall\ u\in \divH^{-\alpha},\ v\in \divH^{\alpha}
.
\end{equation}

\subsubsection{Leray projector}

Since $\Hnull$ is a closed subspace of $L^2$, one has the so-called `Leray projector', \ie the $L^2$-orthogonal projection 
$$
P:L^2\to \Hnull
,
$$
characterized, for all $u\in L^2$, by 
\begin{equation}\label{eq:PCharact}
    \|P[u]-h\|_{L^2}
    \leq 
    \|u-h\|_{L^2}
    ,\spa 
    \forall\ h\in \Hnull 
    .
\end{equation}
Using (\ref{eq:PCharact}), it is straightforward to see that
\begin{equation}\label{eq:ProjectionP}
P[u]
=
\sum_{j\geq 1} \Big(u_j - \frac{(u_j\cdot k_j) k_j}{|k_j|^2} \Big) h_{k_j}
,\spa 
\forall\ u\in L^2
.
\end{equation}
Since (\ref{eq:ProjectionP}) is homogeneous in the coefficients $(u_j)$, one deduces that $P$ commutes with any operator which is diagonal in the basis $(h_{k_j})_{j\geq 0}$, \ie any operator $\Lambda$ for which there exists a sequence $(s_j)_{j\geq 0}\subset \R$ such that $\Lambda[v] = \sum_{j\geq 0} s_j \ps{v}{h_{k_j}}_{L^2} h_{k_j}$ for all $v$ making the last expression well-defined. 
It follows that for all $\alpha\in\R$ we have
\begin{equation}\label{eq:PCommutesLaplacians}
    P(-\Delta)^\alpha[u]
    =
    (-\Delta)^\alpha P[u]
    ,\spa
    \forall\ u \in C^\infty
    .
\end{equation}

\subsection{Estimates for the bilinear form}

We provide estimates on the bilinear form $B$ defined in (\ref{eq:defB}) as 
$$
B[u,v]
\equiv 
P[(u\cdot \nabla)v]
,
$$
where $(u\cdot \nabla)v$ is the vector-field obtained as the matrix product $(\nabla v)u$. The identities stated in Proposition \ref{PropBilinearFormIdentities1}, Proposition \ref{PropBilinearFormEstimates}, and Proposition \ref{PropProdSob} below are key to proving the results of Section \ref{sec:NSEregularity} and Section \ref{sec:Lin}.  A key ingredient in these proofs is the well-known \lad's inequality which, in the form we present next, is valid on a general \textit{two-dimensional} domain: for any $u\in H^1$ one has 
\begin{equation}\label{LadInequality}
\|u\|_{L^4}
\lesssim 
\|u\|^{1/2}_{L^2}\|\nabla u\|^{1/2}_{L^2}
.
\end{equation}
This can be proved by the Sobolev embedding $H^{1/2}\embed L^4$ and interpolation of $H^{1/2}$ between $L^2$ and~$H^1$. \lad's inequality combined with the Cauchy-Schwarz inequality then yields
\begin{equation}\label{eq:ProdL2}
\|uv\|_{L^2}
\leq 
\|u\|_{L^4}\|v\|_{L^4}
\lesssim 
\|u\|^{1/2}_{L^2}\|\nabla u\|^{1/2}_{L^2} \|v\|^{1/2}_{L^2}\|\nabla v\|^{1/2}_{L^2}
.
\end{equation}
The following proposition collects usual estimates on the bilinear form $B$ and can be found, e.g., in Section 9.2 of \cite{Robinson2001}.

\begin{Prop}\label{PropBilinearFormIdentities1}
    Let $u,v,w\in \divH^1$. Then 
    \begin{enumerate}
        \item[(i)] $\ps{B[u,v]}{w}_{L^2}=-\ps{B[u, w]}{v}_{L^2}$, so that $\ps{B[u,v]}{v}_{L^2}=0$;

        \item[(ii)] If $v\in \divH^2$, then $\ps{B[v, v]}{\Delta v}_{L^2}=0$;

        \item[(iii)] If $v\in\divH^2$, then $\|B[u,v]\|_{L^2}\lesssim \|u\|^{1/2}_{L^2}\|\nabla u\|^{1/2}_{L^2} \|\nabla v\|^{1/2}_{L^2}\|\Delta v\|^{1/2}_{L^2}$.
    \end{enumerate}
\end{Prop}
\vspace{1mm}

In the next proposition we generalize the bounds of Proposition \ref{PropBilinearFormIdentities1} to (positive and negative order) Sobolev inner products, which will be needed in  Section \ref{sec:NSEregularity} to establish Sobolev regularity estimates for solutions of the Navier-Stokes equations, and in Section \ref{sec:Lin} to derive similar estimates for the linear approximation of the Navier-Stokes flow. 
We first need to establish the projection bounds (\ref{eq:SobolevProjectionNorm}) and (\ref{eq:BoundBSobolev}), and recall the multiplier inequality for Sobolev norms (\ref{multipSob}). For all $a\in\R$, recall from (\ref{IppSob1}) that $\|u\|^2_{\dH^a} = \|(-\Delta)^{a/2} u\|^2_{L^2}$ for all $u\in\divH^\infty$. Because $\|P[v]\|_{L^2}\leq \|v\|_{L^2}$ for all $v\in L^2$, by (\ref{eq:PCharact}) with $h=0$, and since we established in (\ref{eq:PCommutesLaplacians}) that $P(-\Delta)^{a/2}=(-\Delta)^{a/2}P$ on $H^a$, we then have the projection bound
\begin{equation}\label{eq:SobolevProjectionNorm}
\|P[u]\|_{\dH^{a}}
\leq 
\|u\|_{\dH^a}
,\spa 
\forall\ u\in C^\infty
,\quad 
a\in\R
.
\end{equation}
This obviously extends to any $u\in H^a$ by approximation and further yields 
\begin{equation}\label{eq:BoundBSobolev}
\|B[u,v]\|_{\dH^a}
\leq 
\|(u\cdot \nabla) v\|_{\dH^a}
=
\|(\nabla v)u\|_{\dH^a}
,\spa 
\forall\ u,v\in C^\infty
,\quad
a\in\R
.
\end{equation}
For any $\ell\in\N$, define function spaces $\BB^\ell\equiv C^\ell(\Omega)^2$ if $0\leq \ell \leq 1$ and $\BB^\ell \equiv \dH^\ell$ if $\ell\geq 2$, endowed with their usual norm. Then, we have the multiplier inequality for Sobolev norms
\begin{equation}\label{multipSob}
    \Big(\sum_{i,j=1}^2 \|u_i v_j\|^2_{\dH^a}\Big)^{1/2} 
    \equiv
    \|uv^T\|_{\dH^a}
    \lesssim
    \|u\|_{\BB^{|a|}}\|v\|_{\dH^a}
    ,\spa 
    \forall\ u,v\in \dot{C}^\infty
    ,\quad
    a\in\Z
    .
\end{equation}
This estimate is known to be valid for all $u,v\in C^\infty$ when the $\dH^a$-norm is replaced by the usual $H^a$-norm; see, e.g., Theorem 4.39 in \cite{AF03}. The version given in (\ref{multipSob}) further exploits that $u,v\in\dot{C}^\infty$ and thus follows from the equivalence of the $\dH^a$ and $H^a$ norms for mean-zero vector fields established in (\ref{equivSob}). 

\begin{Lem}\label{PropBilinearFormIdentities2}
    Fix $a\in\Z$, and $u,v,w\in \divH^\infty$. Then,
    \begin{enumerate}
        \item[(i)] $\|B[u,v]\|_{\dH^{a}}\leq \|uv^T\|_{\dH^{a+1}}$;
        
        \item[(ii)] $|\ps{B[u,v]}{w}_{\dH^a}|\lesssim \min\{\|u\|_{\BB^{|a|}}\|v\|_{\dH^a}, \|u\|_{\dH^a}\|v\|_{\BB^{|a|}}\} \| w\|_{\dH^{a+1}}$.
    \end{enumerate}
    In particular, these estimates extend to the case where $u,v,w$ have the Sobolev regularity prescribed by the corresponding inequality.
\end{Lem}

\begin{Proof}{Lemma \ref{PropBilinearFormIdentities2}}
   (i). By an adaptation of Proposition \ref{PropBilinearFormIdentities1}(i) to complex-valued functions, we have $\ps{B[u,v]}{e_j}_{L^2}=-\ps{B[u, \bar e_j]}{v}_{L^2}$ for all $j\geq 1$. Using that $v\in\Hnull$, so that $\ps{P[w]}{v}_{L^2}=\ps{w}{v}_{L^2}$ for any $w\in C^\infty$, and $\nabla \bar e_j = 2i\pi \bar e_j k_j^T$ for all $j\geq 1$, yields 
   \begin{eqnarray*}
       \lefteqn{
       \|B[u,v]\|^2_{\dH^a}
       =
       \sum_{j\geq 1} \lambda_j^a \lvert\ps{B[u,v]}{e_j}_{L^2}\rvert^2
       =
       \sum_{j\geq 1} \lambda_j^a \lvert\ps{B[u,\overline {e_j}]}{v}_{L^2}\rvert^2
       }
       \\
       &&
       =
       \sum_{j\geq 1}\lambda_j^a \lvert\ps{(u\cdot \nabla)\overline {e_j}}{v}_{L^2}\rvert^2 
       =
       \sum_{j\geq 1}\lambda_j^a \lvert\ps{(\nabla \overline {e_j})u}{v}_{L^2}\rvert^2
       =
       \sum_{j\geq 1}\lambda_j^a (2\pi)^2 \lvert\ps{\overline{e_j} k_j^T u}{v}_{L^2}\rvert^2
       .
   \end{eqnarray*}
   Using the Cauchy-Schwarz inequality twice, we find
   \begin{eqnarray*}
       \lefteqn{
       \hspace{-10mm}
       \lvert\ps{\bar e_j k_j^T u}{v}_{L^2}\rvert
       =
       \Big|k_j \cdot  \int_\Om u(x) (v(x)\cdot \bar e_j(x))\, dx\Big|
       }
       \\[2mm]
       &&
       \leq 
       |k_j| \sqrt{\sum_{i=1}^2 \Big( \int_\Om u_i(x)( v(x)\cdot \bar e_j(x) )\, dx \Big)^2}
       \\[2mm]
       &&
       =
       |k_j| \sqrt{\sum_{i=1}^2 \Big( \sum_{\ell=1}^2 [c_j]_\ell \int_\Om u_i(x) v_\ell(x) \bar h_{k_j}(x)\, dx\Big)^2}
       \\[2mm]
       &&\leq 
       |k_j| \sqrt{\sum_{i=1}^2  \sum_{\ell=1}^2 \Big( \int_\Om u_i(x) v_\ell(x) \bar h_{k_j}(x)\, dx\Big)^2} 
       ,
   \end{eqnarray*}
   where we used the fact that $\sum_{\ell=1}^2 [c_j]_\ell^2 = |c_j|^2=1$ for all $j\geq 1$. Recalling that $\lambda_j = (2\pi|k_j|)^2$, we find 
   $$
   \|B[u,v]\|^2_{\dH^a}
   \leq 
   \sum_{i,\ell=1}^2 \sum_{j\geq 1} \lambda_j^{a+1} \lvert\ps{u_iv_\ell}{h_{k_j}}_{L^2}\rvert^2
   =
   \sum_{i,\ell=1}^2 \|u_iv_\ell\|^2_{\dH^{a+1}}
   .
   $$

   (ii). Using the inequality established in Part (i) of the present proof and (\ref{eq:dualProd2}), we have
   $$
   \lvert\ps{B[u,v]}{w}_{\dH^a}\rvert
   \leq 
   \|B[u,v]\|_{\dH^{a-1}}\|w\|_{\dH^{a+1}}
   \leq 
   \|u v^T\|_{\dH^a} \|w\|_{\dH^{a+1}}
   .
   $$
   The multiplier inequality for Sobolev norms (\ref{multipSob}) yields
   \begin{align*}
   \|u v^T\|_{\dH^{a}}
   &\lesssim 
   \min\{\|u\|_{\BB^{|a|}}\|v\|_{\dH^a}, \|u\|_{\dH^a}\|v\|_{\BB^{|a|}}\}
   ,
   \end{align*}
   which concludes the proof.
\end{Proof}
\vspace{3mm}

Building on the previous result, we are to able to prove the following estimate, that will be crucial in the subsequent sections. For any $a\in\Z$, define 
\begin{equation}\label{eq:DefinAStar}
a^*
\equiv
\begin{cases}
    |a|+1 & \textrm{ if } |a|\leq 1,\\[2mm]
    |a| & \textrm{ if } |a|\ge 2,\\[2mm]
\end{cases}
\end{equation}

\begin{Prop}\label{PropBilinearFormEstimates}
    Let $a\in\Z$ and $a^*$ as in (\ref{eq:DefinAStar}). Then, for all $v\in\divH^{a^*}$ and $w\in\divH^{a+1}$ we have
    $$
    \lvert\ps{B[v,w]}{w}_{\dH^a}\rvert
    +
    \lvert\ps{B[w,v]}{w}_{\dH^a}\rvert
    \lesssim 
    \|v\|_{\dH^{a^*}}\|w\|_{\dH^a}\|w\|_{\dH^{a+1}}
    .
    $$
\end{Prop}


\vspace{1mm}

\begin{Proof}{Proposition \ref{PropBilinearFormEstimates}}
    Since $\divH^\infty\equiv C^\infty\cap\Hnull$ is dense in $\divH^{a^*}$ and $\divH^a$, we may assume without loss of generality that $v,w\in \divH^\infty$ to deduce that the desired estimate holds whenever $v\in\divH^{a^*}$ and $w\in\divH^{a+1}$. When $|a|\geq 2$, then $a^* = |a|$ and the result follows from Lemma \ref{PropBilinearFormIdentities2}(ii). It remains to establish the result when $|a|\leq 1$.
    Assume that $a=1$. Then $a^* = 2$. We have by (\ref{IppLaplace}), Proposition \ref{PropBilinearFormIdentities1}(i), and (\ref{eq:dualProd2})
   $$
   \lvert\ps{B[w,v]}{w}_{\dH^1}\rvert
   =
   \lvert\ps{B[w,v]}{\Delta w}_{L^2}\rvert
   =
   \lvert\ps{B[w,\Delta w]}{v}_{L^2}\rvert
   \leq 
   \|B[w,\Delta w]\|_{\dH^{-2}}\|v\|_{\dH^2}
   .
   $$
   The projection bound (\ref{eq:BoundBSobolev}) and the multiplier inequality for Sobolev norms (\ref{multipSob}) provide 
   \begin{eqnarray*}
   \lefteqn{
    \|B[w,\Delta w]\|_{\dH^{-2}}
   \leq
   \|(w\cdot \nabla) \Delta w\|_{\dH^{-2}}
   =
   \|(\nabla \Delta w) w\|_{\dH^{-2}}
    }  
   \\[2mm]
   &&\lesssim 
   \|\nabla\Delta w\|_{\dH^{-2}}\|w\|_{\dH^2}
   =
   \|\Delta w\|_{\dH^{-1}}\|w\|_{\dH^2}
   =
   \|w\|_{\dH^1}\|w\|_{\dH^2}
   ,
   \end{eqnarray*}
   where we used the integration by parts formulae (\ref{IPPSob3}) and (\ref{IppSob1}) in the last two equalities.
   We also have
   $$
   \lvert\ps{B[v,w]}{w}_{\dH^1}\rvert
   =
   \lvert\ps{B[v,w]}{\Delta w}_{L^2}\rvert
   =
   \lvert\ps{B[v,\Delta w]}{w}_{L^2}\rvert
   \leq 
   \|B[v,\Delta w]\|_{\dH^{-2}}\|w\|_{\dH^2}
   .
   $$
   The previous computation provides
   $$
   \|B[v,\Delta w]\|_{\dH^{-2}}
   \lesssim 
   \|v\|_{\dH^2}\| w\|_{\dH^1}
   ,
   $$
   and yields 
   $$
   \lvert\ps{B[v,w]}{w}_{\dH^1}\rvert
   \lesssim 
   \|v\|_{\dH^2} \|w\|_{\dH^1}\| w\|_{\dH^2}
   .
   $$
    Assume now that $a=0$. Then $a^*=1$. Proposition \ref{PropBilinearFormIdentities1}(i) entails that $\ps{B[v,w]}{w}_{L^2}=0$. As for $\ps{B[w,v]}{w}_{L^2}$, we have
    \begin{align*}
     \lvert\ps{B[w,v]}{w}_{L^2}\rvert
     =
     \lvert\ps{B[w, w]}{v}_{L^2}\rvert
     \leq 
     \|v\|_{\dH^1}\|B[w, w]\|_{\dH^{-1}}
     \leq
     \|v\|_{\dH^1}\|w w^T\|_{L^2} 
     ,
    \end{align*}
    where we used Lemma \ref{PropBilinearFormIdentities2}(i) in the last inequality. Since 
    $$
    \|ww^T\|^2_{L^2} 
    =
    \sum_{i,\ell=1}^2 \|w_i w_\ell\|^2_{L^2}
    =
    \int_\Om \sum_{i,\ell=1}^2 (w_iw_\ell)^2\
    =
    \|w\|^4_{L^4}
    ,
    $$
    \lad's inequality (\ref{LadInequality}) combined with (\ref{IPPSob3}) then provides
    $$
    \lvert\ps{B[w, v]}{w}_{L^2}\rvert
    \lesssim 
    \| v\|_{\dH^1}\|w\|_{L^2}\|\nabla w\|_{L^2}
    =
    \|v\|_{\dH^1}\|w\|_{L^2} \|w\|_{\dH^1}
    .
    $$
    Assume now that $a=-1$ with $a^*=2$. We have
    \begin{align*}
    \lvert\ps{B[v,w]}{w}_{\dH^{-1}}\rvert
    \leq 
    \|B[v,w]\|_{\dH^{-2}}\|w\|_{L^2}
    .
    \end{align*}
    The projection bound (\ref{eq:BoundBSobolev}) and the multiplier inequality for Sobolev norms (\ref{multipSob}) entail that 
    $$
    \|B[v,w]\|_{\dH^{-2}}
    \leq 
    \|(\nabla w)v\|_{\dH^{-2}}
    \lesssim 
    \|v\|_{\dH^2} \|\nabla w\|_{\dH^{-2}}
    =
    \|v\|_{\dH^2}\|w\|_{\dH^{-1}}
    ,
    $$
    where we used (\ref{IPPSob3}) in the last equality.
    For $\ps{B[w,v]}{w}_{\dH^{-1}}$, we have by Lemma \ref{PropBilinearFormIdentities2}(i) 
    $$
    \lvert\ps{B[w,v]}{w}_{\dH^{-1}}\rvert 
    \leq 
    \|B[w,v]\|_{\dH^{-1}}\|w\|_{\dH^{-1}}
    \lesssim 
    \|w v^T\|_{L^2} \|w\|_{\dH^{-1}}
    .
    $$
    By the Sobolev embedding $\dH^2\embed H^2\embed L^\infty$, and using (\ref{equivSob}), we have 
    $$
    \|w v^T\|_{L^2}
    \lesssim 
    \|v\|_{L^\infty} \|w\|_{L^2}
    \lesssim 
    \|v\|_{H^2} \|w\|_{L^2}
    \lesssim 
    \|v\|_{\dH^2} \|w\|_{L^2}
    .
    $$
    This establishes the desired estimate when $a=-1$ and concludes the proof.
\end{Proof}
\vspace{3mm
}

We conclude this section with a Sobolev bound on the bilinear form $B$.

\begin{Prop}\label{PropProdSob}
    Fix $a\in\Z$ and $b=\max\{|a|+1,2\}$. For all $v\in \divH^b$ and $w\in\divH^{a+1}$, we have 
    $$
    \|B[v,w]\|_{\dH^a}
    + 
    \|B[w,v]\|_{\dH^a}
    \lesssim 
    \|v\|_{\dH^b}\|w\|_{\dH^{a+1}}
    .
    $$
    In particular if $a\geq 1$, we have
    $$
    \|B[v,w]\|_{\dH^a}
    + 
    \|B[w,v]\|_{\dH^a}
    \lesssim 
    \|v\|_{\dH^{a+1}}\|w\|_{\dH^{a+1}}
    .
    $$
        
        
\end{Prop}



\begin{Proof}{Proposition \ref{PropProdSob}}
Lemma \ref{PropBilinearFormIdentities2}(i) entails that 
$$
\|B[v,w]\|_{\dH^a}
+ 
\|B[w,v]\|_{\dH^a}
\leq
\|v w^T\|_{\dH^{a+1}}+\|w v^T\|_{\dH^{a+1}}
.
$$
The multiplier inequality for Sobolev norms (\ref{multipSob}) provides 
$$
\|v w^T\|_{\dH^{a+1}}+\|w v^T\|_{\dH^{a+1}}
\lesssim 
\|v\|_{\BB^{|a+1|}}\|w\|_{\dH^{a+1}}
.
$$
When $a\geq 1$, we have $\BB^{|a+1|}= \dH^{a+1}=\dH^b$ since $a+1\geq 2$. When $a\leq -3$, we have 
$$
\BB^{|a+1|}
=
\BB^{-a-1}
=
\BB^{|a|-1} 
=
\dH^{|a|-1}
,
$$
since $|a|-1\geq 2$, so that $\|v\|_{\BB^{|a+1|}}\lesssim \|v\|_{\dH^b}$ since $b\geq |a|-1$. When $a=-2$, we have by the projection bound (\ref{eq:BoundBSobolev}) and the multiplier inequality for Sobolev norms (\ref{multipSob})
$$
\|B[v,w]\|_{\dH^{-2}}
\leq 
\|(v\cdot\nabla )w\|_{\dH^{-2}}
\lesssim 
\|v\|_{\dH^2}\|\nabla w\|_{\dH^{-2}}
= 
\|v\|_{\dH^2} \|w\|_{\dH^{-1}}
,
$$
where we used (\ref{IPPSob3}) in the last equality. We also have
$$
\|B[w,v]\|_{\dH^{-2}}
\leq 
\|(w\cdot\nabla )v\|_{\dH^{-2}}
\lesssim 
\|w\|_{\dH^{-2}}\|\nabla v\|_{\dH^2}
= 
\|w\|_{\dH^{-2}} \|v\|_{\dH^3}
.
$$
By the Poincaré inequality (\ref{Poincare}), we deduce that
$$
\|B[v,w]\|_{\dH^{-2}}
+
\|B[w,v]\|_{\dH^{-2}}
\lesssim 
\|v\|_{\dH^3}\|w\|_{\dH^{-1}}
=
\|v\|_{\dH^b} \|w\|_{\dH^{a+1}}
.
$$
When $a=-1$, we use Lemma \ref{PropBilinearFormIdentities2}(i) and the Sobolev embedding $\dH^2 \embed L^\infty$ to the effect that
\begin{eqnarray*}
\lefteqn{
\hspace{-10mm}
\|B[v,w]\|_{\dH^{-1}}
+
\|B[w,v]\|_{\dH^{-1}}
\lesssim 
\|w v^T\|_{L^2}
}
\\[2mm]
&&
\lesssim 
\|v\|_{L^\infty}\|w\|_{L^2}
\lesssim 
\|v\|_{\dH^2}\|w\|_{L^2}
=
\|v\|_{\dH^b}\|w\|_{\dH^{a+1}}
.
\end{eqnarray*}
When $a=0$, we use the Sobolev embedding $\dH^2 \embed L^\infty$ to obtain
$$
\|B[v,w]\|_{L^2}
\lesssim 
\|(v\cdot\nabla)w\|_{L^2}
\lesssim
\|v\|_{L^\infty}\|\nabla w\|_{L^2}
\lesssim
\|v\|_{\dH^2}\|w\|_{\dH^1}
,
$$
where we used (\ref{IPPSob3}) in the last inequality. Lemma \ref{PropBilinearFormIdentities1}(iii), combined with (\ref{IPPSob3}) and Poincaré's inequality, entails that 
$$
\|B[w,v]\|_{L^2}
\lesssim 
\|w\|_{\dH^1}\|v\|_{\dH^2}
.
$$
This yields
$$
\|B[v,w]\|_{L^2}
+
\|B[w,v]\|_{L^2}
\lesssim 
\|v\|_{\dH^2}\|w\|_{\dH^1}
=
\|v\|_{\dH^b}\|w\|_{\dH^{a+1}}
,
$$
and concludes the proof.
\end{Proof}

\subsection{Regularity theory for strong solutions}\label{sec:NSEregularity}

The space $\divH^1$ is commonly denoted $V$ in the literature on Navier-Stokes equations. Then, the topological dual $V'$ is given by $\divH^{-1}$, by virtue of (\ref{dualSob2}). We define the bilinear operator $B:V\times V\to V'$ as
$$
B[u,v]
=
P[(u\cdot \nabla)v]
,\spa 
\forall\ u,v\in V=\divH^1 
.
$$
The fact that $B$ takes values in $V'$ follows from Lemma \ref{PropBilinearFormIdentities2}(i) and \lad's inequality (\ref{eq:ProdL2}), which entail that
$$
\|B[u,v]\|_{\dH^{-1}}
\le 
\|u v^T\|_{L^2}
<
\infty,\spa 
\forall\ u,v\in \divH^1 
.
$$
One obtains a theory of \textit{weak} solutions of the periodic 2D Navier-Stokes system by taking $L^2$-inner products of the first equation with an arbitrary $v\in \divH^1$. Noticing that 
$$
\big\langle\frac{\partial u}{\partial t} , v\big\rangle_{L^2}
=
\frac{d}{dt} \ps{u}{v}_{L^2}
,~~
\ps{(u\cdot\nabla)u}{v}_{L^2} 
= \ps{B[u,u]}{v}_{L^2}
,~~ 
\ps{\nabla p}{v}_{L^2}  = - \int_\Om(\nabla\cdot v) p = 0
$$ 
for all $v\in \divH^1$, and since $f\in \Hnull$, it is enough to solve for $u:[0,T]\to V'$ satisfying the \textit{functional} equation in $V'$
$$
\frac{du}{dt} + A u + B[u,u]
=
f
,\spa
u(0)=\theta 
,
$$
where $A=-P\Delta$ is the Stokes operator. We call a solution $u$ \textit{strong} if the preceding equation holds in $L^2([0,T], \Hnull)$ and is continuous at zero, as in Part i) of the following result, which is proved in Proposition 1 of \cite{NicTit2024}. We note here that any strong solution then also provides a solution to the original equation (\ref{eq:NSequations}) by standard arguments (e.g., using the Helmholtz decomposition theorem), a fact we shall use when considering the vorticity formulation of Navier-Stokes equations in the proof of Proposition \ref{PropNSBound} below.

\begin{Prop}\label{PropNSExist}
	Fix $T>0$, $\nu>0$, and $f\in \Hnull$.
    \begin{enumerate}
        \item[(i)] For all $\theta\in \divH^1$, the 2D periodic Navier-Stokes equations (\ref{eq:NS}) with initial condition $\theta$ admit a unique strong solution 
	$$
	u_{\theta}\in C([0,T], \divH^1)\cap L^2([0,T], \divH^2)
	,\quad 
	\textrm{with } \quad
	\frac{du_\theta}{dt}\in L^2([0,T], \Hnull)
	.
	$$
	In addition, for any $L>0$ there exists a constant $c=c(T,\nu, \|f\|_{L^2},L)$ such that 
	$$
	\sup_{\substack{\theta\in \divH^1\\ \|\theta\|_{\dH^1}\leq L}} \bigg( \sup_{0\leq t\leq T} \|u_{\theta}(t)\|_{\dH^1} + \int_0^T \|u_{\theta}(t)\|_{\dH^2}\, dt 
	\bigg)
	\leq 
	c
	.
	$$

        \item[(ii)] For any $L>0$, there exists a constant $K=K(T,\nu, \|f\|_{L^2}, L)>0$ such that for all~$\theta\in \divH^1$ with $\|\theta\|_{\dH^1}\le L$ and all $\theta'\in \divH^1$ we have
	$$
	\sup_{0\leq t\leq T} \|u_{\theta}(t)-u_{\theta'}(t)\|_{L^2}
	\leq
	K \|\theta-\theta'\|_{L^2}
    .
	$$
    \end{enumerate}
\end{Prop}

\subsubsection{Sobolev bounds on the solution}

In the sequel, we will repeatedly use the following common version of Young's inequality: for all $a,b\geq 0$ and $\ve>0$, we have 
\begin{equation}\label{eq:Young}
ab
\leq 
\frac{1}{4\ve}a^2 + \ve b^2
.
\end{equation}

\begin{Prop}\label{PropNSBound}
    Fix an integer $a\geq 1$ and $f\in \divH^{a-1}$. Then, for all $\theta\in \divH^a$ we have 
        $$
        u_{\theta}\in C([0,T],\divH^a)\cap L^2([0,T], \divH^{a+1})
        .
        $$
        In addition, there exists $C=C(T,\nu,a,\|f\|_{\dH^{a-1}})>0$ such that for all $\theta\in \divH^a$, we have
        $$
        \sup_{0\leq t\leq T} \|u_{\theta}(t)\|^2_{\dH^a} 
        +
        \int_0^T \|u_{\theta}(t)\|^2_{\dH^{a+1}}\, dt 
        \leq 
        C\big(1+ \|\theta\|^{2a}_{\dH^a}\big)
        .
        $$
\end{Prop}
\vspace{1mm}



\begin{Proof}{Proposition \ref{PropNSBound}}
    We prove the result by a Galerkin approximation argument: for each $n$, solve the system of $n$ ODE's for $u_n=\sum_{j=1}^n b_j e_j\in \Hnull$ satisfying
    $$
    \frac{du_n}{dt}-\nu\Delta u_n + P_n[(u_n\cdot\nabla)u_n] = P_n[f]
    ,\spa 
    u_n(0)=P_n[\theta] 
    ,
    $$
    and satisfying uniform-in-$n$ estimates, where $P_n$ is the projection onto the first $n$ eigenfunctions $e_1,\ldots, e_n\in \divH^\infty$ of $A=-\Delta$. In particular, $u_n\in\divH^\infty$ for all $n$. In what follows, the estimates are obtained for such (smooth) $u_n$ and a rigorous limit procedure must be performed to show that the estimates hold for $u$. This will be the case for us since all estimates are obtained in $L^2$, and the limiting argument can be done as in Section 9~in \cite{CF88}. In what follows, we thus assume, without loss of generality, that $u(t)\in \divH^\infty$ for all $t>0$ and that $t\mapsto u(t)$ is continuous from $[0,T]$ to $(\divH^a, \|\cdot\|_{\dH^a})$ for all $a\ge 1$.

    We will prove the result by hand for $a=1,2$ and the conclusion will follow for $a\geq 3$ by induction. To ease notation, we will omit the dependence of $u_{\theta}$ on $\theta$ and simply write $u$ to denote $u_\theta$. Recall that $u$ satisfies 
    \begin{equation}\label{eq:EDPNS}
    \frac{du}{dt}
    -
    \nu \Delta u
    +
    B[u,u]
    =
    f
    ,\spa 
    u(0)=\theta
    .
    \end{equation}
    {\bf The case $a=1$.} Taking the $\dH^1$-inner product of (\ref{eq:EDPNS}) with $u$ yields 
    $$
    \frac12 \frac{d}{dt}\|u(t)\|^2_{\dH^1}
    +
    \nu \| u(t)\|^2_{\dH^2}
    =
    \ps{f(t)}{u(t)}_{\dH^1}
    ,
    $$
    where we used (\ref{IPPdt}), the integration by parts formula (\ref{IPPSob2}), and Proposition \ref{PropBilinearFormIdentities1}(ii). By~(\ref{eq:dualProd2}) and Young's inequality (\ref{eq:Young}), we have 
    $$
    \ps{f(t)}{u(t)}_{\dH^1}
    \leq 
    \|f(t)\|_{L^2}\|u(t)\|_{\dH^2}
    \leq 
    \frac{1}{2\nu}\|f(t)\|^2_{L^2} 
    +
    \frac{\nu}{2}\|u(t)\|^2_{\dH^2}
    .
    $$
    It follows that 
    $$
    \frac12 \frac{d}{dt}\| u(t)\|^2_{\dH^1}
    +
    \frac{\nu}{2} \| u(t)\|^2_{\dH^2}
    \leq 
    \frac{1}{2\nu}\|f(t)\|^2_{L^2} 
    .
    $$
    Integrating the last display in time, and using that $ u(t)\to \theta$ in $\dH^1$ as $t\downarrow 0$, provides 
    $$
    \| u(t)\|^2_{\dH^1}
    +
    \nu \int_0^t \|u(s)\|^2_{\dH^2}\, ds
    \leq 
    \| \theta\|^2_{\dH^1} 
    +
    \frac{1}{\nu}\int_0^t \|f(s)\|^2_{L^2}\, ds
    \leq 
    \| \theta\|^2_{\dH^1} 
    +
    \frac{1}{\nu}\int_0^T \|f(s)\|^2_{L^2}\, ds
    .
    $$
    We deduce that
    $$
    \sup_{0\leq t\leq T} \|u(t)\|^2_{\dH^1}
    +
    \int_0^T \|u(t)\|^2_{\dH^2}\, dt
    \lesssim 
    1+\|\theta\|^2_{\dH^1}
    .
    $$
    {\bf The case $a=2$.} We follow an argument in Section 3.5 in \cite{NicTit2024}. Letting $\nabla^\perp = (-\partial_2, \partial_1)$ denote the skew-gradient, define the corresponding `skew-divergence' $\omega:\R^2\to\R$ of (the smooth map) $u=(u_1,u_2)$ as $\omega\equiv\nabla^\perp\cdot u= -\partial_2 u_1 + \partial_1 u_2$. Using that $\nabla\cdot u = 0$ since $u\in \Hnull$, we deduce that $\partial_1 \partial_2 u_1 = -\partial_2^2 u_2$ and $\partial_1 \partial_2 u_2 = -\partial^2_1 u_1$. These identities lead to 
    $$
    \nabla \omega 
    =
    (\Delta u_2, -\Delta u_1)
    ,
    $$
    from which we deduce that 
    $$
    \ps{\nabla \omega}{h_{k_j}}_{L^2} 
    = 
    \big(\ps{\Delta u_2}{h_{k_j}}_{L^2}, - \ps{\Delta u_1}{h_{k_j}}_{L^2}\big)
    ,\spa 
    \forall\ j\geq 0 
    .
    $$
    In particular, we have $\lvert\ps{\nabla \omega}{h_{k_j}}_{L^2}\rvert=\lvert\ps{\Delta u}{h_{k_j}}_{L^2}\rvert$ for all $j\geq 1$, so that
    $$
    \|\nabla \omega\|_{\dH^\beta}
    =
    \|\Delta u\|_{\dH^\beta}
    =
    \|u\|_{\dH^{\beta+2}}
    ,\spa 
    \forall\ \beta\in\R
    ,
    $$
    by virtue of (\ref{IppSob1}).  Since $\nabla^\perp \cdot (\nabla v)=0$ for any smooth $v:\Om\to\R$, straightforward computations further provide 
    \begin{align*}
    \nabla^\perp\cdot \big( (u\cdot\nabla)u \big) 
    &=
    -(u_1 \partial_1\partial_2 u_1 + u_2 \partial^2_2 u_1) 
    +
    (u_1 \partial^2_1 u_2 + u_2 \partial_1 \partial_2 u_2)
    \\[2mm]
    &=
    u\cdot \big(\nabla(-\partial_2 u_1)\big)
    +
    u\cdot \big(\nabla(\partial_1 u_2)\big)
    \\[2mm]
    &=
    u\cdot\nabla\omega
    .
    \end{align*}
    Applying $\nabla^\perp\cdot$ to the original system of Navier-Stokes equations (\ref{eq:NSequations}) thus leads to the so-called `vorticity' formulation 
    $$
    \frac{\partial \omega}{\partial t} - \nu \Delta \omega 
    +  u\cdot\nabla\omega = (\nabla^\perp\cdot f)
    .
    $$
    Taking the $L^2(\Om)$-inner product of the last display with $-\Delta \omega$ provides 
    $$
    \frac12 \frac{d}{dt}\|\nabla\omega\|^2_{L^2}
    +
    \nu \|\Delta \omega\|^2_{L^2}
    -
    \ps{u\cdot\nabla\omega}{\Delta \omega}_{L^2} 
    =
    -\langle \nabla^\perp\cdot f, \Delta \omega\rangle_{L^2}
    ,
    $$
    where we used (\ref{IPPdt}). Cauchy-Schwarz inequality and Young's inequality (\ref{eq:Young}) provide
    $$
    \lvert\langle \nabla^\perp\cdot f, \Delta \omega\rangle_{L^2}\rvert
    \lesssim
    \|\nabla f\|_{L^2}\|\Delta \omega\|_{L^2}
    \leq 
    \frac{1}{\nu} \|\nabla f\|^2_{L^2}
    +
    \frac{\nu}{4}\|\Delta \omega\|^2_{L^2}
    .
    $$
    An integration by parts yields 
    $$
    \ps{u\cdot\nabla\omega}{\Delta \omega}_{L^2}
        =
        \sum_{i=1}^2 \int_{\Om} (u\cdot \nabla\omega)\partial_i^2 \omega
        =
        - \sum_{i=1}^2 \int_{\Om} \partial_i (u\cdot \nabla\omega)\partial_i \omega
    .
    $$
    Since 
    $
    \partial_i (u\cdot \nabla\omega)
    =
    \partial_i u\cdot \nabla\omega
    + 
    u\cdot \nabla\partial_i\omega
    $
    and $\int_{\Om} (u\cdot\nabla\partial_i \omega)\partial_i\omega = 0$ as in Proposition \ref{PropBilinearFormIdentities1}(i),
    we deduce that 
    $$
    \lvert\ps{u\cdot\nabla\omega}{\Delta \omega}_{L^2}\rvert
    =
    \Big| 
    \sum_{i=1}^2 \int_{\Om} (\partial_i u\cdot \nabla\omega)\partial_i \omega 
    \Big|
    \lesssim
    \|\nabla u\|_{L^2} \|\nabla \omega\|^2_{L^4}  
    .
    $$
    \lad's inequality (\ref{LadInequality}) thus entails that, for some $c>0$, we have
    $$
    \lvert\ps{u\cdot\nabla\omega}{\Delta \omega}_{L^2}\rvert
    \leq 
    c \|\nabla u\|_{L^2} \|\nabla \omega\|_{L^2}\|\Delta \omega\|_{L^2}
    =
    c \| u\|_{\dH^1} \|\nabla \omega\|_{L^2}\|\Delta \omega\|_{L^2}
    ,
    $$
    where we used (\ref{IPPSob3}) in the last equality.
    Consequently, Young's inequality (\ref{eq:Young}) provides
    $$
    \frac12 \frac{d}{dt}\|\nabla\omega\|^2_{L^2}
    +
    \nu \|\Delta \omega\|^2_{L^2}
    \leq 
    \frac{1}{\nu} \|\nabla f\|^2_{L^2}
    +
    \frac{\nu}{4}\|\Delta \omega\|^2_{L^2}
    +
    \frac{c^2}{\nu} \| u\|^2_{\dH^1}\|\nabla \omega\|^2_{L^2}
    +
   \frac{\nu}{4} \|\Delta \omega\|^2_{L^2}
   ,
    $$
    from which we deduce that 
    $$
    \frac{d}{dt}\|\nabla\omega\|^2_{L^2}
    +
    \nu \|\Delta \omega\|^2_{L^2}
    \leq 
    \frac{2}{\nu} \|\nabla f\|^2_{L^2}
    +
    \frac{2c^2}{\nu} \| u\|^2_{\dH^1}\|\nabla \omega\|^2_{L^2}
    .
    $$
    Recalling that $\sup_{0\leq t\leq T} \|u(t)\|^2_{\dH^1}\lesssim 1 + \|\theta\|^2_{\dH^1}$ from the case $a=1$, and $\|\nabla \omega\|_{L^2} = \| u\|_{\dH^2}$, we find
    $$
    \frac{d}{dt}\| u\|^2_{\dH^2}
    +
    \nu \|\Delta \omega\|^2_{L^2}
    \lesssim 
    1
    +
    (1+\|\theta\|^2_{\dH^1})\|u\|^2_{\dH^2}
    .
    $$
    Integrating the last display in time, and using that $u(t)\to \theta$ in $\dH^2$ as $t\downarrow 0$, provides the (uniform in $t\in [0,T]$) estimate
    $$
    \|u(t)\|^2_{\dH^2}
    +
    \nu \int_0^t \|\Delta \omega(s)\|^2_{L^2}\, ds
    \lesssim 
    1+\|\theta\|^2_{\dH^2} + (1+\|\theta\|^2_{\dH^1}) \int_0^t \|u(s)\|^2_{\dH^2}\, ds
    .
    $$
    Using the estimate $\int_0^T \| u(s)\|^2_{\dH^2}\, ds\lesssim 1+\|\theta\|^2_{\dH^1}$ established when $a=1$, and the fact that 
    $$
    \|u\|_{\dH^3}= \|\Delta u\|_{\dH^1} = \|\nabla \omega\|_{\dH^1}\lesssim \|\Delta \omega\|_{L^2}
    ,
    $$
    where the first equality follows from (\ref{IppSob1}) and the third (inequality) from Poincaré's inequality, we find 
    $$
    \sup_{0\leq t\leq T}
    \| u(t)\|^2_{\dH^2}
    +
    \int_0^T \|u(t)\|^2_{\dH^3}\, dt 
    \lesssim 
    1 
    +
    \|\theta\|^4_{\dH^1}
    +
    \|\theta\|^2_{\dH^2}
    .
    $$
    {\bf The case $a\geq 3$.} We prove this final case by induction. Since the conclusion holds for $a=2$, let us assume that for some $a\geq 2$ we have 
    \begin{equation}\label{eq:InductionEstimateNS}
    \sup_{0\leq t\leq T}
    \| u(t)\|^2_{\dH^a}
    +
    \int_0^T \|u(t)\|^2_{\dH^{a+1}}\, dt 
    \lesssim 
    1 
    +
    Q_a(\theta)
    .
    \end{equation}
    where $Q_a$ is a real-valued map with $Q_1(\theta)=\|\theta\|^2_{\dH^1}$ and $Q_2(\theta)=\|\theta\|^4_{\dH^1} + \|\theta\|^2_{\dH^2}$ as we  established when $a=1$ and $a=2$, respectively. Taking the $\dH^{a+1}$-inner product of (\ref{eq:EDPNS}) with $u$ yields,
    $$
    \frac12 \frac{d}{dt}\|u\|^2_{\dH^{a+1}} 
    +
    \nu \|u\|^2_{\dH^{a+2}} 
    +
    \ps{B[u,u]}{u}_{\dH^{a+1}}
    =
    \ps{f}{u}_{\dH^{a+1}}
    .
    $$
    where we used (\ref{IPPdt}) and the integration by parts formula (\ref{IPPSob2}). First observe that Cauchy-Schwarz inequality and Young's inequality (\ref{eq:Young}) entail that  
    $$
    \ps{f}{u}_{\dH^{a+1}}
    \leq 
    \|f\|_{\dH^{a}}\|u\|_{\dH^{a+2}}
    \leq 
    \nu \|f\|^2_{\dH^{a}} + \frac{\nu}{4}\| u\|^2_{\dH^{a+2}}
    .
    $$
    Observe now that the projection bound (\ref{eq:BoundBSobolev}) yields  
    $$
    \lvert\ps{B[u,u]}{u}_{\dH^{a+1}}\rvert
    \leq 
    \|B[u,u]\|_{\dH^{a}}\|u\|_{\dH^{a+2}}
    \leq
    \|(u\cdot\nabla)u\|_{\dH^{a}} \| u\|_{\dH^{a+2}}
    .
    $$
    Because $a\geq 2$, the multiplier inequality for Sobolev norms (\ref{multipSob}) yields 
    $$
    \|(u\cdot\nabla)u\|_{\dH^{a}}
    \leq 
    c_1 \|u\|_{\dH^a}\|\nabla u\|_{\dH^a}
    = 
    c_1 \|u\|_{\dH^a}\| u\|_{\dH^{a+1}}
    ,
    $$
    by virtue of (\ref{IPPSob3}).
    It follows from Young's inequality (\ref{eq:Young}) that 
    $$
    |\ps{B[u,u]}{u}_{\dH^{a+1}}|
    \leq 
    c_2 \|u\|^2_{\dH^a}\|u\|^2_{\dH^{a+1}} 
    + 
    \frac{\nu}{4}\| u\|^2_{\dH^{a+2}}
    .
    $$
    We deduce that 
    $$
    \frac{d}{dt}\|u\|^2_{\dH^{a+1}} 
    +
    \nu \| u\|^2_{\dH^{a+2}} 
    \lesssim 
    \|f\|^2_{\dH^a}
    +
    \|u\|^2_{\dH^a}\|u\|^2_{\dH^{a+1}}
    .
    $$
    Recalling the induction estimate $\sup_{0\leq t\leq T} \|u(t)\|^2_{\dH^a} \lesssim 1+Q_a(\theta)$, we have
    $$
    \frac{d}{dt}\|u\|^2_{\dH^{a+1}} 
    +
    \nu \| u\|^2_{\dH^{a+2}} 
    \lesssim 
    \|f\|^2_{\dH^a}
    +
    (1+Q_a(\theta)) \|u\|^2_{\dH^{a+1}}
    .
    $$
    Integrating the last display in time and using that $u(t)\to \theta$ in $\dH^{a+1}$ at $t\downarrow 0$ yields the (uniform in $t\in [0,T]$) estimate
    $$
    \|u(t)\|^2_{\dH^{a+1}} + \nu \int_0^t \| u(s)\|^2_{\dH^{a+2}}\, ds
    \lesssim 
    1+\|\theta\|^2_{\dH^{a+1}} 
    +
    (1+Q_a(\theta))\int_0^t \|u(s)\|^2_{\dH^{a+1}}\, ds
    .
    $$
    The induction estimate $\int_0^T \|u(s)\|^2_{\dH^{a+1}}\, ds\leq 1+Q_a(\theta)$ thus provides
    \begin{align*}
    \sup_{0\leq t\leq T} \|u(t)\|^2_{\dH^{a+1}} + \nu \int_0^T \| u(s)\|^2_{\dH^{a+2}}\, ds
    &\lesssim 
    1+\|\theta\|^2_{\dH^{a+1}} 
    +
    (1+Q_a(\theta))^2
    \\
    &\lesssim 
    1+\|\theta\|^2_{\dH^{a+1}} 
    +
    Q_a(\theta)^2
    .
    \end{align*}
    By induction, we deduce that (\ref{eq:InductionEstimateNS}) holds for all $a\geq 2$ with $Q_a$ defined recursively through $Q_{a+1}(\theta)=Q_a(\theta)^2 + \|\theta\|^2_{\dH^{a+1}}$. Inspection of the case $a=2$ shows that we also have $Q_2(\theta)=Q_1(\theta)^2+\|\theta\|^2_{\dH^2}$ with $Q_1(\theta)=\|\theta\|^2_{\dH^1}$. We thus have for all $a\geq 1$
    $$
    Q_{a+1}(\theta)
    =
    Q_a(\theta)^2 + \|\theta\|^2_{\dH^{a+1}}
    ,\spa 
    Q_1(\theta) = \|\theta\|^2_{\dH^1}
    .
    $$
    Straightforward computations provide
    $$
    Q_a(\theta)
    \lesssim 
    \sum_{\ell=0}^{a-1} \|\theta\|^{2(\ell+1)}_{\dH^{a-\ell}}
    ,\spa 
    \forall\ a\geq 1
    $$
    Using Poincaré's inequality (\ref{Poincare}), $\|\theta\|_{\dH^\beta}\lesssim \|\theta\|_{\dH^\gamma}$, $0\leq \beta\leq \gamma$, we deduce that for all $a\geq 1$
    $$
    Q_a(\theta)
    \lesssim 
    1
    + \|\theta\|^{2a}_{\dH^a}
    ,
    $$
    which concludes this part of the proof.

    It remains to show the true solution $u:[0,T]\to\divH^a$, as opposed to its Galerkin approximation $u_n$ considered so far in the present proof, is not only a bounded map but also continuous over $[0,T]$. This will follow from Theorem 4 of Section 5.9 in \cite{Evans1998} if we show that $du/dt\in L^2([0,T], \divH^{a-1})$. It is thus enough to establish a uniform-in-$n$ estimate on $\int_0^T \|du_n/dt\|^2_{\dH^{a-1}}\, dt$. For this purpose, still writing $u$ for the Galerkin approximation by abuse of notation, observe that 
    \begin{align}\label{eq:dudt}
    \Big\|\frac{du}{dt}\Big\|_{\dH^{a-1}}
    &\leq 
    \nu \| u\|_{\dH^{a+1}} 
    +
    \|B[u, u]\|_{\dH^{a-1}} 
    +
    \|f\|_{\dH^{a-1}}
    ,
    \end{align}
    where we used (\ref{IPPSob3}). 
    When $a=1$, the projection bound (\ref{eq:BoundBSobolev}), the multiplier inequality for Sobolev norms (\ref{multipSob}), and \lad's inequality (\ref{LadInequality}) yield
    \begin{eqnarray*}
    \lefteqn{
    \|B[u,u]\|_{L^2}
    \leq
    \|(\nabla u)u\|_{L^2}
    \lesssim
    \|u\|^{1/2}_{L^2}\|\nabla u\|_{L^2}\|\Delta u\|^{1/2}_{L^2}
    }
    \\[2mm]
    &&
    \hspace{10mm}
    \lesssim 
    \|u\|^{3/2}_{\dH^1}\|u\|^{1/2}_{\dH^2}
    =
    \|u\|^{3/2}_{\dH^a}\|u\|^{1/2}_{\dH^{a+1}}
    \leq 
    \frac12 \|u\|^3_{\dH^a} + \frac12 \|u\|_{\dH^{a+1}}
    ,
    \end{eqnarray*}
    where we used Poincaré's inequality and (\ref{IPPSob3}) in the third inequality, and Young's inequality (\ref{eq:Young}) in the last one. When $a\geq 2$, then Lemma \ref{PropBilinearFormIdentities2}(i) and the multiplier inequality for Sobolev norms (\ref{multipSob}) provide 
    $$
    \|B[u,u]\|_{\dH^{a-1}}
    \leq 
    \|u u^T\|_{\dH^a}
    \lesssim 
    \|u\|^2_{\dH^a}
    .
    $$
    For all $a\geq 1$, using the (`uniform-in-$n$') estimate 
    $$
    \sup_{0\leq t\leq T} \|u(t)\|^2_{\dH^a} 
    +
    \int_0^T \|u(t)\|^2_{\dH^{a+1}}\, dt
    \lesssim 
    1+\|\theta\|^{2a}_{\dH^a}
    $$ 
    established above, then
    squaring (\ref{eq:dudt}) and integrating in time between $0$ and $T$ yields the `uniform-in-$n$' estimate
    $$
    \int_0^T \Big\| \frac{du}{dt} \Big\|^2_{\dH^{a-1}}\, dt 
    \lesssim 
    1+\int_0^T \|u(t)\|^2_{\dH^{a+1}}\, dt 
    \lesssim 
    1+\|\theta\|^{6a}_{\dH^a}
    .
    $$
    This concludes the proof.
\end{Proof}

\subsubsection{Lipschitz continuity in Sobolev norm}\label{sec:Lip}

We now use Proposition \ref{PropNSBound} to extend the Lipschitz result of Proposition \ref{PropNSExist}(ii) to Lipschitz continuity in any (non-negative) Sobolev norm by adapting the proof of Theorem 9.4 in \cite{Robinson2001}.

\begin{Prop}\label{PropNSLip}
    Fix $a\in\Z$ and $a^*$ as in (\ref{eq:DefinAStar}). Then for any $f\in\divH^{a^*-1}$ and $L>0$, there exists a constant $C=C(T,\nu,\|f\|_{\dH^{a^*-1}},a,L)>0$ such that for all $\theta,\theta'\in \divH^{a^*}$ with $\|\theta\|_{\dH^{a^*}} + \|\theta'\|_{\dH^{a^*}}\le L$,  we have 
        $$
        \sup_{0\leq t\leq T} 
        \|u_{\theta}(t)-u_{\theta'}(t)\|^2_{\dH^a}
        +
        \int_0^T \|u_\theta(t)-u_{\theta'}(t)\|^2_{\dH^{a+1}}\, dt
        \leq 
        K\|\theta-\theta'\|^2_{\dH^a}
        .
        $$
\end{Prop}



\begin{Proof}{Proposition \ref{PropNSLip}}
    Let $w=u_{\theta}-u_{\theta'}$ and observe that $w$ satisfies the equation 
    $$
    \frac{dw}{dt} - \nu \Delta w + B[w, u_{\theta}] + B[u_{\theta'},w]
    =
    0
    ,\spa 
    w(0)=\theta-\theta'
    .
    $$
    Then Proposition \ref{PropLinPSob} with $g=0$ and $\xi = \theta-\theta'$ yields the conclusion.
\end{Proof}

\subsection{Linearization of the Navier-Stokes flow}\label{sec:Lin}

\subsubsection{Existence}

The following Proposition can be found as Theorem 13.20 in \cite{Robinson2001}. It was proved for initial conditions $\theta_0,\theta\in \AA$, where $\AA$ is the global attractor (in $\Hnull$ or $V$, see Lemma 12.6 in the same monograph) of the 2D Navier-Stokes equations. Inspection of the proof, however, shows that these assumptions are only used to ensure that the solutions $u_{\theta_0}$ and $u_{\theta}$ are uniformly bounded in $L^\infty([0,T],\divH^1)$, which holds as soon as we have $f\in \divH$ and $\theta_0$ and $\theta$ belong to a fixed $\divH^1$-ball (Proposition \ref{PropNSExist}(i)). Moreover, note that $\Lambda_{\theta_0,t}$ below necessarily maps into the closed subspace $\Hnull$ of $L^2$ and, hence, in Theorem 13.20 in \cite{Robinson2001} we have $A=-\Delta$.

\begin{Prop}\label{PropExistLin}
	Fix $\theta_0\in \divH^1$ and $f\in \Hnull$, and let $t>0$. 
    \begin{enumerate}
        \item[(i)] There exists a compact linear operator $\Lambda_{\theta_0,t} : \divH^1 \to \Hnull$ such that 
	$$
	\sup_{\substack{\theta\in \divH^1\\ \|\theta-\theta_0\|_{L^2}<\ve}}
	\frac{\|u_{\theta}(t) - u_{\theta_0}(t)-\Lambda_{\theta_0,t}[\theta-\theta_0] \|_{L^2}}{\|\theta-\theta_0\|_{L^2}}
	\to
	0
	,\spa 
	\ve\to 
	0
	.
	$$

        \item[(ii)] For any $\xi\in \divH^1$, the map $U:[0,T]\to \Hnull, t\mapsto\Lambda_{\theta_0,t}[\xi]$ solves the linear equation in $\divH^{-1}$ 
	\begin{equation}\label{LinPDE}
	\frac{dU}{dt}
	-\nu \Delta U
	+ B[u_{\theta_0}, U] 
	+ B[U, u_{\theta_0}]
	=
	0
	,\quad 
	U(0) = \xi 
	.
	\end{equation}

        \item[(iii)] For any $T>0$ and $L>0$ there exists a constant $c=c(\nu, \|f\|_{\Hnull}, T,L)>0$ such that for all $\theta_0\in \divH^1$ with $\|\theta_0\|_{\dH^1}\le L$ and all $\theta\in \divH^1$, we have
	$$
	\sup_{0\leq t\leq T} 
	\|u_{\theta}(t) - u_{\theta_0}(t)-\Lambda_{\theta_0,t}[\theta-\theta_0] \|_{L^2} 
	\leq 
	c\|\theta-\theta_0\|^2_{L^2}
	.
	$$
    \end{enumerate}
\end{Prop}

We will show in Proposition \ref{PropLinSob} below that the operator $\Lambda_{\theta_0,t}:\divH^1\to \Hnull$ extends to a bounded operator on $\Hnull$, which we abusively still denote by $\Lambda_{\theta_0,t}:\Hnull\to \Hnull$, with a uniform (in $t\in [0,T]$) control on the operator of $\Lambda_{\theta_0,t}$. As a consequence, for all $\theta_0\in \divH^1$, we will form a bounded operator 
$$
\I_{\theta_0}:\Hnull\to C([0,T],\Hnull)
$$ 
by letting, for any $\xi\in \Hnull$ and $0\leq t\leq T$
\begin{equation}\label{defI}
\I_{\theta_0}[\xi](t,\cdot) 
\equiv
\Lambda_{\theta_0,t}[\xi](\cdot)
.
\end{equation}
By studying solutions of (\ref{LinPDE}), we will further show that $\I_{\theta_0}$ extends to a bounded operator between $\divH^a$ and $C([0,T],\divH^a)\cap L^2([0,T],\divH^{a+1})$ for any $a\in\Z$. Proposition~\ref{PropExistLin}(i) entails that $\Lambda_{\theta_0,t}:\divH^1\to\Hnull$ is compact for any $\theta_0\in \divH^1$ and $t > 0$, and we will strengthen this result by showing that $\Lambda_{\theta_0,t}$ is in fact compact on $\divH^a$ for all $a\in\Z$: we will establish that $\I_{\theta_0}$ maps any Sobolev space $\divH^a$, $a\in\Z$, into $\cap_{b\in\Z} L^\infty([t_0,T],\divH^b)$ for any $t_0>0$ fixed, thus reflecting the smoothing nature of the associated parabolic flow; see Proposition \ref{PropLinSmooth}. The smoothing property of the linearized flow will play a crucial role to prove a functional Bernstein-von Mises theorem at the forward level in a strong (uniform-type) topology in Theorem \ref{nstokbvm}.

\subsubsection{Sobolev estimates}


The proofs of Proposition \ref{PropLinSob}, Proposition \ref{PropLinStab}, and Proposition \ref{PropLinSmooth} below  follow, respectively, from Proposition \ref{PropLinPSob}, Proposition \ref{PropLinPStab}, and Proposition \ref{PropLinPSmooth}, by noticing that $U\equiv \I_\theta[\xi]$ satisfies (\ref{LinNSsource}) with $g=0$ and $\theta'=\theta$, by virtue of Proposition \ref{PropExistLin}(ii).

\begin{Prop}\label{PropLinSob}
    Fix $a\in\Z$ and $a^*$ as in (\ref{eq:DefinAStar}).  Then, for any $\theta\in \divH^{a^*}$ and $f\in\divH^{a^*-1}$, the operator $\I_{\theta}$ from (\ref{defI}) extends to $\divH^a$ and we have for all $\xi\in \divH^a$
    $$
    \I_{\theta}[\xi]\in C([0,T], \divH^a)\cap L^2([0,T], \divH^{a+1})
    $$
    In addition, for all $f\in\divH^{a^*-1}$ and $L>0$ there exists $C=C(\nu,T,\|f\|_{\dH^{a^*-1}},a,L)>0$ such that for all $\theta\in\divH^{a^*}$ with $\|\theta\|_{\dH^{a^*}}\le L$ and $\xi\in \divH^a$, we have 
    $$
    \sup_{0\leq t\leq T} \|\I_{\theta}[\xi](t)\|^2_{\dH^a}
    +    
    \int_0^T  \|\I_{\theta}[\xi](t)\|^2_{\dH^{a+1}}\, dt 
    \leq  
    C\|\xi\|^2_{\dH^a} 
    $$
\end{Prop}

\vspace{1mm}



\begin{Prop}\label{PropLinStab}
    Fix $a\in\Z$ and $a^*$ as in (\ref{eq:DefinAStar}). Then, for all $f\in\divH^{a^*-1}$ and $L>0$ there exists $C=C(\nu, T, \|f\|_{\dH^{a^*-1}},a, L)>0$ such that for all $\theta\in\divH^{a^*}$ with $\|\theta\|_{\dH^{a^*}}\le L$ and $\xi\in \divH^a$, we have
    $$
    \|\xi\|^2_{\dH^{a}}
    \le 
    C\int_0^T \|\I_\theta[\xi](t)\|^2_{\dH^{a+1}}\, dt
    .
    $$
\end{Prop}

\vspace{1mm}



\begin{Prop}\label{PropLinSmooth}
    Fix $a\in\Z$, $b\in\Z$ such that $b\ge a+1$, and $b^*$ as in (\ref{eq:DefinAStar}). Then for any $t_{\min}\in (0,T)$, $L>0$, and $f\in\divH^{b^*-1}$, there exists $C=C(\nu,T,\|f\|_{\dH^{b^*-1}}, L, t_{\min})>0$ such that for all $\theta\in \divH^{b^*}$ with $\|\theta\|_{\dH^{b^*}}\le L$ and $\xi\in\divH^a$, we have
    $$
    \sup_{t\in [t_{\min}, T]} \|\I_\theta[\xi](t)\|_{\dH^b}
    \leq
    C
    \|\xi\|^2_{\dH^a}
    .
    $$
\end{Prop}

\vspace{3mm}

We now turn to an extension of Proposition \ref{PropExistLin}(iii) to arbitrary Sobolev norms.

\begin{Prop}\label{PropLinQuad}
    Fix $a\in\Z$ with $a\geq 2$. Then, for all $f\in\divH^{a-1}$ and $L>0$ there exists $C=C(\nu, T,\|f\|_{\dH^{a-1}}, a, L)>0$ such that for all $\theta,\theta'\in\divH^a$ with $\|\theta\|_{\dH^a}+\|\theta'\|_{\dH^a}\le L$ we have
    $$
    \sup_{0\leq t\leq T} \|u_{\theta'}(t)-u_{\theta}(t)-\I_{\theta}[\theta'-\theta](t)\|_{\dH^a} 
    \le 
    C\|\theta'-\theta\|^2_{\dH^a}
    .
    $$
\end{Prop}

\vspace{1mm}

\begin{Proof}{Proposition \ref{PropLinQuad}}
    Let $R=u_{\theta'}-u_\theta-\I_{\theta}[\theta'-\theta]$ and $w=u_{\theta'}-u_\theta$. Then one can check that $R$ satisfies 
    $$
    \frac{dR}{dt}-\nu\Delta R + B[u_{\theta},R] + B[R,u_{\theta}] = B[w,w]
    , 
    \quad 
    R(0)=0
    ;
    $$
    see, e.g., the proof of Theorem 13.20 in \cite{Robinson2001}. Proposition \ref{PropLinPSob} entails that
    $$
    \sup_{0\leq t\leq T}
    \|R(t)\|^2_{\dH^a}
    \leq 
    c \int_0^T \|B[w,w](s)\|^2_{\dH^{a-1}}\, ds 
    ,
    $$
    for some constant $c=c(\nu, T,\|f\|_{\dH^{a-1}}, a, L)>0$. 
    Since $a-1\ge 1$, Proposition \ref{PropProdSob} yields 
    $$
    \|B[w,w]\|^2_{\dH^{a-1}}
    \lesssim 
    \|w\|^4_{\dH^{a}}
    .
    $$
   Since $a\geq 2$, we have $a^*=a$ (see (\ref{eq:DefinAStar})), so that
    Proposition \ref{PropNSLip} provides 
    $$
    \sup_{0\leq s\leq T}\|w(s)\|_{\dH^{a}}
    \lesssim 
    \|\theta-\theta'\|_{\dH^{a}}
    ,
    $$
    uniformly in $\|\theta\|_{\dH^a}+\|\theta'\|_{\dH^a}\le L$. It follows that 
    $$
    \sup_{0\leq t\leq T}
    \|R(t)\|^2_{\dH^a}
    \leq 
    C \|\theta-\theta'\|^4_{\dH^a}
    ,
    $$
    holds for some constant $C=C(\nu, T, \|f\|_{\dH^{a-1}}, a, L)$, and all $\theta,\theta'\in\dH^a$ with $\|\theta\|_{\dH^a}+\|\theta'\|_{\dH^a}\le L$. This concludes the proof.
\end{Proof}

\section{Linear parabolic estimates}

In this section, we provide Sobolev estimates for the solution of linear parabolic equations that will be of importance to establish the results of Section \ref{sec:InvertI}, and the Sobolev estimates of Section \ref{sec:Lip} and Section \ref{sec:Lin}.   

\subsection{Estimates for the periodic heat equation}\label{sec:ParaEstim}

In this section we recall some standard facts about solutions $u:[0,T]\times\Om\to\R^2$ to the PDE
\begin{align}\label{heateq}
    \frac{du}{dt} - \nu \Delta u &= g \spa \text{on } (0,T]\times\Om
    \\[2mm]\nonumber
    u(0)&=h \spa \text{on } \Om 
    ,
\end{align}
that further ensure that $u(t)$ lies in $ \Hnull$ when the initial condition $h$ and the heat source $g$ do so. The solution can be expressed in the orthonormal basis $\{e_j : j\geq 1\}$ from (\ref{defej}), consisting of eigenfunctions of the periodic negative Laplacian with corresponding eigenvalues $\lambda_j\ge 0$, as
\begin{equation}\label{repheat}
u(t) 
=
\sum_{j\geq 1} e^{-\nu\lambda_j t}\Big( h_j + \int_0^t g_j(s) e^{\nu\lambda_j s}\, ds\Big)e_j
,
\end{equation}
when $g(t)=\sum_{j\geq 1} g_j(t) e_j$ and $h = \sum_{j\geq 1} h_j e_j$.
Arguing as in the proof of Proposition \ref{PropLinPSob} below (but setting all terms involving $B$ to zero) one shows that for all $a\in\Z$, $g\in L^2([0,T],\divH^{a-2})$, and $h\in\divH^{a-1}$, the solution $u=u_{g,h}$ to (\ref{heateq}) satisfies
\begin{equation}\label{heatestimate}
\|u\|_{L^2([0,T],\divH^a)}
\leq 
C\big( \|g\|_{L^2([0,T],\dH^{a-2})} + \|h\|_{\dH^{a-1}} \big)
,
\end{equation}
for some constant $C=C(T,a)>0$.

\subsection{Regularity theory for a class of linear parabolic equations}

In this section, we establish regularity estimates for solutions to the linear parabolic equation
\begin{align}\label{LinNSsource}
\frac{dU}{dt}
-
\nu \Delta U 
+
B[u_{\theta},U]
+
B[U,u_{\theta'}]
=
g\quad
&\textrm{ on } [0,T]\times \Om,
\\[2mm]\nonumber
U(0)=\xi\quad
&\textrm{ on } \Om
,
\end{align}
where $u_{\theta}$ and $ u_{\theta'}$ are the unique solutions of the Navier-Stokes equation (\ref{eq:NS}) associated with source term $f$ and initial conditional $\theta$ and $\theta'$, respectively, and where the equation holds in $\divH^{a-1}$ for general $a \in \Z$. 

\subsubsection{Sobolev bounds on the solution}

\begin{Prop}\label{PropLinPSob}
    Fix $a\in\Z$ and $a^*$ as in (\ref{eq:DefinAStar}). Then for any $\theta,\theta'\in \divH^{a^*}$, $f\in\divH^{a^*-1}$, $g\in L^2([0,T],\divH^{a-1})$, and~$\xi\in\divH^a$, there exists a unique solution $U=U_{\xi,g}:[0,T]\times \Om\to\R^2$ to the parabolic equation~(\ref{LinNSsource}), and we have 
    $$
    U\in C([0,T], \divH^a)\cap L^2([0,T], \divH^{a+1}), \frac{dU}{dt} \in L^2([0,T], \divH^{a-1}).
    $$
    In addition, for any $f\in\divH^{a^*-1}$ and $L>0$ there exists $C=C(\nu,T, \|f\|_{\dH^{a^*-1}}, a, L)>0$ such that for all $\xi\in \divH^a$, $g\in L^2([0,T],\divH^{a-1})$,  $\theta,\theta'\in\divH^{a^*}$ with $\|\theta\|_{\dH^{a^*}}+\|\theta'\|_{\dH^{a^*}}\le L$, and all $t\in [0,T]$
    $$
    \|U(t)\|^2_{\dH^a}
    +    
    \int_0^t  \|U(t)\|^2_{\dH^{a+1}}\, dt 
    \leq  
    C\Big( \|\xi\|^2_{\dH^a} + \int_0^t \|g(s)\|^2_{\dH^{a-1}}\, ds\Big)
    .
    $$
\end{Prop}


\begin{Proof}{Proposition \ref{PropLinPSob}}
As in the proof of Proposition \ref{PropNSBound}, the following proof should be understood in the sense of a Galerkin approximation where the estimates below are established first for $U_n\equiv \sum_{j=1}^n b_j e_j\in\Hnull$, leading to a system of $n$ ODE's 
\begin{align*}
\frac{dU_n}{dt}
-
\nu \Delta U_n
+
P_n[(u_{\theta} \cdot\nabla)U_n]
+
P_n[(U_n\cdot \nabla) u_{\theta'}]
=
P_n[g]\quad
&\textrm{ on } [0,T]\times \Om,
\\[2mm]\nonumber
U_n(0)=P_n[\xi]\quad
&\textrm{ on } \Om
,
\end{align*}
determining the (time-dependent) coefficients $b_j$, where $P_n:L^2\to \Hnull$ is the orthogonal projection onto the first $n$ eigenfunctions $e_1,\ldots, e_n\in\divH^\infty$ of the Stokes operator $A=-\Delta$. Since the estimates we obtain on $U_n$ are uniform in $n$ in the strong $L^2$-sense, one can show that they transfer to the true solution $U$ of (\ref{LinNSsource}) by a limiting procedure. Since $e_j\in \divH^\infty$ for all $j\geq 1$, in particular we have $U_n\in\divH^\infty$ for all $n$ and the following estimates (to be understood as being performed on $U_n$ for each $n$) are valid. In the remainder of this proof, we 
write $U$ for $U_n$.


    Taking the $\divH^a$-inner product of (\ref{LinNSsource}) with $U$ yields
    $$
    \frac12 \frac{d}{dt} \|U\|^2_{\dH^a}
    +
    \nu \| U\|^2_{\dH^{a+1}} 
    +
    \ps{B[u_{\theta},U]}{U}_{\dH^a} 
    +
    \ps{B[U,u_{\theta'}]}{U}_{\dH^a}
    =
    \ps{g}{U}_{\dH^a}
    ,
    $$
    where we used (\ref{IPPdt}) and the integration by parts formula (\ref{IPPSob2}). Notice that, when $a=0$, then $\ps{B[u_\theta,U]}{U}_{\divH^a}=0$ by virtue of Proposition \ref{PropBilinearFormIdentities1}(i). On the one hand, the Cauchy-Schwarz inequality and Young's inequality (\ref{eq:Young}) entail that
    $$
    \lvert\ps{g}{U}_{\dH^a}\rvert
    \leq 
    \|g\|_{\dH^{a-1}}\|U\|_{\dH^{a+1}}
    \leq 
    \nu^{-1} \|g\|^2_{\dH^{a-1}} + \frac{\nu}{4}\| U\|^2_{\dH^{a+1}}
    .
    $$
    On the other hand, Proposition \ref{PropBilinearFormEstimates} entails that 
    $$
    \lvert \ps{B[u_{\theta},U]}{U}_{\dH^a} \rvert
    +
    \lvert \ps{B[U,u_{\theta'}]}{U}_{\dH^a} \rvert
    \le c_1 \big(\|u_{\theta}\|_{\dH^{a^*}}+\|u_{\theta'}\|_{\dH^{a^*}}\big)
    \|U\|_{\dH^a}\| U\|_{\dH^{a+1}}
    ,
    $$
    for some $c_1=c_1(a)>0$. Because $a^*\geq 1$, $f\in \divH^{a^*-1}$, and $\theta,\theta'\in\divH^{a^*}$, then Proposition \ref{PropNSBound} entails that 
    $$
    \sup_{0\leq t\leq T}
    \Big(
    \|u_{\theta}(t)\|_{\dH^{a^*}}
    +
    \|u_{\theta'}(t)\|_{\dH^{a^*}}
    \Big)
    \leq
    c_2\Big( 
    1 + \|\theta\|^{2a^*}_{\dH^{a^*}}
    + \|\theta'\|^{2a^*}_{\dH^{a^*}}
    \Big)
    <
    \infty
    ,
    $$
    for some constant $c_2=c_2(T,\nu, a, \|f\|_{\dH^{a^*-1}})>0$.
    Consequently, we deduce from Young's inequality (\ref{eq:Young}) that
    $$
    \frac12 \frac{d}{dt} \|U\|^2_{\dH^a}
    +
    \nu \| U\|^2_{\dH^{a+1}} 
    \leq 
    \nu^{-1} \|g\|^2_{\dH^{a-1}} + \frac{c_3}{2} \|U\|^2_{\dH^a} + \frac{\nu}{2}\| U\|^2_{\dH^{a+1}}
    ,
    $$
    for some constant $c_3=c_3(T,\nu, a, \|f\|_{\dH^{a^*-1}}, L)=2c_1 c_2 \nu (1+2 L^{2a^*})>0$. Rearranging terms
    \begin{equation}\label{eq:ParabolicEnergyEstimate}
    \frac{d}{dt} \|U\|^2_{\dH^a}
    +
    \nu \| U\|^2_{\dH^{a+1}} 
    \leq 
    2\nu^{-1} \|g\|^2_{\dH^{a-1}} + c_3 \|U\|^2_{\dH^a}
    .
    \end{equation}    
    In particular, we have
    \begin{equation}\label{eq:GronwallU}
    \frac{d}{dt} \|U\|^2_{\dH^a}
    \leq 
    2\nu^{-1} \|g\|^2_{\dH^{a-1}} + c_3 \|U\|^2_{\dH^a}
    ,
    \end{equation}
    which, after integration in time and recalling that $U(0)=\xi$, provides for all $t\in [0,T]$
    $$
    \|U(t)\|^2_{\dH^a}
    \leq 
    \|\xi\|^2_{\dH^a} + 2\nu^{-1} \int_0^t \|g(s)\|^2_{\dH^{a-1}}\, ds + c_3\int_0^t \|U(s)\|^2_{\dH^a}\, ds 
    .
    $$
    It follows from Grönwall's inequality that, uniformly in $t\in [0,T]$, we have
    $$
    \|U(t)\|^2_{\dH^a}
    \leq
    e^{c_3 t} \Big( \|\xi\|^2_{\dH^a} + 2\nu^{-1} \int_0^t \|g(s)\|^2_{\dH^{a-1}}\, ds \Big) 
    .
    $$
    Integrating (\ref{eq:ParabolicEnergyEstimate}) and neglecting the term $\|U(t)\|^2_{\dH^a}$ now yields 
    \begin{align*}
    \nu
    \int_0^t \| U(s)\|^2_{\dH^{a+1}}\, ds
    &\leq 
    \|\xi\|^2_{\dH^a} + 2\nu^{-1} \int_0^t \|g(s)\|^2_{\dH^{a-1}}\, ds + c_3 \int_0^t \|U(s)\|^2_{\dH^a}\, ds
    \\[2mm]
    &\le 
    (1+c_3 e^{c_3t})
    \Big(
    \|\xi\|^2_{\dH^a}
    +
    2\nu^{-1} \int_0^t \|g(s)\|^2_{\dH^{a-1}}\, ds
    \Big)
    .
    \end{align*}
    Uniqueness of the solution $U$ to (\ref{LinNSsource}) follows by noticing that the difference $w$ of any two such solutions satisfies (\ref{LinNSsource}) with $g=0$ and $\xi=0$, so that the estimates established previously (applied to $w$ withg $g=0$ and $\xi=0$) entail that $w=0$.
    
    It remains to show that $U:[0,T]\to\divH^a$ is a continuous map. This is done as in the proof of Proposition \ref{PropNSBound}, and we first estimate 
    $$
    \bar{B}_a
    \equiv 
    \|B[u_{\theta}, U]\|_{\dH^{a-1}} 
    +
    \|B[U, u_{\theta'}]\|_{\dH^{a-1}}
    .
    $$
    When $a\leq -2$, Lemma \ref{PropBilinearFormIdentities2}(i) combined with the multiplier inequality for Sobolev norms (\ref{multipSob}) entails that 
    \begin{eqnarray*}
    \lefteqn{
        \bar B_a 
        \leq 
    \|u_\theta U^T\|_{\dH^a} + \|U u_{\theta'}^T\|_{\dH^a}
    \lesssim 
    \|u_\theta\|_{\dH^{|a|}}\|U\|_{\dH^a} 
    + 
    \|U\|_{\dH^a} \|u_{\theta'}\|_{\dH^{|a|}}
    }
    \\[2mm]
    &&
    =
    (\|u_{\theta}\|_{\dH^{a^*}}+\|u_{\theta'}\|_{\dH^{a^*}}) \|U\|_{\dH^a}
    \lesssim
    (\|u_{\theta}\|_{\dH^{a^*}}+\|u_{\theta'}\|_{\dH^{a^*}}) \|U\|_{\dH^{a+1}}
    ,
    \end{eqnarray*}
    where we used (\ref{IPPSob3}) in the last inequality. When $a=-1$, the Poincaré inequality (\ref{Poincare}) and Lemma \ref{PropBilinearFormIdentities2} entail that 
    \begin{eqnarray*}
    \lefteqn{
        \bar B_a
        \lesssim 
        \bar B_{a+1}
        \leq 
        \|u_\theta U^T\|_{L^2} + \|U u_{\theta'}^T\|_{L^2}
        \lesssim 
        (\|u_\theta\|_{L^\infty}+\|u_{\theta'}\|_{L^\infty}) \|U\|_{L^2}
    }
    \\[2mm]
    &&\lesssim 
    (\|u_\theta\|_{L^\infty}+\|u_{\theta'}\|_{L^\infty}) \|U\|_{L^2}
    \lesssim 
    (\|u_\theta\|_{\dH^2}+\|u_{\theta'}\|_{\dH^2}) \|U\|_{L^2}
    \\[2mm]
    &&
    \hspace{10mm}=
    (\|u_\theta\|_{\dH^{a^*}}+\|u_{\theta'}\|_{\dH^{a^*}}) \|U\|_{\dH^{a+1}}
    ,
    \end{eqnarray*}
    where we used the Sobolev injection $\dH^2\embed L^\infty$.
    When $a=0$, we have by Lemma \ref{PropBilinearFormIdentities2}(i), \lad's inequality (\ref{eq:ProdL2}), and (\ref{IPPSob3}) 
    \begin{align*}
    \bar B_0 
    &\leq 
    \|u_\theta U^T\|_{L^2} + \|U u_{\theta'}^T\|_{L^2} 
    \\
    &\lesssim 
    \|u_\theta\|_{\dH^1} \|U\|_{\dH^1} + \|u_{\theta'}\|_{\dH^1} \|U\|_{\dH^1}
    \\
    &=
    (\|u_\theta\|_{\dH^{a^*}} + \|u_{\theta'}\|_{\dH^{a^*}}) \|U\|_{\dH^{a+1}}
    .
    \end{align*}
    When $a\geq 1$, we have $\max\{|a-1|+1,2\} = \max\{|a|,2\}=a^*$, so that Proposition \ref{PropProdSob} entails 
    $$
    \bar B_a
    \lesssim 
    (\|u_\theta\|_{\dH^{a^*}} + \|u_{\theta'}\|_{\dH^{a^*}}) \|U\|_{\dH^{a}}
    \lesssim 
    (\|u_\theta\|_{\dH^{a^*}} + \|u_{\theta'}\|_{\dH^{a^*}}) \|U\|_{\dH^{a+1}}
    ,
    $$
    where we used the Poincaré inequality (\ref{Poincare}).
    Using the previous estimates on $\|u_\theta\|_{\dH^{a^*}} + \|u_{\theta'}\|_{\dH^{a^*}}$, we deduce that, irrespective of $a\in\Z$, we have
    $\bar B_a 
    \lesssim
    \|U\|_{\dH^{a+1}}
    .$
    We then have
    \begin{align*}
    \Big\|\frac{dU}{dt}\Big\|_{\dH^{a-1}}
    &\leq 
    \nu \|\Delta U\|_{\dH^{a-1}} 
    +
    \|B[u_{\theta}, U]\|_{\dH^{a-1}} 
    +
    \|B[U, u_{\theta'}]\|_{\dH^{a-1}}
    +
    \|g\|_{\dH^{a-1}}
    \\[2mm]
    &\lesssim
    1+
    \|U\|_{\dH^{a+1}}
    +
    \|g\|_{\dH^{a-1}}
    ,
    \end{align*}
    where we used (\ref{IPPSob3}). Using the previous estimate on $\int_0^T \|U\|^2_{\dH^{a+1}}\, dt$, we deduce that 
    $$
    \int_0^T \Big\|\frac{dU}{dt}\Big\|^2_{\dH^{a-1}}
    <
    \infty 
    ,
    $$
    so that $U\in L^2([0,T],\divH^{a+1})$ with $dU/dt\in L^2([0,T], \divH^{a-1})$. The conclusion now follows from Theorem 4 of Section 5.9 in \cite{Evans1998}.
\end{Proof}

\subsubsection{Sobolev stability estimates}


In the next proposition, we establish a stability estimate for the solution $U$ of (\ref{LinNSsource}) of the form 
\begin{equation}\label{eq:Stab}
\|\xi\|_{\dH^{-1}}
\lesssim
\|U\|_{L^2([0,T],\Hnull)}
,\spa 
\forall\ \xi\in \Hnull 
,
\end{equation}
when the inhomogeneous term $g$ in (\ref{LinNSsource}) is zero. 
The proof of (\ref{eq:Stab}) is rather elementary but crucially exploits the fact that $U(t)\to \xi$ as $t\downarrow 0$, established in Proposition \ref{PropLinPSob}.

\begin{Prop}\label{PropLinPStab}
    Let $a\in\Z$, and $a^*$ be as in (\ref{eq:DefinAStar}). Then, for all $f\in \divH^{a^*-1}$ and $L>0$ there exists~$C=C(\nu, T, \|f\|_{\dH^{a^*-1}},a, L)>0$ such that for all $\theta,\theta'\in\divH^{a^*}$ with $\|\theta\|_{\dH^{a^*}}+\|\theta'\|_{\dH^{a^*}}\le L$ and $\xi\in \divH^a$, the solution $U$ to (\ref{LinNSsource}) associated with $g=0$ satisfies 
    $$
    \|\xi\|^2_{\dH^{a}}
    \le 
    C\int_0^T \|U(t)\|^2_{\dH^{a+1}}\, dt
    .
    $$
\end{Prop}



\begin{Proof}{Proposition \ref{PropLinPStab}}
    Proposition \ref{PropLinPSob} with $g=0$ entails that $U\in C([0,T], \divH^a)$. In particular, we have $U(t)\to \xi$ in $\divH^a$ as $t\downarrow 0$. From the proof of Proposition \ref{PropLinPSob} (still with $g=0$) 
    $$
    \frac12 \frac{d}{dt}\|U\|^2_{\dH^{a}}
    +
    \nu \|U\|^2_{\dH^{a+1}}
    +
    \ps{B[u_{\theta}, U]}{U}_{\dH^{a}}
    +
    \ps{B[U,u_{\theta'}]}{U}_{\dH^{a}}
    =
    0
    ,
    $$
    with 
    $$
    \lvert 
    \ps{B[u_{\theta}, U]}{U}_{\dH^{a}}
    +
    \ps{B[U,u_{\theta'}]}{U}_{\dH^{a}}
    \rvert 
    \leq 
    c \|U\|_{\dH^a}\|U\|_{\dH^{a+1}}
    ,
    $$
    for $c=c(T,\nu, a, \|f\|_{\dH^{a^*-1}}, L)>0$. The Poincaré inequality $\|U\|_{\dH^a}\leq \lambda_1^{-1} \|U\|_{\dH^{a+1}}$ from (\ref{Poincare}) yields 
    $$
    0
    \le
    \frac12 \frac{d}{dt}\|U\|^2_{\dH^{a}}
    +
    (\nu+c\lambda_1^{-1})\|U\|^2_{\dH^{a+1}}
    .
    $$
    Integrating the last display in time between $0$ and an arbitrary $t\in [0,T]$, and using that $U(t)\to \xi$ in $\dH^a$ as $t\to 0$, yields 
    $$
    \|\xi\|^2_{\dH^a} 
    \le
    \frac12 \|U(t)\|^2_{\dH^a} 
    + 
    (\nu+c\lambda_1^{-1}) \int_0^t \|U(s)\|^2_{\dH^{a+1}}\, ds 
    ,\spa 
    \forall\ t\in [0,T]
    .
    $$
    Integrating again and using again the Poincaré inequality (\ref{Poincare}) provides 
    \begin{align*}
    \|\xi\|^2_{\dH^a} 
    &\leq 
    \frac{1}{2T} \int_0^T \|U(t)\|^2_{\dH^a}\, dt
    +
    \frac{\nu+c\lambda_1^{-1}}{T} \int_0^T (T-t) \|U(t)\|^2_{\dH^{a+1}}\, dt 
    \\[2mm]
    &\leq 
    \frac{\lambda_1^{-2}}{2T} \int_0^T \|U(t)\|^2_{\dH^{a+1}}\, dt
    +
   (\nu+c\lambda_1^{-1}) \int_0^T \|U(t)\|^2_{\dH^{a+1}}\, dt 
    \\[2mm]
    &=
    \Big(\frac{\lambda_1^{-2}}{2T}+\nu+c\lambda_1^{-1}\Big) \int_0^T \|U(t)\|^2_{\dH^{a+1}}\,dt,
    \end{align*}
    which concludes the proof.
\end{Proof}

\subsubsection{Smoothing of the parabolic flow}

\begin{Prop}\label{PropLinPSmooth}
    Let $a\in\Z$, $b\in\Z$ such that $b\ge a+1$, and $b^*$  as in (\ref{eq:DefinAStar}). Then for every $f\in\divH^{b^*-1}$, $t_{\min}\in (0,T)$, and $L>0$, there exists $C=C(\nu,T,\|f\|_{\dH^{b^*-1}}, L, t_{\min})>0$ such that for all $\theta,\theta'\in \divH^{b^*}$ with $\|\theta\|_{\dH^{b^*}}+\|\theta'\|_{\dH^{b^*}}\le L$ and $\xi\in\divH^a$ the solution $U$ to (\ref{LinNSsource}) associated with $g=0$ satisfies
    $$
    \sup_{t\in [t_{\min}, T]} \|U(t)\|_{\dH^b}
    \leq
    C
    \|\xi\|^2_{\dH^a}
    .
    $$
\end{Prop}

\begin{Proof}{Proposition \ref{PropLinPSmooth}}
    Since $\theta,\theta'\in \divH^{b^*}$ and $f\in \divH^{b^*-1}$ with $b\geq a+1$, we have from (\ref{eq:GronwallU}) in the proof of Proposition \ref{PropLinPSob} with $\dH^{a+1}$-norms instead of $\dH^a$
    $$
    \frac{d}{dt}\|U\|^2_{\dH^{a+1}}
    \leq 
    c\|U\|^2_{\dH^{a+1}}
    .
    $$
    Fix $s\in (0,t)$ and integrate the last inequality between $s$ and $t$ to obtain 
    $$
    \|U(t)\|^2_{\dH^{a+1}}
    \leq
    \|U(s)\|^2_{\dH^{a+1}}
    +
    c \int_s^t \|U(\tau)\|^2_{\dH^{a+1}}\, d\tau 
    .
    $$
    Integrating the last display in $s\in [0,t]$ provides
    $$
    t \| U(t)\|^2_{\dH^{a+1}}
    \leq 
    \int_0^t\| U(s)\|^2_{\dH^{a+1}}\, ds
    +
    c\int_0^t s \| U(s)\|^2_{\dH^{a+1}}\, ds 
    .
    $$
    Since $b\geq a$, Proposition \ref{PropLinPSob} entails that we have, uniformly in $t\in [t_{\min}, T]$
    $$
    \|U(t)\|^2_{\dH^{a+1}}
    \lesssim 
    \|\xi\|^2_{\dH^a}
    .
    $$
    Because the map $t\mapsto U(t+\tau)$ satisfies (\ref{LinNSsource}) translated in time by $\tau$, an induction argument yields for any integer $k$ and $\tau \in [t_{\min}/k, T/k]$
    $$
    \|U(k\tau)\|_{\dH^{a+k}}
    \lesssim 
    \|U((k-1)\tau)\|_{\dH^{a+(k-1)}}
    \lesssim 
    \ldots 
    \lesssim 
    \|U(0)\|_{\dH^a}
    =
    \|\xi\|_{\dH^a}
    .
    $$
    Taking $k=b-a$ and $\tau\in [t_{\min}/k, T/k]$ yields the conclusion.
\end{Proof}

\begin{acks}[Acknowledgments]
The authors thank Gustav R\o{}mer for many discussions.
\end{acks}
\begin{funding}
The authors gratefully acknowledge funding from an ERC Advanced Grant (UKRI G116786) as well as by EPSRC programme grant EP/V026259.
\end{funding}

\bibliographystyle{imsart-nameyear}
\bibliography{mybib}       


\end{document}